\documentclass[10pt]{amsart}

\usepackage{amssymb,amsmath}
\usepackage{amsthm,mathabx}
\usepackage{qsymbols}
\usepackage{graphicx}

\usepackage[utf8]{inputenc}

\usepackage[T1]{fontenc}
\usepackage[pdftex,pdfborder={0 0 0},colorlinks,pagebackref]{hyperref}
\usepackage{hyperref}

\newtheorem{thm}{Theorem}[section]
\newtheorem{lem}{Lemma}[section]
\newtheorem{prop}{Proposition}[section]
\newtheorem{cor}{Corollary}[section]
\newtheorem{defi}{Definition}[section]
\newtheorem{rem}{Remark}[section]

\newcommand{\X}{\mathcal{X}}
\newcommand{\Y}{\mathcal{Y}}
\newcommand{\Z}{\mathcal{Z}}
\newcommand{\R}{\mathbb{R}}
\newcommand{\N}{\mathbb{N}}
\newcommand{\ep}{\varepsilon}
\newcommand{\E}{\mathbb{E}}

\usepackage{color}

\begin{document}

\title{Weak optimal transport with unnormalized kernels}

\author{Philippe Chon\'e, Nathael Gozlan and Francis Kramarz}

\date{\today}

\thanks{P. Choné gratefully acknowledges support from the French National Research Agency (ANR), (Labex ECODEC No. ANR-11-LABEX-0047) and F. Kramarz from the European Research Council (ERC), advanced grant FIRMNET (project ID 741467). N. Gozlan is supported by a grant of the Simone and Cino Del Duca Foundation.}

\address{Philippe Choné: CREST-ENSAE, Institut Polytechnique de Paris, France } \email{Philippe.Chone@ensae.fr}

\address{Nathael Gozlan: Universit\'e Paris Cité, MAP5, UMR 8145, 45 rue des Saints Pères, 75270 Paris Cedex 06}
\email{natael.gozlan@parisdescartes.fr}

\address{Francis Kramarz: CREST-ENSAE, Institut Polytechnique de Paris, France and Department of Economics, Uppsala University, Sweden}
\email{Francis.Kramarz@ensae.fr}


\begin{abstract}
We introduce a new variant of the weak optimal transport problem where mass is distributed from one space to the other through unnormalized kernels. We give sufficient conditions for primal attainment and prove a dual formula for this transport problem. We also obtain dual attainment conditions for some specific cost functions. As a byproduct, we obtain a transport characterization of the stochastic order defined by convex positively $1$-homogenous functions, in the spirit of Strassen theorem for convex domination.
\end{abstract}
\maketitle

\section*{Introduction}

In this paper, we study the mathematical aspects of a new variant of the optimal transport problem, related to the weak optimal transport problem introduced in \cite{GRST17}, recently considered by the first and third authors in \cite{CK21} in an economic context.

In what follows $\X$ and $\Y$ are compact metrizable spaces, $\mathcal{P}(\X)$ (resp. $\mathcal{P}(\Y)$) denotes the set of all Borel probability measures on $\X$ (resp. $\Y$) and $\mu \in \mathcal{P}(\X)$ and $\nu \in \mathcal{P}(\Y)$ are fixed probability measures.

In the usual Monge-Kantorovich transport problem,  given a cost function $\omega : \X \times \Y \to \R$, assumed to be measurable and bounded from below, the optimal transport cost between $\mu$ and $\nu$ is defined as
\begin{equation}\label{eq:MK}
\mathcal{T}_\omega(\mu,\nu) = \inf_{\pi \in \Pi(\mu,\nu)} \iint \omega(x,y)\,\pi(dxdy),
\end{equation}
where $\Pi(\mu,\nu)$ denotes the set of all couplings between $\mu$ and $\nu$, that is the set of all probability measures $\pi$ on $\X \times \Y$ such that the $\X$-marginal of $\pi$ is $\mu$ and the $\Y$-marginal of $\pi$ is $\nu$. The quantity $\pi(dxdy)$ represents the amount of mass taken at $x$ and sent to $y$. The marginal constraints on $\pi$ thus ensure that all the mass of $\mu$ has been recombined to reconstruct $\nu$. We refer to the textbooks \cite{Vil03,Vil09, Gal16,PC19} for a panorama of applications.

To motivate the introduction of weak optimal transport, recall that any coupling $\pi \in \Pi(\mu,\nu)$ can be disintegrated as follows
\[
\pi(dxdy) = \mu(dx)p^x(dy),
\]
where $p=(p^x)_{x\in \X}$ is a probability kernel from $\X$ to $\Y$ (which is $\mu$ almost surely unique). In an informal way, for all $x \in \X$ the probability $p^x \in \mathcal{P}(\Y)$ contains all the information about how the mass taken at $x$ is distributed over $\Y$.
Using this notation, one sees in particular that
\[
 \iint \omega(x,y)\,\pi(dxdy) = \int \left(\int \omega (x,y)\,p^x(dy)\right)\mu(dx),
\]
which highlights that, in the Monge-Kantorovich optimal transport problem, the mass transfers from $\X$ to $\Y$ are penalized only through their mean costs $\int \omega (x,y)\,p^x(dy)$, $x\in \X.$ By contrast, the Weak Optimal Transport (WOT) framework allows to consider more general penalizations on the probability kernel $p$. Given a cost function $c : \X \times \mathcal{P}(\Y) \to \R$ assumed to be measurable and bounded from below, the weak optimal transport cost between $\mu$ and $\nu$ is defined as
\begin{equation}\label{eq:WOT}
\mathcal{T}_c(\mu,\nu) = \inf_{p \in \mathcal{P}(\mu,\nu)} \int c(x,p^x)\,\mu(dx),
\end{equation}
where $\mathcal{P}(\mu,\nu)$ denotes the set of all probability kernels $p=(p^x)_{x\in \X}$ transporting $\mu$ onto $\nu$, i.e. such that
\[
\mu p (dy):= \int p^x(dy)\,\mu(dx) = \nu(dy).
\]
This definition, which finds its origin in the works by Marton \cite{Mar96b,Mar96a} on transport-entropy inequalities and their relations to the concentration of measure phenomenon \cite{GL10,Led01}, stricly extends the setting of the Monge-Kantorovich transport problem (which corresponds to $c(x,p) = \int \omega(x,y)\,p(dy)$, $x\in \X$, $p\in \mathcal{P}(\Y)$). It turns out that the WOT setting includes several interesting variants of the optimal transport problem such as the Schr\"odinger / entropic regularized transport problem \cite{Leo14,Cut13,BCC15} or the martingale transport problem \cite{BHLP13, GHLT14, BJ16}. General tools such as a Kantorovich type duality formula \cite{GRST17, ABC19, BBP19} and a cyclical monotonicity criterium \cite{GJ20, BBP19} have been developed within the WOT framework. We refer to the nice survey paper \cite{BP22} for a general panorama of recent results and applications of WOT.  Among the new WOT problems recently considered, the class of barycentric transport problems attracted a particular attention. These barycentric transport problems correspond to cost functions of the form
\[
c(x,p) = \theta\left(x- \int y\,p(dy)\right),\qquad x\in \X, p \in \mathcal{P}(\Y),
\]
with $\X, \Y \subset \R^d$ and $\theta: \R^d \times \R^d \to \R$ a convex function bounded from below. The introduction of such barycentric optimal transport costs was initially motivated by their applications in concentration of measure \cite{GRST17, GRSST18}. In dimension $1$ and for a general convex function $\theta$, the structure of optimal plans has been settled in \cite{GRSST18, ACJ20, BBP20}. For the quadratic cost function $\theta = \|\,\cdot\,\|_2^2$ on the Euclidean space $(\R^d, \|\,\cdot\,\|_2)$ with $d\geq 1$, the structure of optimal plans has been described in  \cite{GJ20, BBP19}, and yields a new characterization of the couples $(\mu,\nu)$ for which the Brenier transport map \cite{Br87, Br91} is a contraction and to a new proof \cite{FGP20} of the Caffarelli contraction theorem \cite{Caf00}. These barycentric cost functions also recently found applications in machine learning \cite{CTF22}.

The new setting studied in the present paper consists in relaxing the assumption that $p=(p^x)_{x\in \X}$ appearing in \eqref{eq:WOT} is a \emph{probability} kernel. To state a formal definition, we need to introduce additional notions and notations. We will denote by $\mathcal{M}(\Y)$ the set of all finite nonnegative measures on $\Y$. This set will always be equipped with the usual weak topology, and with the cylindric $\sigma$-field (the $\sigma$-field generated by the maps $\mathcal{M}(\Y) \ni m \mapsto m(A)$ with $A$ a Borel subset of $\Y$). A \emph{nonnegative kernel from $\X$ to $\Y$} is a collection $q=(q^x)_{x\in \X}$ of elements of $\mathcal{M}(\Y)$ such that the map $\X \to \R_+ : x\mapsto q^x(A)$ is measurable for all Borel set $A \subset \Y$. Given $\mu \in \mathcal{P}(\X)$ and $\nu \in \mathcal{P}(\Y)$ and a measurable cost function $c:\X\times \mathcal{M}(\Y) \to \R$ such that there exist $r_0,r_1 \in \R$ such that
\begin{equation}\label{eq:cbound}
\tag{LB} c(x,m) \geq r_0 + r_1 m(\Y),\qquad \forall x \in \X,\forall m \in \mathcal{M}(\Y),
\end{equation} we consider
\begin{equation}\label{eq:TPintro}
\mathcal{I}_c(\mu,\nu) = \inf_{q\in \mathcal{Q}(\mu,\nu)} \int c(x,q^x)\,\mu(dx),
\end{equation}
where $\mathcal{Q}(\mu,\nu)$ denotes the set of all nonnegative kernels from $\X$ to $\Y$ such that $\mu q = \nu$, where as above $\mu q(dy) = \int q^x(dy)\,\mu(dx)$.

\begin{rem}\label{rem:unbalanced}\
The same transport problem can be stated for measures $\mu$, $\nu$ with different masses. It is not difficult to see that this unbalanced problem can be reduced to the one above simply by redefining the cost function. Hence, we reduce attention to the balanced case.
\end{rem}

The transport problem \eqref{eq:TPintro} has been introduced in \cite{CK21}  with the following economic motivation. The space $\X$ represents firms' technologies in a given industry and the space $\Y$ represents workers' skills. The probability measures $\mu$ and $\nu$ represent the distributions of firms and of workers respectively. A firm of type $x\in \X$ which recruits a distribution of workers $m \in \mathcal {M}(\mathcal{\Y})$ produces output given by $-c(x,m)$. The problem \eqref{eq:TPintro} consists in maximizing total output in the industry over all possible assignments of workers to firms. More precisely, a nonnegative kernel $q \in \mathcal{Q}(\mu,\nu)$ represents a particular hiring policy, with $q^x(dy)$ giving the distribution of workers hired by firms of type $x \in \X$. The condition $\mu q=\nu$ expresses that all workers are employed. The mass $N(x):=q^x(\mathcal{\Y})$ represents the total number of workers hired by a firm $x \in \X$, i.e., the size of firms with technology~$x$. Importantly, these firms' sizes are an outcome of the optimization process, whereas OT models restrict to \emph{probability} kernels and hence cannot accommodate this question.

The main difficulty in dealing with the transport problem \eqref{eq:TPintro} is that, unlike Problem \eqref{eq:WOT}, assuming that the cost function is jointly lower semicontinuous and convex in its second variable is not enough to ensure existence of a minimizer. To obtain existence of a minimizer, one needs to introduce additional conditions:
\begin{itemize}
\item We will say that the cost function $c:\X \times \mathcal{M}(\Y) \to \R$ satisfies Assumption \eqref{eq:A-intro} if there exists a family of continuous functions $(a_k)_{k\geq0}$ on $\X$ and a family of continuous functions $(b_k)_{k\geq0}$ on $\X\times \Y$ such that
\begin{equation} \label{eq:A-intro}
\tag{A} c(x,m) = \sup_{k\geq0} \left\{\int b_k(x,y)\,m(dy) + a_k(x)\right\},\qquad x\in \X, m \in \mathcal{M}(\Y).
\end{equation}
Note that this condition implies in particular that $c$ is jointly lower semicontinuous, convex with respect to its second variable and satisfies \eqref{eq:cbound}.
\item We will say that $c$ satisfies Assumption \eqref{eq:B-intro} if
\begin{equation} \label{eq:B-intro}
\tag{B} \frac{c(x,\lambda m )}{\lambda} \underset{\lambda \to \infty}{\longrightarrow} +\infty, \qquad\forall x\in \X,\ \forall m \in \mathcal{M}(\Y) \setminus \{0\}.
\end{equation}
\end{itemize}

Let us now present the main results of this paper.

Our first main contribution, is a primal attainment result for the transport problem \eqref{eq:TPintro}. Under Assumptions \eqref{eq:A-intro} and \eqref{eq:B-intro}, we show that for all probability measures $\mu \in \mathcal{P}(\X)$ and $\nu \in \mathcal{P}(\Y)$, there exists a nonnegative kernel $q \in \mathcal{Q}(\mu,\nu)$ such that $\mathcal{I}_c(\mu,\nu) = \int c(x,q^x)\,\mu(dx)$ (see Theorem \ref{thm:strongsol}). In a nutshell, Assumption \eqref{eq:B-intro} prevents mass accumulation on sets of $\mu$ measure $0$. Note that existence of solutions can also hold under other types of conditions on $c$ (see in particular Theorem \ref{thm:duality-eco} dealing with nonpositive cost functions $c$ having a moderate growth).

Our second main result is a Kantorovich type duality formula for the transport problem \eqref{eq:TPintro}. Under assumptions \eqref{eq:A-intro} and \eqref{eq:B-intro}, for all probability measures $\mu \in \mathcal{P}(\X)$ and $\nu \in \mathcal{P}(\Y)$, then
\[
\mathcal{I}_c(\mu,\nu) = \sup_{f\in \mathcal{C}_b(\Y)}\left\{  \int K_{c} f(x)\,\mu(dx) - \int f(y)\,\nu(dy)\right\},
\]
where $\mathcal{C}_b(\Y)$ denotes the set of (bounded) continuous functions on $\Y$ and the operator $K_c$ is defined by
\[
K_{c} f(x) = \inf_{m \in \mathcal{M}(\Y)} \left\{\int f\,dm + c(x,m)\right\},\qquad x\in \X.
\]
Note that, at least formally, if one allows $c$ to take the value $+\infty$ and $c(\,\cdot\,,m)=+\infty$ when $m$ is not of mass $1$, then $\mathcal{I}_c(\mu,\nu) = \mathcal{T}_c(\mu,\nu)$ and one recovers the duality formula for WOT \cite{BBP19}. As we shall see in Theorem \ref{thm:duality2}, the duality formula for $\mathcal{I}_c$ actually holds under a more general condition \eqref{eq:Approx}, implied in particular by Assumption \eqref{eq:B-intro}.

The third main contribution of this paper is a general investigation of transport problems involving cost functions of the following form
\begin{equation}\label{eq:conical-intro}
c(x,m) = F\left(x, \int y \,dm\right),\qquad x\in \X, m\in \mathcal{M}(\Y)
\end{equation}
where $\Y$ is a compact subset of $\R^d$ whose conical hull is denoted by $\Z$ and $F:\X \times \Z \to \R$ is convex with respect to its second variable. Such cost functions will be called \emph{conical} in all the paper because such cost functions are naturally related to \emph{positively $1$-homogenous convex} functions (their epigraphs being \emph{cones}). Recall that a function $\varphi : \R^d \to \R$ is said positively $1$-homogenous (or positively homogenous of degree $1$) if $\varphi(tx)=t\varphi(x)$ for all $t\geq0$ and $x\in \R^d.$ This link between conical cost functions and positively $1$-homogenous convex appears in the duality formula for the transport problem \eqref{eq:TPintro}. More precisely, we will prove that if $c$ is a conical cost function satisfying \eqref{eq:A-intro} and the convex hull of the support of $\nu$ does not contain $0$, then, under some mild integrability condition on $F$, the following reduced duality formula holds
\begin{equation}\label{eq:dual-intro}
\mathcal{I}_c(\mu,\nu) = \sup_{\varphi}\left\{  \int Q_{F} \varphi(x)\,\mu(dx) - \int \varphi(y)\,\nu(dy)\right\},
\end{equation}
where $\varphi$ runs over the set of positively $1$-homogenous convex functions and the operator $Q_F$ is defined by
\[
Q_F\varphi(x) = \inf_{z\in \mathcal{Z}} \{\varphi(z) + F(x,z)\}, \qquad x\in \X.
\]
Moreover, under the same assumptions, we will show the existence of \emph{dual optimizers}; see Theorem \ref{thm:duality-conical} for a precise statement. These conical cost functions are precisely those that were considered in \cite{CK21} and motivated the present paper. In the economic model of \cite{CK21}, a dual optimizer $\varphi$ represents a ``wage schedule’’ ($\varphi(y)$ is the wage paid to workers with skills $y$), while $Q_F\varphi(x)$ represents the opposite of the profit earned by firms with technology $x$ (the profit is the produced output $-F(x,z)$ minus the firm’s wage bill $\varphi(z)$, with $z$ being the sum of the skills of firm~$x$’ employees).

A byproduct of our primal attainment and duality results for conical cost functions is what we believe to be a new variant of the Strassen's theorem \cite{Str65}, that we shall now present. Recall that if $\mu,\nu$ are probability measures on $\R^d$ having finite first moments, one says that $\mu$ is dominated by $\nu$ in the convex order, which is denoted by $\mu \leq_c \nu$, if
\begin{equation}\label{eq:convorder-intro}
\int f\,d\mu \leq \int f\,d\nu
\end{equation}
for all convex function $f : \R^d \to \R.$
Strassen's theorem provides the following useful probabilistic characterization of the convex order: $\mu \leq_c \nu$ if and only if there exists a couple of random variables $(X_0,X_1)$ such that $X_0\sim \mu$, $X_1 \sim \nu$ and $(X_0,X_1)$ is a martingale:
\[
\E[X_1 \mid X_0] = X_0\qquad \text{a.s}.
\]
If $\pi(dxdy) = \mu(dx)p^x(dy)$ denotes the law of $(X_0,X_1)$, the martingale condition is equivalent to the following centering condition on the probability kernel $p$ : for $\mu$ almost all $x$,
\[
\int y\,p^x(dy) = x.
\]
Our generalization of Strassen's theorem deals with a weaker variant of the convex order defined as follows:  if \eqref{eq:convorder-intro} holds for all positively $1$-homogenous convex functions $f$, we will say that $\mu$ is dominated by $\nu$ in the \emph{positively $1$-homogenous convex order} and write $\mu \leq_{phc} \nu$. As we will see in Theorem \ref{thm:convorder}, if $\nu$ is a compactly supported probability measure on $\R^d$ such that the convex hull of the support of $\nu$ does not contain $0$, then $\mu \leq_{phc} \nu$ if and only if there exists a nonnegative kernel $q \in \mathcal{Q}(\mu,\nu)$ such that, for $\mu$ almost all $x$,
\begin{equation}\label{eq:martgen-intro}
\int y\,q^x(dy) = x.
\end{equation}
See Theorem \ref{thm:convorder} for the case where $0$ belongs to the convex hull of the support of $\mu$. Let us briefly explain how this Strassen type result is connected to conical costs. Consider $\mu,\nu$ two compactly supported probability measures on $\R^d$ and denote by  $\Y$ the support of $\nu$. For $p>1$, the conical cost function
\[
c(x,m) = \left\|x - \int y\,dm\right\|^p,\qquad x \in \R^d, m \in \mathcal{M}(\Y),
\]
where $\|\,\cdot\,\|$ is an arbitrary norm on $\R^d$, satisfies Assumption \eqref{eq:B-intro} if and only if $0$ does not belong to the convex hull of $\Y$. According to our primal existence result, we thus have $\mathcal{I}_c(\mu,\nu) = 0$ if and only if there is some $q \in \mathcal{Q}(\mu,\nu)$ satisfying \eqref{eq:martgen-intro}. Using the dual formulation \eqref{eq:dual-intro}, one can then show with some extra work, that $\mu \leq_{phc} \nu$ implies that $\mathcal{I}_c(\mu,\nu) = 0$, thus completing the proof.  In our proof, we actually follow a slightly different route, since we use the cost $c$ above with $p=1$ which will allow us to relax the assumption on the support of $\nu$.

The new version of Strassen's theorem will also enable us to describe optimal transport plans for conical transport costs, in the spirit of \cite{GJ20}.
As we will see in Theorem \ref{thm:primal-structure-bis}, as soon as a conical cost function $c$ of the form \eqref{eq:conical-intro} satisfies Assumption \eqref{eq:A-intro} and the convex hull of the support of $\nu$ does not contain $0$, then
\begin{equation}\label{eq:distance-intro}
\mathcal{I}_c(\mu,\nu)= \inf_{\gamma \leq_{phc} \nu} \mathcal{T}_{F}(\mu,\gamma),
\end{equation}
where $\mathcal{T}_F$ denotes the Monge-Kantorovich optimal transport cost associated to the cost function $F$:
\[
\mathcal{T}_F(\mu,\gamma) = \inf_{\pi \in \Pi(\mu,\gamma)} \iint F(x,z)\,\pi(dxdz),\qquad \forall \mu \in \mathcal{P}(\X),\forall \gamma \in \mathcal{P}( \Z).
\]
Moreover, if $q$ is a kernel minimizer for $\mathcal{I}_c(\mu,\nu)$, then the map $S(x) = \int y\,q^x(dy)$, $x\in \X$, does not depend on the particular choice of the optimizer $q$ and provides an optimal transport for the cost $\mathcal{T}_F$ between $\mu$ and a probability measure $\bar{\nu} \leq_{phc} \nu$ that achieves the infimum in \eqref{eq:distance-intro}. The map $S$ can also be related to dual optimizers (see Theorem \ref{thm:articulation} and Corollary \ref{cor:barS}). In the particular case where $\X$ is a compact subset of $\R^d$ and $F(x,z) = \frac{1}{2}\|x-z\|_2^2$, $x,z\in \R^d$, with $\|\,\cdot\,\|_2$ the standard Euclidean norm, more can be said about the form of the transport map $S$. Namely, we show in Theorem \ref{thm:Brenier} that there exists some closed convex set $C$ such that for $\mu$ almost every $x$, $S(x) = x - p_C(x)$ where $p_C$ is the orthogonal projection onto the set $C$.

Let us point out that during the preparation of this work, we learned about the recent paper \cite{KR21}, devoted to the study of
\[
\inf_{\gamma \leq_\mathcal{A} \nu}  \mathcal{T}_{F}(\mu,\gamma),
\]
where $F:\X \times \R^d \to \R_+$ is some lower semicontinuous function, $\mathcal{A}$ is some cone of continuous functions on $\R^d$ and $\gamma \leq_\mathcal{A} \nu$ means that $\int f\,d\gamma \leq \int f\,d\nu$ for all $f \in \mathcal{A}$. A general duality formula has been obtained in \cite{KR21} for these distance functionals (called backward projection there): under adequate assumptions, 
\[
\inf_{\gamma \leq_\mathcal{A} \nu}  \mathcal{T}_{F}(\mu,\gamma) = \sup_{\varphi \in \mathcal{A}} \int Q_F\varphi\,d\mu - \int \varphi\,d\nu,
\]
with $Q_F$ defined as above (with $\Z = \R^d$).
We refer to \cite[Theorem 4.3]{KR21} for a precise statement. Applying this result to the class $\mathcal{A}$ of all convex positively $1$-homogenous functions together with \eqref{eq:distance-intro}, gives back the duality formula \eqref{eq:dual-intro}. Importantly, the two papers complement each other, since the identity \eqref{eq:distance-intro} crucially requires the variant of Strassen theorem for the convex positively $1$-homogenous order proved here. It would be very interesting to see if other backward projections admit representations in terms of weak transport costs $\mathcal{T}_c$ or $\mathcal{I}_c$ for some special classes of cost functions $c$, but this question will not be considered here.

In the conical case described above, a basic feature of the corresponding transport problem is that it admits, in general, more than one solution. This non-uniqueness of solutions is no longer true for other classes of cost functions, also considered in \cite{CK21}, which we describe now. Suppose that $c : \X \times \mathcal{M}(\Y) \to \R$ is given by
\begin{equation}\label{eq:cGF}
c(x,m) = G\left(\int F(x,y)\,m(dy)\right),\qquad x\in \X, m \in \mathcal{M}(\Y),
\end{equation}
where $F:\X \times \Y \to (0,\infty)$ is some continuous function and $G: [0,\infty) \to \R$ is a convex differentiable function. We establish in Theorem \ref{thm:uniqueness1d}, that when $\X,\Y$ are compact subsets of $\R$, $\mu$ has no atoms, $G$ is monotonic and $F:\R^2 \to \R$ is twice continuously differentiable and satisfies the following condition
\[
\frac{\partial^2 \ln F}{\partial x\partial y} \neq 0
\]
then the transport problem \eqref{eq:TPintro} associated to a cost function of the form \eqref{eq:cGF} admits at most one kernel solution (and exactly one whenever $G'(x) \to +\infty$ as $x \to \infty$ for instance). Moreover, this kernel solution is of the following form
\[
\bar{q}^x(dy) =\bar{N}(x)\delta_{\bar T(x)}(dy),
\]
for $\mu$ almost every $x\in \X$, where $\bar{N}$ is a density with respect to $\mu$ and $\bar T$ is a monotonic function. This uniqueness result is obtained as a consequence of a general result of independent interest establishing a relation between the support of primal solutions and dual optimizers (see Proposition \ref{prop:dualprimal} for details).

The paper is organized as follows. In Section \ref{sec:1}, we introduce a modified formulation of the transport problem \eqref{eq:TPintro} involving couplings $\pi$ with a first marginal absolutely continuous with respect to $\mu$ and second marginal equal to $\nu$. This class of couplings being not closed in general, primal attainment is not always true (when it holds, we call such coupling a strong solution). To compensate this non-attainment issue, we introduce the notion of \emph{weak solution}. These weak solutions are defined as limit points of minimizing sequences, and as such always exist. We conclude Section 1 by giving several explicit examples admitting a) only strong solutions; b) only weak (but not strong) solutions; c) solutions of both types. In Section \ref{sec:2}, we show that under adequate assumptions, weak solutions can be interpreted as couplings minimizing a lower semicontinuous functional (on its domain of definition), denoted $\bar{I}_c^\mu$. One of the main result of this Section is Theorem \ref{thm:strongsol} which shows that, under Assumption \eqref{eq:B-intro}, all weak solutions are strong. The main result of Section \ref{sec:duality}, Theorem \ref{thm:duality2}, provides the dual formulation of the transport problem presented above. Section \ref{sec:monotonicity} deals with cost functions of the form \eqref{eq:cGF}. We prove in Theorem \ref{thm:existencedualopt} that the dual problem admits at least one solution. Then, we establish in Proposition \ref{prop:dualprimal} a general link between supports of primal solutions and this dual optimizer, on which relies the proof of the uniqueness result presented above (Theorem \ref{thm:uniqueness1d}). Section \ref{sec:4} is entirely devoted to the study of the transport problem \eqref{eq:TPintro} for conical cost functions. We prove in particular in Theorem \ref{thm:duality-conical} duality and dual attainment under conditions that are weaker than in Theorem \ref{thm:duality2}. This Section also contains the Strassen type result presented above characterizing the positively $1$-homogenous convex order (Theorem \ref{thm:convorder}) and, as a corollary, the identity \eqref{eq:distance-intro}. We also prove in Theorem \ref{thm:duality-eco} a primal attainment result for a special class of nonpositive conical cost functions. Finally, the paper ends with an Appendix containing the proofs of some technical results stated in Section \ref{sec:2}.

\pagebreak

\tableofcontents

\section{A new transport problem}\label{sec:1}
In this section, we first introduce an alternative, albeit equivalent, formulation of the transport problem \eqref{eq:TPintro} which involves couplings and is thus closer to the usual optimal transport framework. Then we introduce the notion of weak solutions which compensate the fact that the transport problem \eqref{eq:TPintro} does not always admit minimizers. Finally, we study several explicit examples of transport problems \eqref{eq:TPintro} in dimension one.

\subsection{Definitions, equivalent formulation and first properties}
First, we introduce some notations. If $E$ is some Polish metric space, we will denote by $\mathcal{P}(E)$ the set of all Borel probability measures on $E$ and by $\mathcal{M}(E)$ (resp. $\mathcal{M}_s(E)$) the set of all nonnegative finite measures (resp. finite signed measures) on $E$. The space $\mathcal{M}_s(E)$ will be equipped with the topology of weak convergence, that is the coarsest topology that makes the maps $\mathcal{M}_s(E) \to \R : m \mapsto \int f\,dm$ continuous for all $f \in \mathcal{C}_b(E)$, the space of all bounded continuous functions on $E$.

In all what follows, $\X$ and $\Y$ will be two compact metrizable spaces and $c : \X \times \mathcal{M}(\Y) \to \R$ will be a cost function which will always be assumed to be convex with respect to its second variable, jointly lower semicontinuous on $\X\times \mathcal{M}(\Y)$ and to satisfy the lower bound \eqref{eq:cbound}. Given $\mu \in \mathcal{P}(\X)$, for any nonnegative kernel $q$ from $\X$ to $\Y$ such that $\mu q (\Y)=1$ we will set
\[
I_c^\mu[q] = \int c(x,q^x)\,\mu(dx).
\]
With this notation, the transport problem \eqref{eq:TPintro} can be restated as
\begin{equation}\label{eq:TP}
\mathcal{I}_c(\mu,\nu) = \inf_{q \in \mathcal{Q}(\mu,\nu)} I_c^\mu[q],
\end{equation}
where we recall that $\mathcal{Q}(\mu,\nu)$ denotes the set of all nonnegative kernels $q$ from $\X$ to $\Y$ such that $\mu q = \nu$. We will say that $\bar{q}\in \mathcal{Q}(\mu,\nu)$ is a \emph{kernel solution} to the transport problem \eqref{eq:TP} if
\[
\mathcal{I}_c(\mu,\nu) = \int c(x,\bar{q}^x)\,\mu(dx).
\]

First, observe that $\mathcal{I}_c$ is jointly convex.
\begin{prop}\label{prop:convexity}
The functional $\mathcal{P}(\X) \times \mathcal{P}(\Y) \to \R \cup\{+\infty\} : (\mu,\nu) \mapsto \mathcal{I}_c(\mu,\nu)$ is convex.
\end{prop}
\proof
Take $\mu_0,\mu_1 \in \mathcal{P}(\X)$, $\nu_0,\nu_1 \in \mathcal{P}(\Y)$ and for $t \in ]0,1[$ let $\mu_t = (1-t)\mu_0 + t\mu_1$ and $\nu_t = (1-t)\nu_0+t\nu_1$.
It will be convenient to work with a reference probability measure $m$ such that $\mu_t \ll m$ for all $t\in [0,1]$. One can take for instance $m = \mu_{1/2}$, since $\mu_{1/2}(A) = 0$ implies that $\mu_0(A)=\mu_1(A)=0$ which implies that $\mu_t(A)=0$.
 Let $q_0$ and $q_1$ be nonnegative kernels such that $\mu_0q_0=\nu_0$ and $\mu_1q_1=\nu_1$.
 Let $q_t$ be the nonnegative kernel defined by
 \[
 q_t^x(dy) = \frac{(1-t) h_0(x)}{(1-t) h_0(x)+th_1(x)} q_0^x(dy) + \frac{t h_1(x)}{(1-t) h_0(x)+th_1(x)} q_1^x(dy),
 \]
 where $h_0$ and $h_1$ are the densities of $\mu_0,\mu_1$ with respect to $m.$
 We have
\begin{align*}
 \mu_tq_t (dy) &= \int q^x_{t}(dy)\mu_t(dx) = \int \left[ (1-t)h_0(x) q^x_{0}(dy)+ th_1(x)q^x_{1}(y) \right] m(dx) \\
 &= (1-t) \nu_0(dy) + t\nu_1(dy) = \nu_t(dy),
 \end{align*}
 hence $\mu_t q_t=\nu_t$.
 By convexity of $c(x,\,\cdot\,)$, we have
 \begin{align*}
 \int c(x,q_t^x)\,\mu_t(dx) &\leq \int (1-t)h_0(x)c(x,q_0^x)+th_1(x)c(x,q_1^x)\,m(dx) \\&= (1-t)\int c(x,q_0^x)\,\mu_0(dx) + t\int c(x,q_1^x)\,\mu_1(dx).
 \end{align*}
 The result then follows by minimizing over $q_0$ and $q_1$.
\endproof

Let us now derive an alternative, but equivalent, formulation of the transport problem \eqref{eq:TP} closer to the classical viewpoint in optimal transport. We will denote by $\Pi(\eta,\nu)$ the set of all transport plans between two probability measures $\eta \in \mathcal{P}(\X)$ and $\nu \in \mathcal{P}(\Y)$, that is the set of probability measures on $\X \times \Y$ having $\eta$ and $\nu$ as marginals. For any $\mu \in \mathcal{P}(\X)$, we will consider
\[
\Pi(\ll \mu,\nu) = \bigcup_{\eta \ll \mu} \Pi(\eta,\nu) \qquad \text{and}\qquad \Pi(\ll \mu,\,\cdot\,) = \bigcup_{\nu \in \mathcal{P}(\Y)} \Pi(\ll \mu,\nu)
\]
where $\eta \ll \mu$ means that $\eta$ is absolutely continuous with respect to $\mu$. In other words, $\Pi(\ll \mu,\,\cdot\,)$ is the set of all probability measures on $\X\times \Y$ whose first marginal is absolutely continuous with respect to $\mu$.
Observe that if $q \in \mathcal{Q}(\mu,\nu)$ then the function $N$ defined by $N(x) = q^x(\Y)$ is such that
\[
\int N(x)\,\mu(dx) = 1.
\]
Therefore, $N$ is a probability density with respect to $\mu$. Moreover, $\pi(dxdy) = \mu(dx)q^x(dy)$ is a transport plan between $\eta(dx) := N(x)\,\mu(dx)$ and $\nu(dy)$. Conversely, if $\eta \in \mathcal{P}(\X)$ is absolutely continuous with respect to $\mu$ and $\pi \in \Pi(\eta,\nu)$ with $\pi(dxdy) = \eta(dx)p^x(dy)$, then the nonnegative kernel $q$ defined by $q^x(dy) = \frac{d\eta}{d\mu}(x)p^x(dy)$, $x \in \X$, belongs to $ \mathcal{Q}(\mu,\nu)$.

With a slight abuse of notation, let us also denote by $I_c^\mu$ the function defined on $\Pi(\ll \mu,\,\cdot\,)$ by
\[
I_c^\mu[\pi]= \int c\left(x,\frac{d\pi_1}{d\mu}(x)p^x\right)\,\mu(dx),\qquad \pi \in \Pi(\ll \mu,\,\cdot\,),
\]
where $\pi_1$ is the first marginal of $\pi$ and $p$ is the probability kernel such that $\pi(dxdy) = \pi_1(dx)p^x(dy)$.
With this notation, it thus holds
\begin{equation}\label{eq:TPpi}
\mathcal{I}_c(\mu,\nu) = \inf_{\pi \in \Pi(\ll \mu,\nu)} I_c^\mu[\pi].
\end{equation}

\begin{defi}[Strong solutions] Let $\mu \in \mathcal{P}(\X)$ and $\nu \in \mathcal{P}(\Y)$; a probability measure $\bar{\pi} \in \mathcal{P}(\X\times \Y)$ is called a \emph{strong} solution of the transport problem \eqref{eq:TP} if $\bar{\pi} \in \Pi(\ll \mu,\nu)$ and $\mathcal{I}_c(\mu,\nu) = I_c^\mu[\bar{\pi}]$.
\end{defi}
Note that if $\bar{q}$ is a kernel solution to the transport problem \eqref{eq:TP} then the transport plan $\bar{\pi}(dxdy) = \mu(dx)\bar{q}^x(dy) \in \Pi(\ll \mu,\nu)$ is a strong solution of the transport problem \eqref{eq:TP} and, conversely, any strong solution defines a kernel solution.

Since the set $\Pi(\ll \mu,\nu)$ is not closed in general, the infimum in the transport problem \eqref{eq:TP} is not always attained and strong solutions may not always exist. This technical issue motivates the introduction of weak solutions.

\begin{defi}[Weak solutions]
Let $\mu \in \mathcal{P}(\X)$ and $\nu \in \mathcal{P}(\Y)$; a probability measure $\bar{\pi} \in \mathcal{P}(\X\times \Y)$ is called a \emph{weak} solution of the transport problem \eqref{eq:TP} if there exists a sequence of transport plans $\pi_n \in \Pi(\ll\mu,\nu)$ such that $\pi_n \to \bar{\pi}$ for the weak topology and $I_c^\mu[\pi_n] \to \mathcal{I}_c(\mu,\nu)$.
\end{defi}
Of course, a strong solution is also a weak solution. Under adequate conditions on the cost function $c$, weak solutions will be interpreted in Section \ref{sec:weaksol} as solutions of a related minimization problem.

Weak solutions always exist as shows the following elementary result.
\begin{prop}\label{prop:ex1}
For any $\mu \in \mathcal{P}(\X)$ and $\nu \in \mathcal{P}(\Y)$, the transport problem \eqref{eq:TP} admits at least one weak solution.
\end{prop}
\proof
Let $(\pi_n)_{n\in \N}$ be a minimizing sequence in $\Pi(\ll \mu,\nu)$, that is $\lim_{n\to \infty}I_c^\mu[\pi_n] = \mathcal{I}_c(\mu,\nu)$.
Since $\X\times \Y$ is compact, the space $\mathcal{P}(\X\times \Y)$ is also compact. Therefore, the sequence $(\pi_n)_{n\in \N}$ admits at least one converging subsequence, and any limit point $\bar{\pi}$ is a weak solution of the transport problem \eqref{eq:TP}.
\endproof

 At least in the simple case when $\X$ is a finite set, strong solutions always exist.
\begin{thm}\label{thm:strongfinite}
Suppose that $\X$ is a finite then, for any $\mu \in \mathcal{P}(\X)$ and $\nu \in \mathcal{P}(\Y)$, every weak solution of the transport problem \eqref{eq:TP} is a strong solution.
\end{thm}
We will need the following lemma.
\begin{lem}\label{lem:Iscifinite}
If $\X$ is a finite set, then for any $\mu \in \mathcal{P}(\X)$ the set $\Pi(\ll \mu,\,\cdot\,)$ is a closed subset of $\mathcal{P}(\X\times \Y)$ and the functional $\Pi(\ll \mu,\,\cdot\,) \to \R : \pi \mapsto I_c^\mu[\pi]$ is lower semicontinuous.
\end{lem}
\proof
Fix $\mu \in\mathcal{P}(\X)$. The fact that $\Pi(\ll \mu,\,\cdot\,)$ is a closed subset of $\mathcal{P}(\X\times \Y)$ is easy to see and left to the reader.
Let $(\pi_n)_{n\in \N}$ be a sequence of  $\Pi(\ll \mu,\,\cdot\,)$ converging to some $\pi$. Since $\Pi(\ll \mu,\,\cdot\,)$ is closed, it follows that $\pi$ belongs to $\Pi(\ll \mu,\,\cdot\,)$. Denote by $\eta_n$ (resp. $\eta$) the first marginal of $\pi_n$ (resp. $\pi$). For all $x \in \X$ such that $\eta_n(x)>0$, $\pi_n^x(dy)=\frac{\pi_n(x,dy)}{\eta_n(x)}$. If $\eta_n(x)=0$, set $\pi_n^x(dy) = \nu(dy)$ (say).
Then $\eta_n \to \eta$ and $\frac{d\eta_n}{d\mu}(x) \to \frac{d\eta}{d\mu}(x)$ for all $x \in \X$ such that $\mu(x)>0$. Also, it is clear that $\pi_n^x(dy) = \frac{\pi_n(x,dy)}{\eta_n(x)} \to \frac{\pi(x,dy)}{\eta(x)}= \pi^x(dy)$ as $n\to \infty$, for all $x$ such that $\eta(x)>0$.
So, using the lower semicontinuity of $c$, one gets
\begin{align*}
\liminf_{n\to \infty}I_c^\mu[\pi_n] = \liminf_{n\to \infty} \sum_{x\in \X} c\left(x,\frac{d\eta_{n}}{d\mu}(x)\pi_n^x\right) \,\mu(x) &\geq \sum_{x\in \X}  \liminf_{n\to \infty}c\left(x,\frac{d\eta_{n}}{d\mu}(x)\pi_n^x\right) \,\mu(x)\\ &\geq  \sum_{x\in \X}  c\left(x,\frac{d\eta}{d\mu}(x)\pi^x\right) \,\mu(x)
=I_c^\mu[\pi],
\end{align*}
which completes the proof.
\endproof

\proof[Proof of Theorem \ref{thm:strongfinite}]
Let $\pi$ be some weak solution of the transport problem \eqref{eq:TP} and $(\pi_n)_{n\in \N}$ be a minimizing sequence converging to $\pi$. Since $\X$ is finite, it follows from Lemma \ref{lem:Iscifinite} that $\Pi(\ll \mu,\nu)$ is closed, and so $\pi \in \Pi(\ll \mu,\nu)$.  According to Lemma \ref{lem:Iscifinite}, it follows that $\mathcal{I}_c(\mu,\nu) = \liminf_{n\to \infty}I_c^\mu[\pi_n] \geq I_c^\mu[\pi]$, and so $\pi$ is a strong solution.
\endproof

\subsection{Examples}
We study below particular cases of the transport problem \eqref{eq:TP} and we describe their set of solutions. These explicit examples show that all the possibility in terms of uniqueness or non-uniqueness or existence of strong solutions can occur.

\subsubsection{An example without strong solution} Suppose that $\mu$ is the uniform measure on $\X = [0, 1]$ and $\nu$ is an arbitrary probability measure on $\Y=[2,3]$ and define
\[
c(x,m) = \int |x-y|^2\,m(dy), \quad x \in [0,1],\qquad m \in \mathcal{M}(\Y).
\]
Then,
\[
\mathcal{I}_c(\mu,\nu) = \inf_{\mu q = \nu} \iint |y-x|^2\mu(dx)q^x(dy).
\]
Since for all $x\in [0,1]$ and $y\in [2,3]$, $|y-x|^2 \geq |y-1|^2$ 
\[
\mathcal{I}_c(\mu,\nu)\geq \int_2^3 |y-1|^2\,\nu(dy)
\]
holds. This lower bound is not reached. Indeed, suppose by contradiction that there is some $q \in \mathcal{Q}(\mu,\nu)$ such that $ \iint |y-x|^2\mu(dx)q^x(dy) = \int |y-1|^2\,\nu(dy)$. Then, denoting by $\pi(dxdy) = \mu(dx)q^x(dy)$ the associated transport plan, it would hold $\pi(\{1\}\times [2,3]) = 1$ and so $\pi = \delta_1 \otimes \nu$ and $\pi_1 = \delta_1$. Since $\delta_1$ is not absolutely continuous with respect to $\mu$ this is not possible. So this problem does not admit strong solutions.

On the other hand, define for all $n\geq 2$, $\pi_n(dxdy) = \eta_n(dx)\otimes \nu(dy)$, with $\eta_n$ the uniform probability measure on $[1-1/n,1]$.
The associated kernel is given by
\[
q_n^x(dy)= n\mathbf{1}_{[1-1/n,1]}(x) \nu(dy),\qquad x \in [0,1],
\]
and, then,
\begin{align*}
\int c(x,q_n^x)\,\mu(dx) &= \int_0^1 \int_2^3 |y-x|^2q_n^x(dy)\mu(dx) = n\int_2^3\int_{1-1/n}^1|y-x|^2dx\nu(dy)\\
& \leq \int_2^3 |y-(1-1/n)|^2\,\nu(dy) \to  \int_2^3 |y-1|^2\,\nu(dy)
\end{align*}
as $n\to \infty$.
This shows that $\mathcal{I}_c(\mu,\nu) = \int_2^3 |y-1|^2\,\nu(dy)$ and that $\pi=\delta_1 \otimes \nu$ is a weak solution of this transport problem.

\subsubsection{An example with a unique strong solution} In this paragraph, we modify the definition of the cost function of the first example and observe the consequence in terms of existence of strong solutions. Let $\mu$ be the uniform distribution on $\X = [0,1]$ and $\nu = \delta_{2}$ on $\Y= [2,3]$ and consider now the cost
\[
c(x,m) = \left(\int |y-x|\,dm\right)^2,\quad x \in [0,1],\qquad m \in \mathcal{M}(\Y).
\]
If $q \in \mathcal{Q}(\mu,\delta_2)$ then $q^x(\R \setminus\{2\}) = 0$ for almost all $x$.
Therefore, denoting $N(x) = q^x(\{2\})$, it holds
\[
I_c^\mu[q] = \int_0^1 (2-x)^2 N^2(x)\,dx.
\]
By Cauchy-Schwarz,
\[
1 = \int_0^1 N(x)\,dx \leq \left(\int_0^1(2-x)^2N^2(x)\,dx\right)^{1/2}  \left(\int_0^1\frac{1}{(2-x)^2}\,dx\right)^{1/2}.
\]
So, letting $C = \left(\int_0^1\frac{1}{(2-x)^2}\,dx\right)^{-1}$, we get that
\[
I_c^\mu[q]  \geq C
\]
and there is equality if and only if $N(x) = \frac{C}{(2-x)^2}$, $x \in [0,1]$. So $q^x(dy) = N(x)\delta_2(dy)$ is the unique nonnegative kernel achieving the minimum in $\mathcal{I}_c(\mu,\delta_2)$. Equivalently $\pi(dxdy) = N(x)\,dx \otimes \delta_2$ is the unique strong solution of the transport problem. It will follow from Theorem \ref{thm:strongsol} below that all weak solutions are actually strong in this example, so $\pi$ is also the unique weak solution.

\subsubsection{An example exhibiting both strong and weak solutions}
Let $\mu$ be the uniform distribution on $[0,1]$ and $\nu(dy) = 2y^2 \mathbf{1}_{[0,1]}(dy) + \frac{1}{3} \delta_{0}(dy)$.
Consider the cost function
\[
c(x,m) = \left|x - \int y\,dm\right|^p,\qquad x\in [0,1], m\in  \mathcal{M}([0,1]),
\]
with $0<p$. We refer to Section \ref{sec:strassen}  for more insights about this type of costs and the construction of the weak solution below.

Let us first show that the transport problem \eqref{eq:TP} between $\mu$ and $\nu$ admits a strong solution.
Consider the nonnegative kernel $\bar{q}^x(dy) = 2x \nu(dy)$. Since $\int x\,\mu(dx) = \int y\,\nu(dy)=\frac{1}{2}$, it is clear that $\mu q = \nu$ and that $\int y\,\bar{q}^x(dy) = x$. Therefore $\int c(x,\bar{q}^x)\,\mu(dx) = 0$, which shows that $\mathcal{I}_c(\mu,\nu)=0$ and $\bar{q}$ is a strong solution.

Now let us construct a weak (but not strong) solution. Define $\eta(dx) = \sqrt{x}\mu(dx) + \frac{1}{3} \delta_0(dx)$ and let $\pi(dxdy) = \eta(dx)\delta_{\sqrt{x}}(dy)$. We claim that $\pi$ is a weak solution of the transport problem \eqref{eq:TP} between $\mu$ and $\nu.$ First, it is easy to check that the second marginal of $\pi$ is $\nu$, in other words $\nu$ is the push-forward of $\eta$ under the map $x\mapsto \sqrt{x}$.
Let us now construct a minimizing sequence converging to $\pi$.
Define $\pi_\ep$, for $0<\ep<1/2$ as the law of $((1-\ep)X+\ep U, \sqrt{X})$, where $(X,U)$ is a couple of independent random variables such that $X \sim \eta$ and $U\sim \mu$. A simple calculation shows that $\pi_\ep(dxdy) = \mu(dx)q^x_\ep(dy)$, where $q$ is the nonnegative kernel defined by
\[
q^x_\ep(dy):=\frac{1}{3\ep} \mathbf{1}_{x\leq \ep} \delta_0(dy) + \frac{2y^2}{\ep} \mathbf{1}_{\left[\sqrt{\max\left(\frac{x-\ep}{1-\ep},0\right)} ; \sqrt{\min\left(\frac{x}{1-\ep},1\right)} \right]}(y)\,dy.
\]
Therefore, for all $x \in [0,1]$,
\begin{align*}
b_\ep(x):=\int y \,q^x_\ep(dy) &= \frac{1}{2\ep} \left[\min\left(\frac{x}{1-\ep},1\right)^2 -\max\left(\frac{x-\ep}{1-\ep},0\right)^2\right]\\
& = \frac{1}{(1-\ep)^2}\left\{\begin{array}{ll} \frac{1}{2}\frac{x^2}{\ep}  & \text{ if } 0\leq x\leq \ep   \\ (x-\frac{\ep}{2})  & \text{ if } \ep \leq x\leq 1-\ep   \\ \frac{1}{2}\frac{(1-x)}{\ep}(1+x-2\ep)  &    \text{ if } 1-\ep \leq x\leq 1\end{array}\right..
\end{align*}
Thus, one sees that for all $x\in [0,1]$, $b_\ep(x) \to x$ as $\ep \to 0$ and that $\sup_{0<\ep<1/2} \sup_{x\in [0,1]} b_\ep(x) <+\infty$.
So, applying the dominated convergence theorem yields
\[
\int c(x,q_\ep^x)\,\mu(dx) \to 0
\]
as $\ep \to 0$. Since $\pi^\ep \to \pi$ in the weak sense, this shows that $\pi$ is a weak solution (and obviously not strong).
\subsubsection{A particular case of a one dimensional nonpositive conical cost function}\label{sec:CK}
Consider the following cost function $c : [\alpha,\beta]\times \mathcal{M}([\gamma,\delta]) \to \R_-$ where $\alpha,\beta,\gamma,\delta \geq0$ .
\begin{equation}\label{eq:CKdim1}
c(x,m) = -x \left(\int y\,dm\right)^{\eta},\qquad x\in [\alpha,\beta], m\in \mathcal{M}([\gamma,\delta]),
\end{equation}
where $0<\eta<1$.
This cost function is a particular case of the cost functions considered in \cite{CK21} (in arbitrary dimensions).
Note that $c$ satisfies Assumption \eqref{eq:cbound}.
Indeed, by concavity of the function $y \mapsto y^\eta$ on $\R_+$, it holds
\[
y^\eta \leq 1 + \eta (y-1),\qquad \forall y\geq0.
\]
Therefore,
\[
c(x,m) \geq - x\left(1-\eta+\eta \int y\,dm\right)\geq -\beta\left(1-\eta+\eta \delta m([\gamma,\delta])\right) := r_0 +r_1 m([\gamma,\delta]).
\]

The following result provides information on strong solutions of the transport problem associated to the cost function $c$ defined above.
\begin{prop}
Let $\mu \in \mathcal{P}([\alpha,\beta])$ and $\nu \in \mathcal{P}([\gamma,\delta])$; the transport problem \eqref{eq:TP} between $\mu$ and $\nu$ with respect to the cost function $c$ defined by \eqref{eq:CKdim1} admits strong solutions.
For instance, denoting by $\bar{\mu}(dx) = \frac{1}{Z}x^{\frac{1}{1-\eta}}\,\mu(dx)$, where $Z$ is a normalizing constant, then the coupling $\bar{\pi}$ given by
\[
\bar{\pi}= \bar{\mu}\otimes \nu
\]
is a strong solution.
More generally,  $q \in \mathcal{Q}(\mu,\nu)$ is a nonnegative kernel solution if and only if there exists some constant $C>0$ such that
\begin{equation}\label{eq:cond-ex-dim1}
\int y\,q^x(dy) = C x^{\frac{1}{1-\eta}}
\end{equation}
for $\mu$ almost all $x \in [\alpha,\beta]$.
In particular, if $\mu$ has a positive density $f$ on $[\alpha,\beta]$ and $\nu$ a positive density $g$ on $[\gamma,\delta]$, and $T:[\alpha,\beta] \to [\gamma,\delta]$ is a continuously differentiable bijection such that for some constant $C>0$ it holds
\begin{equation}\label{eq:cond2a-ex-dim1}
N(x)T(x) = Cx^{\frac{1}{1-\eta}},
\end{equation}
for Lebesgue almost all $x\in [\alpha,\beta]$, where $N$ is the density (with respect to $\mu$) defined by
\begin{equation}\label{eq:cond2b-ex-dim1}
N(x)=\frac{g(T(x))|T'(x)|}{f(x)},\qquad  \forall x \in [\alpha,\beta]
\end{equation}
then $q^x(dy) = N(x)\delta_{T(x)}$, $x \in [\alpha,\beta]$, is a strong solution.
\end{prop}
We will see in Theorem \ref{thm:duality-eco} below that for such cost function, all weak solutions are actually strong.
\proof
Let $\pi \in \Pi(\ll \mu,\nu)$, then $\frac{d\pi_1}{d\mu}(x) = \frac{d\pi_1}{d\bar{\mu}}(x)  \frac{1}{Z}x^{\frac{1}{1-\eta}}$ and so
\begin{multline*}
-I_c^\mu[\pi] =  \int_{\alpha}^\beta x \left(\frac{d\pi_1}{d\mu}\right)^{\eta} \left(\int_\gamma^\delta y\,\pi^x(dy)\right)^{\eta}\,\mu(dx) =   \frac{1}{Z^\eta} \int_{\alpha}^\beta x^{\frac{1}{1-\eta}} \left(\frac{d\pi_1}{d\bar{\mu}}\right)^{\eta} \left(\int_\gamma^\delta y\,\pi^x(dy)\right)^{\eta}\,\mu(dx) \\
 =  Z^{1-\eta} \int_{\alpha}^\beta  \left(\frac{d\pi_1}{d\bar{\mu}}\right)^{\eta} \left(\int_\gamma^\delta y\,\pi^x(dy)\right)^{\eta}\,\bar{\mu}(dx) \leq  Z^{1-\eta} \left(\int_{\alpha}^\beta  \frac{d\pi_1}{d\bar{\mu}}(x) \int_\gamma^\delta y\,\pi^x(dy)\,\bar{\mu}(dx)\right)^{\eta}
 =   Z^{1-\eta} \left( \int y\,\nu(dy)\,\right)^{\eta},
\end{multline*}
where the inequality follows from the concavity of the function $u\mapsto u^{\eta}$.
Note that if $\pi_1=\bar{\mu}$ and $\pi^x(dy)=\nu(dy)$ for all $x$, there is equality. In other words, $\bar{\pi} = \bar{\mu}\otimes \nu$ is a strong solution of the transport problem between $\mu$ and $\nu$.
Moreover, according to the equality case in Jensen's inequality and the strict concavity of $u\mapsto u^{\eta}$, we see that there is equality above if and only if the function $x\mapsto  \frac{d\pi_1}{d\bar{\mu}}(x) \int_\gamma^\delta y\,\pi^x(dy)$ is constant $\bar{\mu}$ almost surely. Writing $\pi(dxdy) = \mu(dx)q^x(dy)$, we see that this condition is equivalent to the existence of $C>0$ such that \eqref{eq:cond-ex-dim1} holds $\mu$ almost everywhere.
Now, let us assume that $\mu$ has a positive density $f$ on $[\alpha,\beta]$ and $\nu$ a positive density $g$ on $[\gamma,\delta]$, and let us look for solutions of the form $q^x(dy) = N(x)\delta_{T(x)}$, where $x\mapsto T(x)$ is a continuously differentiable bijection from $[\alpha,\beta]$ to $[\gamma,\delta]$. First of all,
if $N$ satisfies \eqref{eq:cond2b-ex-dim1}, then for any bounded measurable function $h$ on $[\gamma,\delta]$, it holds
\[
\int_\alpha^\beta h(T(x))N(x)f(x)\,dx = \int_\alpha^\beta h(T(x))g(T(x))|T'(x)|\,dx  =\int_\gamma^\delta h(y) g(y)\,dy,
\]
by the change of variable formula, which shows that $\mu q = \nu$.
Now, according to \eqref{eq:cond2a-ex-dim1}, it holds
\[
\int y\,q^x(dy) = N(x)T(x) = Cx^{\frac{1}{1-\eta}},
\]
for $\mu$ almost every $x \in [\alpha,\beta]$, which shows that $q$ satisfies \eqref{eq:cond-ex-dim1} and completes the proof.
\endproof

In the following result we consider the particular case where $\mu=\nu$ is the uniform measure on $[0,1]$.
\begin{cor}
If $\mu$ and $\nu$ are both the uniform distribution on $[0,1]$, the three following kernels are strong solutions of the problem:
\begin{itemize}
\item \emph{Random sorting:} $q_0^x(dy)= N_0(x) \mu(dy)$ with $N_0(x)=C x^{a_0}$, for all $x \in [0,1]$, $a_0=1/(1-\eta)$ and $C=(2-\eta)/(1-\eta)$ is such that $\int_{0}^{1} N_0(x)\,\mu(dx)=1$;
\item \emph{Positive Assortative Matching:} $q_1^x(dy) = N_1(x) \, \delta_{T_1(x)}$, where for all $x\in [0,1]$,
\[
T_1(x)= x^{a_1}, \ \ N_1(x)= T_1'(x), \ \  a_1 = \frac{2-\eta}{2(1-\eta)}=\frac C2;
\]
\item \emph{Negative Assortative Matching:} $q_2^x(dy) = N_2(x) \, \delta_{T_2(x)}$, where for all $x\in [0,1]$,
\[
T_2(x) = \sqrt{1-x^{a_2}}, \ \ N_2(x) =-T'_2(x), \ \  a_2=\frac{2-\eta}{1-\eta}=C.
\]
\end{itemize}
\end{cor}
\proof
The verification that $q_1$ and $q_2$ are strong solutions is left to the reader.
\endproof

\section{Weak solutions as minimizers of an extended functional}\label{sec:2}
As explained above, the difficulty in dealing with the minimization problem \eqref{eq:TPpi} is that the set $\Pi(\ll \mu,\,\cdot\,)$ is not closed in general, and so the optimal value of the problem can be reached at the boundary. In this Section, we first identify the closure of $\Pi(\ll \mu,\nu)$, using a simple approximation technique from \cite{LV09}. Then, we introduce an explicit functional $\bar{I}_c^\mu$ which is a lower semicontinuous extension of $I_c^\mu$, and we introduce a condition (see \eqref{eq:Approx} below) under which $\bar{I}_c^\mu$ coincides with the lower semicontinuous envelope of $I_c^\mu$. With this condition in force, we can interpret weak solutions as minimizers of $\bar{I}_c^\mu$ on the closure of $\Pi(\ll \mu,\nu)$. Finally, when Assumption \eqref{eq:B-intro} is satisfied, we will see that every weak solution is strong.

\subsection{Closure of $\Pi(\ll \mu,\nu)$}
Let us introduce a general mollifying approximation technique from \cite[Theorem C.5]{LV09}, which will be very useful in the next paragraphs.

\begin{lem}[Lott-Villani \cite{LV09}]\label{lem:LV}
Let $(S,d)$ be an arbitrary compact metric space and $\mu$ be a Borel probability measure on $S$.
There exists a family of kernels $(K_n)_{n\geq0}$ such that
\begin{enumerate}
\item[$(i)$] For all $n\geq 0$, $K_n:S\times S \to \R_+$ is a continuous and symmetric function such that, for all $x \in \mathrm{Supp}(\mu)$, $\int K_n(x,y)\mu(dy)=1$.
\item[$(ii)$] For all continuous function $f : \mathrm{Supp} (\mu) \to \R$, the functions $K_nf$, $n\geq0$, defined by
\begin{equation}\label{eq:K_nf}
K_nf(y) := \int K_n(x,y) f(x)\,\mu(dx),\qquad y\in S,
\end{equation}
are continuous on $\mathrm{Supp} (\mu)$ and such that $K_nf \to f$ uniformly on $\mathrm{Supp} (\mu)$ as $n \to \infty$.
\item[$(iii)$] For all probability measure $\eta \in \mathcal{P}(S)$ such that $\eta\left( \mathrm{Supp}(\mu) \right)=1$, the probability measures $K_n\eta$, $n\geq0$, defined by
\[
K_n\eta (dy) := \int K_n(x,y)\,\eta(dx)\mu(dy)
\]
is such that $K_n\eta \to \eta$ as $n \to \infty$ for the weak convergence.
\end{enumerate}
\end{lem}

For a fixed $\mu \in \mathcal{P}(\X)$, in what follows we denote by $\Pi( \mathrm{Supp}(\mu),\nu)$ the set of probability measures $\pi$ on $\X\times \Y$ such that $\pi_1\left(\mathrm{Supp}(\mu)\right)=1$ and $\pi_2 = \nu$, where $\pi_1$ and $\pi_2$ denote the marginals of  $\pi$ on $\X$ and $\Y$, respectively.
\begin{lem}\label{lem:descr}
For any $\mu \in \mathcal{P}(\X)$ and $\nu \in \mathcal{P}(\Y)$, it holds
\[
\mathrm{cl}\,\Pi(\ll \mu,\nu) = \Pi( \mathrm{Supp}(\mu),\nu),
\]
where $\mathrm{cl}\,\Pi(\ll \mu,\nu)$ denotes the closure of $\Pi(\ll \mu,\nu)$ for the weak topology. More precisely, for any $\pi \in \Pi( \mathrm{Supp}(\mu),\nu)$ with $\pi(dxdy) = \eta(dx)\pi^x(dy)$, the sequence $(\pi_n)_{n\geq0}$ defined for all  $n\geq 0$ by
\[
\pi_n(dxdy) =  \int K_n(x,z)\pi^z(dy)\eta(dz)\mu(dx),
\]
where $(K_n)_{n\geq0}$ is the sequence of kernels given by Lemma \ref{lem:LV} (applied to $S=\X$ and $\mu$) is such that $\pi_n \in \Pi(\ll \mu,\nu)$ for all $n\geq0$ and $\pi_n\to \pi$ for the weak topology as $n\to \infty.$
\end{lem}

\proof[Proof of Lemma \ref{lem:descr}]
The inclusion $\subset$ is clear. Let us show the other inclusion. Let $\pi \in \Pi( \mathrm{Supp}(\mu),\nu)$ and set $\eta = \pi_1$.
We claim that the first marginal of $\pi_n$ is $K_n\eta$ and the second marginal is $\nu$.
Indeed, if $f:\X\to \R$ is a continuous function, then
\begin{align*}
\iint f(x)\,\pi_n(dxdy) &= \iiint  f(x)K_n(z,x)\pi^z(dy)\mu(dx)\eta(dz)\\
&= \iint   f(x)K_n(z,x)\mu(dx)\eta(dz)\\
& =  \int f(x)(K_n\eta)(dx)
\end{align*}
and if $g:\Y \to \R$ is a continuous function, then
\begin{align*}
\iint g(y)\,\pi_n(dxdy) &= \iiint g(y)K_n(z,x)\pi^z(dy)\mu(dx)\eta(dz)\\
&= \iint g(y)\pi^z(dy)\eta(dz)\\
&=  \int g(y)\nu(dy).
\end{align*}
If $f : \X\times \Y\to \R$ is a continuous function, then denoting by $f_y$ the function $x\mapsto f(x,y)$, it follows from Item $(ii)$ of Lemma \ref{lem:LV} that
\[
\iint f(x,y)\,\pi_n(dxdy) = \iiint  f(x,y)K_n(z,x)\mu(dx)\pi^z(dy)\eta(dz) = \iint K_nf_y(z)\,\pi^z(dy)\eta(dz) \to \iint f(y,z)\pi(dydz),
\]
as $n\to \infty$. In other words $\pi_n \to \pi$ in the weak topology.
Also, since $K_n\eta \ll \mu$, $\pi_n$ belongs to $\Pi(\ll \mu,\nu)$ which completes the proof.
\endproof

\subsubsection{Lower semicontinuous extensions of $I_c^\mu$}
For any fixed $\mu \in \mathcal{P}(\X)$, consider the functional
\[
\bar{I}_c^\mu : \mathcal{P}(\X\times \Y) \to \R\cup\{+\infty\}
\]
defined by
\[
\bar{I}_c^\mu[\pi] = \int c\left(x, \frac{d\pi_1^{ac}}{d\mu}(x)\pi^x(dy) \right)\,\mu(dx) +  \int c'_\infty\left(x, \pi^x \right)\,\pi_1^{s}(dx), \qquad \forall \pi \in \mathcal{P}(\X\times \Y),
\]
where $\pi_1 = \pi_1^{ac}+\pi_1^s$ is the decomposition of $\pi_1$ into an absolutely continuous part and a singular part with respect to $\mu$ and
\[
c'_\infty (x,m) = \lim_{\lambda \to \infty} \frac{c(x,\lambda m)}{\lambda},\qquad x\in \X, m\in \mathcal{M}(\Y)
\]
is the recession function of $c(x,\,\cdot\,)$. Note that this limit is always well-defined since, by convexity of $c(x,\,\cdot\,)$, the function $\lambda \mapsto \frac{c(x,\lambda m) - c(x,0)}{\lambda}$ is non-decreasing on $(0,\infty).$

The following proposition shows in particular that, for a fixed $\mu$ and under Assumption \eqref{eq:A-intro}, the functional $\bar{I}_c^\mu$ is a lower semicontinuous extension of $I_c^\mu$.

\begin{prop}\label{prop:barI_csci}
Under Assumption \eqref{eq:A-intro}, the function $\mathcal{P}(\X)\times \mathcal{P}(\X\times \Y) : (\mu,\pi) \mapsto \bar{I}_c^\mu[\pi]$ is lower semicontinuous and such that $\bar{I}_c^\mu = I_c^\mu$ on $\Pi(\ll\mu,\,\cdot\,)$. More generally, the function  $\R_+\times \mathcal{P}(\X)\times \mathcal{P}(\X\times \Y) : (\lambda,\mu,\pi) \mapsto \bar{I}_{c_\lambda}^\mu[\pi]$ is lower semicontinuous, with $c_\lambda :=c(x,\lambda m)$, $\lambda\in \R_+$, $x \in \X$, $m \in \mathcal{M}(\Y)$.
\end{prop}
The proof of Proposition \ref{prop:barI_csci} (which is adapted from  \cite{AFP00}) is presented in Section \ref{app:1} of Appendix.

For a fixed $\mu \in \mathcal{P}(\X)$, let us now introduce the lower semicontinuous envelope of $I_c^\mu$, denoted $\tilde{I}_c^\mu$ and defined as follows: for all $\pi \in  \mathcal{P}(\X\times \Y)$
\[
\tilde{I}_c^\mu[\pi] = \sup_{V \in \mathcal{V}(\pi)} \inf_{\gamma \in V \cap \Pi(\ll \mu,\,\cdot\,)} I_c^\mu[\gamma],
\]
where $\mathcal{V}(\pi)$ denotes the class of all open neighborhoods of $\pi$. By convention $\inf \emptyset = +\infty$, so in particular, $\tilde{I}_c^\mu = +\infty$ outside $ \mathrm{cl}\, \Pi(\ll\mu,\,\cdot\,)= \Pi(\mathrm{Supp}(\mu),\,\cdot\,)$.

At this level of generality, it is not clear whether $\bar{I}_c^\mu$ and $\tilde{I}_c^\mu$ always coincide on $\Pi(\mathrm{Supp}(\mu),\,\cdot\,)$. A necessary and sufficient condition for this to happen is given in the next proposition:

\begin{prop}\label{prop:compbartilde}
For a fixed $\mu \in \mathcal{P}(\X)$ and under Assumption \eqref{eq:A-intro}, it holds
\begin{itemize}
\item[$(i)$] for all $\pi \in \mathcal{P}(\X\times \Y)$, $\bar{I}_c^\mu[\pi] \leq \tilde{I}_c^\mu[\pi]$,
\item[$(ii)$] for all $\pi \in \Pi(\ll \mu,\,\cdot\,)$, $I_c^\mu[\pi] = \bar{I}_c^\mu[\pi] = \tilde{I}_c^\mu[\pi]$,
\item[$(iii)$]  the functionals $\bar{I}_c^\mu$ and $\tilde{I}_c^\mu$ coincide on $\Pi(\mathrm{Supp}(\mu),\,\cdot\,)$ if and only if for all $\pi \in \Pi(\mathrm{Supp}(\mu),\,\cdot\,)$ there exists a sequence $\pi_n \in \Pi(\ll \mu,\,\cdot\,)$  converging to $\pi$ for the weak topology and such that $I_c^\mu[\pi_n] \to \bar{I}_c^\mu[\pi]$.
\end{itemize}
\end{prop}

\proof
Since $\bar{I}_c^\mu$ is lower semicontinuous, for all $\pi \in \mathcal{P}(\X\times \Y)$ it holds
\[
\bar{I}_c^\mu[\pi] = \sup_{V \in \mathcal{V}(\pi)} \inf_{\gamma \in V} \bar{I}_c^\mu[\gamma] \leq \sup_{V \in \mathcal{V}(\pi)} \inf_{\gamma \in V \cap \Pi(\ll\mu,\,\cdot\,)} \bar{I}_c^\mu[\gamma] =  \sup_{V \in \mathcal{V}(\pi)} \inf_{\gamma \in V \cap \Pi(\ll\mu,\,\cdot\,)} I_c^\mu[\gamma] =\tilde{I}_c^\mu[\pi],
\]
and so $\bar{I}_c^\mu \leq \tilde{I}_c^\mu$, which proves $(i)$.  On the other hand, if $\pi \in \Pi(\ll\mu,\,\cdot\,)$, then
\[
\inf_{\gamma \in V \cap \Pi(\ll \mu,\,\cdot\,)} I_c^\mu[\gamma] \leq I_c^\mu[\pi] = \bar{I}_c^\mu[\pi]
\]
and so, optimizing over $V\in \mathcal{V}(\pi)$, $\tilde{I}_c^\mu[\pi] \leq \bar{I}_c^\mu[\pi]$ which proves $(ii)$.
Let us prove $(iii).$ Suppose that $\pi \in \Pi(\mathrm{Supp}(\mu),\,\cdot\,)$ is such that there exists a sequence $\pi_n \in \Pi(\ll \mu,\,\cdot\,)$ for which $I_c^\mu[\pi_n] \to \bar{I}_c^\mu[\pi]$. Then, since $\tilde{I}_c^\mu$ is lower semicontinuous, the following holds
\[
\tilde{I}_c^\mu[\pi] \leq \liminf_{n\to \infty}\tilde{I}_c^\mu[\pi_n] = \liminf_{n\to \infty}I_c^\mu[\pi_n] = \bar{I}_c^\mu[\pi].
\]
Since the inequality $\bar{I}_c^\mu[\pi] \leq \tilde{I}_c^\mu[\pi]$ is always true, it is in fact an equality.
Conversely, suppose that $\tilde{I}_c^\mu=\bar{I}_c^\mu$ on $ \Pi(\mathrm{Supp}(\mu,\,\cdot\,)$.  If $\pi \in \Pi(\mathrm{Supp}(\mu),\,\cdot\,)$, then according to Lemma \ref{lem:descr}, $\pi \in \mathrm{cl}\, \Pi(\ll\mu,\,\cdot\,)$. Therefore, for any open neighborhood $V$ of $\pi$, the set $V\cap \Pi(\ll\mu,\,\cdot\,)$ is non-empty. Now, it easily follows from the definition of  $\tilde{I}_c^\mu$, that there exists some sequence $\pi_n \in \Pi(\ll\mu,\,\cdot\,)$ such that $\pi_n \to \pi$ for the weak topology with $I_c^\mu[\pi_n] \to \tilde{I}_c^\mu[\pi]$ and so $I_c^\mu[\pi_n] \to \bar{I}_c^\mu[\pi]$.
\endproof

For a fixed $\mu \in \mathcal{P}(\X)$, let us introduce the following variants of problem \eqref{eq:TP}: for $\nu \in \mathcal{P}(\Y)$,
\begin{equation}\label{eq:TPbar}
\widebar{\mathcal{I}}_c(\mu,\nu) = \inf_{\pi \in \Pi(\mathrm{Supp}(\mu),\nu)} \bar{I}_c^\mu[\pi]
\end{equation}
and
\begin{equation}\label{eq:TPtilde}
\widetilde{\mathcal{I}}_c(\mu,\nu) = \inf_{\pi \in \Pi(\mathrm{Supp}(\mu),\nu)} \tilde{I}_c^\mu[\pi].
\end{equation}
Unlike transport problem \eqref{eq:TP}, the transport problems \eqref{eq:TPbar} and \eqref{eq:TPtilde} always admit solutions.
\begin{lem}\label{lem:soltildebar}
Under Assumption \eqref{eq:A-intro}, for any $\mu \in \mathcal{P}(\X),\nu \in \mathcal{P}(\Y)$ such that $\widebar{\mathcal{I}}_c(\mu,\nu) <\infty$, there exists $\pi \in \Pi(\mathrm{Supp}\,(\mu),\nu)$ such that $\bar{I}_c^\mu[\pi]=\widebar{\mathcal{I}}_c(\mu,\nu).$ The same is true for the transport problem \eqref{eq:TPtilde}.
\end{lem}
\proof
The functional $\bar{I}_c^\mu$ is lower semicontinuous on the compact set $\Pi(\mathrm{Supp}(\mu),\nu)$ so it attains its lower bound.
\endproof

Finally, the following result will be very useful in Section \ref{sec:duality} where we deal with duality.
\begin{prop}\label{prop:tildeIsci}
Under Assumption \eqref{eq:A-intro}, the functional $\R_+\times \mathcal{P}(\X)\times \mathcal{P}(\Y)\to \R\cup\{+\infty\} : (\lambda, \mu,\nu) \mapsto \widebar{\mathcal{I}}_{c_\lambda}(\mu,\nu)$ (with $c_\lambda$ defined as in Proposition \ref{prop:barI_csci}) is lower semicontinuous at any point $(\lambda, \mu,\nu)$ with $\mathrm{Supp}\,(\mu)=\X$.
\end{prop}
\proof
Let $(\lambda, \mu,\nu) \in \mathcal{P}(\X) \times \mathcal{P}(\Y)$ be such that $\mathrm{Supp}\,(\mu)=\X$ and consider $(\lambda_n,\mu_n,\nu_n)\in \mathcal{P}(\X)\times \mathcal{P}(\Y)$ a sequence converging to $(\lambda,\mu,\nu)$. According to Lemma \ref{lem:soltildebar}, for all $n\geq0$ there exists $\pi_n \in  \Pi( \mathrm{Supp}(\mu_n),\nu_n)$ such that
\[
\bar{I}_{c_{\lambda_n}}^{\mu_n}[\pi_n] = \widebar{\mathcal{I}}_{c_{\lambda_n}}(\mu_n,\nu_n).
\]
Let $\ell = \liminf_{n\to \infty} \bar{I}_{c_{\lambda_n}}^{\mu_n}[\pi_n]$. Extracting a subsequence if necessary, one can assume without loss of generality that
$\bar{I}_{c_{\lambda_n}}^{\mu_n}[\pi_n] \to \ell$ as $n \to \infty$. Since $\mathcal{P}(\X\times \Y)$ is compact, the sequence $\pi_n$ admits a converging subsequence, which we denote again by $\pi_n$. Let $\bar{\pi}$ be the limit of $\pi_n$. Since $\nu_n\to\nu$, $\bar{\pi} \in \Pi(\X,\nu) =\Pi(\mathrm{Supp}\,(\mu),\nu) $. Then, by semicontinuity of $\bar{I}_{c_\cdot}^{\cdot}[\,\cdot\,]$, one has
\[
\ell = \lim_{n\to \infty} \bar{I}_{c_{\lambda_n}}^{\mu_n}[\pi_n] \geq \bar{I}_{c_{\lambda}}^{\mu}[\bar{\pi}] \geq \inf_{\pi \in  \Pi( \mathrm{Supp}(\mu),\nu)} \bar{I}_{c_{\lambda}}^\mu[\pi] = \widebar{\mathcal{I}}_{c_{\lambda}}(\mu,\nu).
\]
\endproof

\subsection{Weak solutions as minimizers of $\bar{I}_c^\mu$}\label{sec:weaksol}
The following inequality is always true
\begin{equation}\label{eq:ineqbartilde}
\widebar{\mathcal{I}}_c(\mu,\nu) \leq \widetilde{\mathcal{I}}_c(\mu,\nu)  \leq \mathcal{I}_c(\mu,\nu) .
\end{equation}
Indeed, the first inequality comes from $\bar{I}_c^\mu \leq \tilde{I}_c^\mu$ (Item $(i)$ of Proposition \ref{prop:compbartilde}) and the second from $\tilde{I}_c^\mu = I_c^\mu$ on $\Pi(\ll \mu,\nu) \subset \Pi(\mathrm{Supp}(\mu),\nu)$. Because \eqref{eq:ineqbartilde} may not always hold with equality, we provide a sufficient condition for such equality to hold now.

We will say that $c$ satisfies Assumption \eqref{eq:Approx} if for all $\mu \in \mathcal{P}(\X)$ and $\nu \in \mathcal{P}(\Y)$,
\begin{align}\label{eq:Approx}
\tag{Approx}& \forall \pi \in \Pi(\mathrm{Supp}(\mu),\nu) \text{ there exists a sequence } \pi_n \in \Pi(\ll \mu,\nu) \text{ converging to } \pi\\ 
& \notag \text{ and such that }I_c^\mu[\pi_n] \to \bar{I}_c^\mu[\pi].
\end{align}
Of course, when $\pi \in \Pi(\ll \mu,\nu)$, one can choose the constant sequence $\pi_n=\pi$, $n\geq0$. Only the case $\pi \in  \Pi(\mathrm{Supp}(\mu),\nu)\setminus \Pi(\ll \mu,\nu)$ is non trivial in the above condition. Indeed, this condition is trivially satisfied when $\X$ is finite. We will see below more general sufficient conditions for \eqref{eq:Approx}.

\begin{thm}\label{thm:weaksol}
Let $c:\X \times \mathcal{M}(\Y) \to \R$ be a cost function satisfying condition \eqref{eq:A-intro} and \eqref{eq:Approx} and $\mu \in \mathcal{P}(\X)$.
\begin{itemize}
\item[$(i)$] For any $\pi \in \Pi(\mathrm{Supp}(\mu),\,\cdot\,)$, it holds $\bar{I}_c^\mu[\pi] = \tilde{I}_c^\mu[\pi]$.
\item[$(ii)$]  For any $\nu \in \mathcal{P}(\Y)$, it holds $\mathcal{I}_c(\mu,\nu) = \widetilde{\mathcal{I}}_c(\mu,\nu)= \overline{\mathcal{I}}_c(\mu,\nu)$.
\item[$(iii)$] Let $\nu \in \mathcal{P}(\Y)$ be such that $\mathcal{I}_c(\mu,\nu)<+\infty$. A coupling $\pi \in \Pi(\mathrm{Supp}(\mu),\nu)$ is a weak solution of the transport problem \eqref{eq:TP} if and only if $\pi$ minimizes $\bar{I}_c^\mu$ on $\Pi(\mathrm{Supp}(\mu),\nu)$.
\end{itemize}
\end{thm}

\proof
Item $(i)$ follows from Proposition \ref{prop:compbartilde} (Item $(iii)$) and Assumption \eqref{eq:Approx}.
Let us show Item $(ii)$. Let $\pi \in \Pi(\mathrm{Supp}(\mu),\nu)$ and consider a sequence $\pi_n \in \Pi(\ll \mu,\nu)$ converging to $\pi$ such that $I_c^\mu[\pi_n] \to \bar{I}_c^\mu[\pi]$. For all $n$, it holds $\mathcal{I}_c(\mu,\nu) \leq I_c^\mu[\pi_n]$, and so letting $n \to \infty$ gives $\mathcal{I}_c(\mu,\nu) \leq \bar{I}_c^\mu[\pi]$. Optimizing over all $\pi \in \Pi(\mathrm{Supp}(\mu),\nu)$ yields $\mathcal{I}_c(\mu,\nu) \leq \widebar{\mathcal{I}}_c(\mu,\nu)$ which together with \eqref{eq:ineqbartilde} proves the claim.
Let us finally show Item $(iii)$. Since $\widebar{\mathcal{I}}_c(\mu,\,\cdot\,) = \mathcal{I}_c(\mu,\,\cdot\,)$, it follows that any minimizer of $\bar{I}_c^\mu$ on  $\Pi(\mathrm{Supp}(\mu),\nu)=\mathrm{cl}\, \Pi(\ll\mu,\nu)$ is a weak solution. Conversely, if $\pi_n$ is a sequence of $\Pi(\ll \mu,\nu)$ converging to some $\pi$ and such that $I_c^\mu[\pi_n] \to \mathcal{I}_c(\mu,\nu)$, then by lower semicontinuity of $\bar{I}_c^\mu$, it holds $\bar{I}_c^\mu[\pi] \leq \mathcal{I}_c(\mu,\nu) =\widebar{\mathcal{I}}_c(\mu,\nu)$, and so $\pi \in \Pi(\mathrm{Supp}(\mu),\nu)$ is a minimizer of $\bar{I}_c^\mu$ on $\Pi(\mathrm{Supp}(\mu),\nu)$, which completes the proof.
\endproof
Let us now give some concrete conditions on $c$ ensuring \eqref{eq:Approx}.
We will say that a cost function $c:\X \times \mathcal{M}(\Y)\to \R$ satisfies Assumption \eqref{eq:C} if
\begin{equation}\label{eq:C}\tag{C}
\left\{\begin{aligned}
& - \text{ for all } m \in \mathcal{M}(\Y),\text{  the functions } c(\,\cdot\,,m)\text{  and }c'_\infty(\,\cdot\,,m)\text{  are continuous on } \X,\\
& \text{and} \\
& \notag - \text{  there exists }a\geq0\text{  such that }c'_\infty(x,p)\leq a\text{  for all }x \in \X\text{  and }p \in \mathcal{P}(\Y).
\end{aligned}\right.
\end{equation}
For instance, the cost function introduced in Section \ref{sec:CK}:
\[
c(x,m) =  - x \left(\int y\,dm\right)^{\eta}, \qquad x\in [\alpha,\beta], m\in \mathcal{M}([\gamma,\delta]),
\]
where $\alpha,\beta,\gamma,\delta \geq0$, $\eta \in (0,1)$, is such that
\[
c'_\infty(x,p) = 0,\qquad \qquad x\in [\alpha,\beta], p\in \mathcal{P}([\gamma,\delta]),
\]
and so $c$ satisfies Assumption \eqref{eq:C}.

\begin{lem}\label{lem:recovery}
If $c: \X \times \mathcal{M}(\Y) \to \R$ satisfies Assumptions \eqref{eq:A-intro} and \eqref{eq:C}, then it satisfies Assumption \eqref{eq:Approx}. More precisely, for any $\mu \in \mathcal{P}(\X)$ and $\pi \in \Pi(\mathrm{Supp}(\mu),\nu)$, the sequence $\pi_n \in \Pi(\ll\mu,\nu)$, $n\geq0$, defined in Lemma \ref{lem:descr} is such that $I_c^\mu(\pi_n) \to \bar{I}_c^\mu(\pi)$.
\end{lem}
The proof of Lemma \ref{lem:recovery} is presented in Section \ref{App:2} of Appendix.

\subsection{A criterion for the existence of strong solutions}
Recall Assumption \eqref{eq:B-intro} given in the introduction, which can be restated as follows:
\[
c'_\infty(x,m)=+\infty,\qquad \forall m \in \mathcal{M}(\Y)\setminus \{0\},\qquad \forall x\in \X.
\]
Under Assumption \eqref{eq:B-intro}, one gets $\bar{I}_c^\mu[\pi] = I_c^\mu[\pi]$, if $\pi \in \Pi(\ll \mu,\,\cdot\,)$ and $+\infty$ otherwise.
\begin{lem}\label{lem:BAprox}
If $c:\X \times \mathcal{M}(\Y) \to \R$ is a cost function satisfying Assumptions \eqref{eq:A-intro} and  \eqref{eq:B-intro}, then it satisfies \eqref{eq:Approx}.
\end{lem}
\proof
If $\pi_n \in \Pi(\ll \mu,\nu)$ is any sequence converging to $\pi \in  \Pi(\mathrm{Supp}(\mu),\nu)\setminus \Pi(\ll \mu,\nu)$ (such sequences always exist according to Lemma \ref{lem:descr}), then since $\bar{I}_c^\mu$ is lower semicontinuous, one gets
\[
\liminf_{n\to \infty}I_c^\mu[\pi_n] =\liminf_{n\to \infty}\bar{I}_c^\mu[\pi_n]  \geq \bar{I}_c^\mu[\pi]=+\infty
\]
and so $I_c^\mu[\pi_n] \to  \bar{I}_c^\mu[\pi]$.
\endproof

The following result shows that strong solutions always exist under Assumptions \eqref{eq:A-intro} and  \eqref{eq:B-intro}.
\begin{thm}\label{thm:strongsol}
Let $c:\X \times \mathcal{M}(\Y) \to \R$ be a cost function satisfying Assumptions \eqref{eq:A-intro} and  \eqref{eq:B-intro}. If $\mu \in \mathcal{P}(\X)$, $\nu \in \mathcal{P}(\Y)$ are such that $\mathcal{I}_c(\mu,\nu)<+\infty$, then any weak solution of the transport problem \eqref{eq:TP} is a strong solution.
\end{thm}

\proof
According to Lemma \ref{lem:BAprox} and Theorem \ref{thm:weaksol}, if  $\pi$ is a weak solution, then
\[
\bar{I}_c^\mu[\pi]=\widebar{\mathcal{I}}_c(\mu,\nu)=  \mathcal{I}_c(\mu,\nu) <\infty.
\]
Therefore, $\bar{I}_c^\mu[\pi] <+\infty$ and so $\pi \in \Pi(\ll \mu,\nu)$ and $\bar{I}_c^\mu[\pi] = I_c^\mu[\pi] =  \mathcal{I}_c(\mu,\nu)$, which shows that $\pi$ is a strong solution.
\endproof

Condition \eqref{eq:B-intro} applies for instance if there exists $\phi : \R_+ \to \R$ a function such that $\phi(u)/u \to +\infty$, when $u \to + \infty$, such that
\begin{equation}\label{eq:bound_below_c}
c(x,m) \geq \phi (m(\Y)),\qquad \forall x\in \X, \forall m \in \mathcal{M}(\Y).
\end{equation}
If $\X,\Y \subset \R^d$ and the convex hull of $\Y$ does not contain $0$, this assumption is satisfied, for instance, by the following conical cost functions
\begin{equation}\label{eq:coutquadratique}
c(x,m) = \left\| x - \int y\,dm\right\|^p,\qquad x\in \X,m \in \mathcal{M}(\Y),
\end{equation}
where $\|\,\cdot\,\|$ is an arbitrary norm on $\R^d$ and $p>1$.

\begin{rem}
Let us briefly indicate another possible method for proving existence of strong solutions when $c$ satisfies \eqref{eq:bound_below_c}. Let $\nu \in \mathcal{P}(\Y)$ be such that $\mathcal{I}_c(\mu,\nu)<+\infty$ and assume $(\pi_n)_{n\geq 0}$ is a sequence in $\Pi(\ll \mu,\nu)$ such that $I_c^\mu[\pi_n] \to \mathcal{I}_c(\mu,\nu)$ and $\pi_n\to \pi$. Then, from \eqref{eq:bound_below_c}, $\sup_{n\in \N}\int \phi(N_n)\,d\mu <+\infty$, denoting by $N_n$ the density of the first marginal of $\pi_n$ with respect to $\mu$. Therefore, the sequence $(N_n)_{n\geq 0}$ is uniformly integrable and so, according to the Dunford-Pettis theorem, admits a converging subsequence for the topology $\sigma(L^1(\mu),L^\infty(\mu))$. Hence, $\pi \in \Pi(\ll \mu,\nu)$ and is therefore a strong solution.
\end{rem}

\section{Dual formulation}\label{sec:duality}
In this Section, we establish a Kantorovich-type dual formula for the transport problem \eqref{eq:TPintro}.
The derivation of these dual forms will make use of the following abstract Fenchel-Moreau biconjugation theorem (see e.g \cite[Theorem 2.3.4]{Zal02}).
\begin{thm}\label{thm:FM}
Let $E$ be a Hausdorff locally convex topological vector space and $E'$ its topological dual space. If $F:E\to (-\infty,\infty]$ is a convex function such that $\mathrm{dom}\,(F)=\{x\in E : F(x)<\infty\} \neq \emptyset$, then for any $x \in \mathrm{dom}\,(F)$ where $F$ is lower semicontinuous, one has
\[
F(x) = \sup_{\ell \in E'} \left\{ \ell(x) - F^*(\ell)\right\}
\]
where
\[
F^*(\ell) = \sup_{x\in E} \left\{\ell(x) - F(x)\right\},\qquad \ell \in E'.
\]
\end{thm}

As noticed in Remark \ref{rem:unbalanced}, the definition of $\mathcal{I}_c$ also makes sense for arbitrary finite non-negative measures : more precisely, if $\alpha \in \mathcal{M}(\X)$ and $\beta \in \mathcal{M}(\Y)$, one defines
\[
\mathcal{I}_c(\alpha,\beta) = \inf \int c(x,q^x)\,\alpha(dx),
\]
where the infimum runs over the set of all non-negative kernels $q=(q^x)_{x\in \X}$ such that $q^x \in \mathcal{M}(\Y)$, for all $x\in X$, and $\int q^x(dy)\alpha(dx) = \beta(dy)$. As follows immediately from this definition, the functional $\mathcal{I}_c$ is positively $1$-homogenous: $\mathcal{I}_c(\lambda \alpha,\lambda \beta) = \lambda \mathcal{I}_c(\alpha,\beta)$, for all $\lambda\geq0$. One can define accordingly $\overline{\mathcal{I}}_c$.
Note that, if $\alpha(\X)>0$ and $\beta(\Y)>0$, we have the relation
\[
\overline{\mathcal{I}}_c (\alpha,\beta) = \alpha(\X) \overline{\mathcal{I}}_{c_{\lambda}} \left(\frac{\alpha}{\alpha(\X)},\frac{\beta}{\beta(\Y)}\right)
\]
with $\lambda = \frac{\beta(\Y)}{\alpha(\X)}$ and $c_\lambda(x,m) = c(x,\lambda m)$, $x\in \X$, $m\in \mathcal{M}(\Y)$.

In what follows, we apply Theorem \ref{thm:FM} using the following setting: $E = \mathcal{M}_s(\X) \times \mathcal{M}_s(\Y)$ equipped with the product weak topology whose topological dual is $E'=\mathcal{C}_b(\X)\times \mathcal{C}_b(\Y)$ and $F: \mathcal{M}_s(\X) \times \mathcal{M}_s(\Y) \to \R\cup\{+\infty\}$ defined as follows
\begin{equation}\label{eq:F}
F(\alpha,\beta) =\left\{\begin{array}{lll} \overline{\mathcal{I}}_c(\alpha, \beta) & \text{ if } \alpha,\beta \geq 0    \\
 + \infty  & \text{ otherwise}     \end{array}\right.,
\end{equation}
for all $\alpha \in \mathcal{M}_s(\X)$ and $\beta \in \mathcal{M}_s(\Y)$.
\begin{lem}\label{lem:F}
The functional $F$ is convex on $\mathcal{M}_s(\X) \times \mathcal{M}_s(\Y)$. Moreover, under Assumption \eqref{eq:A-intro}, the functional $F$ is lower semicontinuous at any point $(\alpha,\beta)$ such that $\alpha,\beta\geq0$ and $\mathrm{Supp}\,(\alpha)=\X.$
\end{lem}
\proof
The first statement easily follows from Proposition \ref{prop:convexity}. Let $\alpha_n,\beta_n$ be sequences converging respectively to finite nonnegative measures $\alpha,\beta$ such that $\mathrm{Supp}\,(\alpha)=\X$ and $\beta(\Y)>0$ let us show that $\liminf_{n\to \infty} F(\alpha_n,\beta_n) \geq F(\alpha,\beta)$.  As $\alpha$ has full support, $\alpha(\X)>0$ and since $\alpha_n(\X) \to \alpha(\X)$, it follows that $\alpha_n(\X)>0$ for all $n$ large enough. Since $\alpha_n/\alpha_n(\X) \to \alpha/\alpha(\X)$ (which has full support),  $\beta_n/\beta_n(\Y) \to \beta/\beta(\Y)$, and $\lambda_n = \beta_n(\Y)/\alpha_n(\X) \to  \lambda = \beta(\Y)/\alpha(\X)$, it follows from Proposition \ref{prop:tildeIsci} that
\[
\liminf_{n\to \infty} \overline{\mathcal{I}}_{c_{\lambda_n}}\left(\frac{\alpha_n}{\alpha_n(\X)},\frac{\beta_n}{\beta_n(\Y)} \right) \geq \overline{\mathcal{I}}_{c_\lambda}\left(\frac{\alpha}{\alpha(\X)},\frac{\beta}{\beta(\Y)} \right)
\]
which proves the claim when $\beta(\Y)>0$. The case $\beta(\Y)=0$ is left to the reader.
\endproof

\begin{thm}\label{thm:duality2}
Under Assumptions \eqref{eq:A-intro} and \eqref{eq:Approx}, the following holds
\begin{equation}\label{eq:duality-2}
\mathcal{I}_c(\mu,\nu) = \sup_{f\in \mathcal{C}_b(\Y)}\left\{  \int K_{c} f(x)\,\mu(dx) - \int f(y)\,\nu(dy)\right\},\qquad \forall \nu \in \mathcal{P}(\Y)
\end{equation}
where
\[
K_{c} f(x) = \inf_{m \in \mathcal{M}(\Y)} \left\{\int f\,dm + c(x,m)\right\},\qquad x\in \X.
\]
In particular, \eqref{eq:duality-2} holds whenever $c$ satisfies Assumption \eqref{eq:A-intro} and Assumption \eqref{eq:B-intro} or  \eqref{eq:C}.
\end{thm}
\begin{rem}
Obtaining general sufficient conditions for dual attainment in Theorem \ref{thm:duality2} would be useful and lead to a cyclical monotonicity criterium characterizing optimality of transport plans, in the spirit of the $C$-monotonicity criterium obtained for WOT \cite{BBP19,BP22b}.
\end{rem}

The proof below is adapted from the proof of \cite[Theorem 4.2]{ABC19}.

\textbf{Note.} In the published version of this paper (SIAM Journal of Mathematical Analysis, 2023) the proof of Theorem \ref {thm:duality2} contained an error, which is now fixed in the present version.

\proof
Let $\mu \in \mathcal{P}(\X)$, $\nu\in \mathcal{P}(\Y)$ and assume that $\mu$ has full support.
According to Lemma \ref{lem:F}, the function $F$ defined by \eqref{eq:F} is convex on $\mathcal{M}_s(\X)\times \mathcal{M}_s(\Y)$ and lower semicontinuous at $(\mu,\nu)$. Therefore, according to Theorem \ref{thm:weaksol} (Item $(ii)$) and Theorem \ref{thm:FM}, it holds
\[
\mathcal{I}_c(\mu,\nu) = F(\mu,\nu) = \sup_{(\varphi,\psi) \in \mathcal{C}_b(\X)\times \mathcal{C}_b(\Y)} \left\{\int \varphi\,d\mu + \int \psi\,d\nu - F^*(\varphi,\psi)\right\},
\]
with, for all $(\varphi,\psi) \in \mathcal{C}_b(\X)\times \mathcal{C}_b(\Y)$,
\begin{align*}
F^*(\varphi,\psi) &= \sup_{(\alpha,\beta)\in \mathcal{M}_s(\X)\times \mathcal{M}_s(\Y)}\left\{\int \varphi\,d\alpha + \int \psi\,d\beta - F(\alpha,\beta)\right\}\\
& = \sup_{(\alpha,\beta)\in \mathcal{M}(\X)\times \mathcal{M}(\Y)} \sup_{\lambda \geq0}\left\{\lambda \int \varphi\,d\alpha + \lambda \int \psi\,d\beta - \lambda \overline{\mathcal{I}}_c(\alpha,\beta)\right\}\\
& = \chi_{\mathcal{K}_c} (\varphi,\psi),
\end{align*}
where
\[
\mathcal{K}_c = \left\{(\varphi,\psi)\in \mathcal{C}_b(\X)\times \mathcal{C}_b(\Y) : \int \varphi\,d\alpha + \int \psi\,d\beta \leq \overline{\mathcal{I}}_c(\alpha,\beta), \forall \alpha \in \mathcal{M}(\X),\forall \beta \in \mathcal{M}(\Y) \right\}
\]
and $\chi_{\mathcal{K}_\alpha} (\varphi,\psi) = 0$ if $(\varphi,\psi)\in \mathcal{K}_c$ and $+\infty$ otherwise.
Thus, we get
\[
\mathcal{I}_c(\mu,\nu) = F(\mu,\nu) = \sup_{(\varphi,\psi) \in \mathcal{K}_c} \left\{\int \varphi\,d\mu + \int \psi\,d\nu \right\},
\]
Now, observe that if $(\varphi,\psi) \in \mathcal{K}_c$, then (by choosing $\alpha = \delta_x$, with $x\in \X$) one gets
\[
\varphi(x) \leq \inf_{\beta \in \mathcal{M}(\Y)} \left\{-\int \psi\,d\beta + \overline{\mathcal{I}}_c(\delta_x,\beta)\right\} \leq \inf_{\beta \in \mathcal{M}(\Y)} \left\{-\int \psi\,d\beta + c(x,\beta)\right\} = K_c(-\psi)(x),
\]
where we used that $\overline{\mathcal{I}}_c(\delta_x,\beta) \leq \mathcal{I}_c(\delta_x,\beta) = c(x,\beta)$.
Thus,
\[
\mathcal{I}_c(\mu,\nu) \leq  \sup_{f\in \mathcal{C}_b(\Y)}\left\{  \int K_{c} f(x)\,\mu(dx) - \int f(y)\,\nu(dy)\right\}.
\]
The converse inequality is always true. Indeed, if $\pi(dxdy) = N(x)\mu(dx)\pi^x(dy) \in \Pi(\ll \mu,\nu)$, then
\[
\int K_cf(x)\,\mu(dx) \leq \int \left(\int f(y) d(N(x)\pi^x)(dy)+c(x,N(x)\pi^x)\right)\,\mu(dx) = \int f\,d\nu + I_c^\mu[\pi].
\]
Formula \eqref{eq:duality-2} is thus proved when $\mu$ has full support. When $\mu$ does not have full support, letting $\tilde{\X}=\mathrm{Supp}(\mu)$ and applying the preceding reasoning in the space $\mathcal{M}_s(\tilde{\X})\times \mathcal{M}_s(\Y)$ gives the desired duality formula.
\endproof

\section{Monotonicity properties and uniqueness of primal solutions}\label{sec:monotonicity}
In this Section, we consider cost functions of the following form
\begin{equation}\label{eq:cost2}
c(x,m) = G\left(\int F(x,y)\,m(dy)\right),\qquad x\in \X, m\in \mathcal{M}(\Y),
\end{equation}
where $F : \X \times \Y \to (0,+\infty)$ is a continuous function and $G: \R^+ \to \R$ is a convex function, assumed to be differentiable on $(0,+\infty)$ and we denote by $G'(0)= \lim_{x \to 0^+} G'(x) \in \R \cup \{-\infty\}$ and $G'(+\infty) =  \lim_{x \to +\infty} G'(x) \in \R \cup \{+\infty\}$. We establish below that the dual problem admits a solution (Theorem \ref{thm:existencedualopt}) and then use this dual optimizer to  get information on the support of primal solutions (Proposition \ref{prop:dualprimal}). Finally, we consider the particular case when $\X$ and $\Y$ are subsets of $\R$ and prove uniqueness of primal solutions under suitable assumptions on $F, G$ and $\mu$ (Theorem \ref{thm:uniqueness1d}).

First, let us check that $c$ satisfies the assumptions introduced in the preceding Sections. Writing $G$ as a countable supremum of affine functions, one easily sees that $c$ satisfies Assumption \eqref{eq:A-intro} and in particular \eqref{eq:cbound}. Thus $\mathcal{I}_c(\mu,\nu)$ makes sense for any $\mu \in \mathcal{P}(\X)$ and $\nu \in \mathcal{P}(\Y)$. In addition, for any $m \in \mathcal{M}(\Y)\setminus \{0\}$ and $x\in \X$
\[
c'_\infty(x,m) = G'(+\infty).
\]
So $c$ satisfies Assumption \eqref{eq:B-intro} if $G'(+\infty) = +\infty$ and Assumption \eqref{eq:C} otherwise. Therefore, Theorem \ref{thm:duality2} applies and it is easily seen that
\begin{equation}\label{eq:dualityL1}
\mathcal{I}_c(\mu,\nu) =  \sup_{f\in \mathcal{L}^1(\nu)}\left\{  \int K_{c} f(x)\,\mu(dx) - \int f(y)\,\nu(dy)\right\},\qquad \forall \mu \in \mathcal{P}(\X), \forall \nu \in \mathcal{P}(\Y),
\end{equation}
with
\[
\mathcal{L}^1(\nu) = \left\{f : \Y \to \R : f \text{ measurable and } \int |f|\,d\nu<+\infty\right\}
\]
and
\begin{equation}\label{eq:K_cgen}
K_cf(x) = \inf_{m \in \mathcal{M}(\Y) \text{ s.t } f \in L^1(m)} \left\{\int f\,dm + c(x,m)\right\}, \qquad x \in \X, f \in \mathcal{L}^1(\nu).
\end{equation}
Because $K_cf$ is upper semicontinuous on $\X$ as an infimum of continuous functions and $K_c f(x) \leq G(0)$ for all $x \in \X$, thus $\int K_cf(x)\,\mu(dx)$ always makes sense.

The following result establishes dual attainment.

\begin{thm}\label{thm:existencedualopt}
If the function $G$ in \eqref{eq:cost2} is such that $G'(0) >-\infty$, then for every $\mu \in \mathcal{P}(\X)$ and $\nu \in \mathcal{P}(\Y)$, there exists a bounded function $\bar{f}$ on $\Y$ such that
\[
\mathcal{I}_c(\mu,\nu) =   \int K_{c} \bar f(x)\,\mu(dx) - \int \bar f(y)\,\nu(dy).
\]
The same conclusion holds if $G'(0)=-\infty$ and $\Y$ is a finite set.
Moreover, if $G$ is non-decreasing (resp. non-increasing) on $\R^+$, then $\bar{f}$ can be chosen nonpositive (resp. nonnegative).
\end{thm}
\proof
Let us show that the supremum in \eqref{eq:dualityL1} can be restricted to the class of measurable functions $f$ such that $f \geq a$, with $a = \mathcal{I}_c(\mu,\nu) -\sup_{u\in \X, v\in \Y} G(2F(u,v))-1$.
Note that if $f$ is not bounded from below, then $K_c(f)(x)=-\infty$ for all $x\in \X$, so the supremum in \eqref{eq:dualityL1} can be restricted to functions $f$ bounded from below. If $f$ is such a function, then for all $y \in \Y$, one gets
\[
K_cf(x) \leq 2 f(y) + G(2 F(x,y)) \leq 2 f(y) + \sup_{u\in \X, v\in \Y} G(2F(u,v)).
\]
So optimizing over $y$, one obtains
\[
K_cf(x) \leq 2 \inf f + \sup_{u\in \X, v\in \Y} G(2F(u,v))
\]
and so
\[
\int K_c f(x)\,\mu(dx) - \int f(y)\,\nu(dy) \leq  \inf f + \sup_{u\in \X, v\in \Y} G(2F(u,v))
\]
Therefore, if $\inf f < \mathcal{I}_c(\mu,\nu) -\sup_{u\in \X, v\in \Y} G(2F(u,v))-1$, then
\[
\int K_c f(x)\,\mu(dx) - \int f(y)\,\nu(dy) <  \mathcal{I}_c(\mu,\nu) -1,
\]
and so $f$ can be dropped from the supremum in \eqref{eq:dualityL1}. We thus conclude that the supremum in \eqref{eq:dualityL1} can be restricted to functions $f$ bounded from below by $a$.
In the case, where $G$ is non-increasing, this lower bound can be improved. Indeed, if $f(y_0)<0$ for some $y_0 \in \Y$, then for all $\lambda >0$ one has
\[
K_cf(x) \leq \lambda f(y_0) + G(\lambda F(x,y_0)) \to - \infty
\]
as $\lambda \to +\infty$. So the supremum in \eqref{eq:dualityL1} can be restricted in this case to nonnegative functions.

Now, let us show that the supremum in \eqref{eq:dualityL1} can be further restricted to functions $f$ such that $f \leq b$, where $b = [G'(0)]_-\sup_{x\in \X,y\in \Y} F(x,y)$ with $[x]_- = \max(-x ; 0)$. Let $f \in \mathcal{L}^1(\nu)$; define $A = \{y \in \Y : f(y) \leq b\}$ and for all $m \in \mathcal{M}(\Y)$ write $m_A(dy) = \mathbf{1}_A(y)\,m(dy)$ and $m_{A^c}(dy) = \mathbf{1}_{A^c}(y)\,m(dy)$. Since $u\mapsto G(u)+[G'(0)]_-u$ is non-decreasing, for all $x \in \X$, one has
\begin{align*}
\int  \min(f,b)\,dm + c(x,m) &= \int f\,dm_A +bm(A^c)+ G\left(\int F(x,y)\,m_A(dy) + \int F(x,y)\,m_{A^c}(dy)\right)\\
& \geq \int f\,dm_A + bm(A^c)+ G\left(\int F(x,y)\,m_A(dy) \right) - [G'(0)]_-\int F(x,y)\,m_{A^c}(dy)\\
& \geq  \int f\,dm_A +  G\left(\int F(x,y)\,m_A(dy) \right) +\left(b- [G'(0)]_-\sup_{x\in \X,y\in \Y} F(x,y)\right)m(A^c)\\
& \geq K_cf(x)
\end{align*}
and so, letting $\hat f = \min(f,b)$, one gets $K_c \hat f  \geq K_c f$. On the other hand, since $\hat f  \leq f$, it also holds that $K_c \hat f \leq K_c f$ and so $K_c f = K_c \hat f $. Since,
\[
\int K_c f\,d\mu - \int f\,d\nu \leq \int K_c \hat f\,d\mu - \int \hat f\,d\nu
\]
one concludes that the supremum in \eqref{eq:dualityL1} can be restricted to functions bounded from above by $b$. In particular, when $G$ is non-decreasing, then $[G'(0)]_- = 0$ and one can restrict attention to nonpositive functions.

Let us now show dual attainment. Consider a sequence $(g_n)_{n\geq0}$ of elements of $\mathcal{B} = \{ f \in \mathcal{L}^1(\nu): a\leq f\leq b\}$ such that $\int K_cg_n\,d\mu - \int g_n\,d\nu \to \mathcal{I}_c(\mu,\nu)$.
According to the Dunford-Pettis theorem (or the Banach-Alaoglu-Bourbaki theorem), one can extract from $(g_n)_{n\geq0}$ a subsequence (still denoted $(g_n)_{n\geq0}$) converging to some $g_\infty \in \mathcal{B}$ for the weak topology $\sigma(L^1,L^\infty)$: for all $h \in L^\infty(\nu)$, $\int g_n h\,d\nu \to \int g_\infty h\,d\nu$. Moreover, according to Mazur's lemma, there exists a sequence $(f_n)_{n\geq0}$ of the form $f_n = \sum_{i=0}^{N_n} \lambda_i^{(n)} g_{n+i}$ with $N_n \geq 0$, $\lambda_0^{(n)},\ldots,\lambda_{N_n}^{(n)}\geq0$ and $\sum_{i=0}^{N_n} \lambda_i^{(n)}=1$ such that $f_n$ converges strongly in $L^1(\nu)$ to $g_\infty$, as $n\to \infty.$ Extracting a subsequence if necessary, one can further assume that $(f_n)_{n\geq0}$ converges $\nu$ almost everywhere to $g_\infty$. Let $\Phi : \mathcal{B} \to \R : f \mapsto  \int K_cf\,d\mu - \int f\,d\nu$. One sees that
\[
K_c((1-t)f+tg) \geq (1-t)K_cf+tK_cg,\qquad \forall t\in [0,1],\forall f,g \in \mathcal{B}
\]
and so $\Phi$ is concave. Therefore
\[
\Phi(f_n) \geq \sum_{i=0}^{N_n} \lambda_i^{(n)} \Phi(g_{n+i}) \geq \inf_{k\geq n} \Phi(g_k) \to \mathcal{I}_c(\mu,\nu),
\]
as $n\to \infty$, and so $\Phi(f_n) \to \mathcal{I}_c(\mu,\nu)$ as $n\to \infty.$ Since $K_c f_n \leq G(0)$ for all $n\geq0$, one can apply Fatou's Lemma
\begin{align*}
\mathcal{I}_c(\mu,\nu) = \limsup_{n\to +\infty} \int K_cf_n\,d\mu - \lim_{n\to +\infty}\int f_n\,d\nu \leq \int \limsup_{n\to +\infty} K_cf_n\,d\mu - \int g_\infty\,d\nu.
\end{align*}
For all $m\in \mathcal{M}(\Y)$ and $x\in \X$, it holds
\[
K_cf_n(x) \leq \int f_n\,dm + c(x,m)
\]
and so, applying Fatou's Lemma again, one gets
\[
\limsup_{n\to +\infty} K_cf_n(x) \leq \int \limsup_{n\to +\infty}f_n\,dm + c(x,m),\qquad x\in \X
\]
and so, optimizing over $m$ yields $\limsup_{n\to +\infty} K_cf_n \leq K_c(\bar{f})$, with $\bar{f} = \limsup_{n\to +\infty}f_n \in \mathcal{B}$, and so
\[
\mathcal{I}_c(\mu,\nu) \leq  \int K_c \bar{f}\,d\mu - \int g_\infty\,d\nu.
\]
Finally, since $f_n$ converges $\nu$ almost everywhere to $g_\infty$, it holds $g_\infty = \bar{f}$ $\nu$ a.e. and so
\[
\mathcal{I}_c(\mu,\nu) \leq  \int K_c \bar{f}\,d\mu - \int \bar{f}\,d\nu,
\]
which shows that $\bar{f}$ is a dual optimizer.
Finally, let us show dual attainment when $\Y$ is a finite set and $G'(0) \geq - \infty$. According to what precedes, the supremum in \eqref{eq:dualityL1} can be restricted to functions $f \geq a$. On the other hand, since $\int K_c f\,d\mu \leq G(0)$, one can further restrict the supremum in \eqref{eq:dualityL1} to functions $f$ such that $\int f\,d\nu \leq G(0)+1-\mathcal{I}_c(\mu,\nu) := a'$. Since $\Y$ is finite, the set $\mathcal{C} = \{f\in \mathcal{L}^1(\nu) : a \leq f, \int f\,d\nu \leq a' \}$ is compact. Reasoning as above, one shows that any maximizing sequence of the dual problem admits a subsequence converging to a dual optimizer, which completes the proof.
\endproof

The following result relates primal and dual optimizers (provided they exist).

\begin{prop}\label{prop:dualprimal}
Let  $\mu \in \mathcal{P}(\X)$ and $\nu \in \mathcal{P}(\Y)$ be such that $\mathcal{I}_c(\mu,\nu) <+\infty$ and suppose that $\bar{q}$ is a kernel minimizer of $\mathcal{I}_c(\mu,\nu)$ and that $\bar{f} \in \mathcal{L}^1(\nu)$ is a dual optimizer:
\[
\mathcal{I}_c(\mu,\nu) = I_c^\mu[\bar{q}] = \int K_c\bar{f}\,d\mu - \int \bar{f}\,d\nu.
\]
Then, the following relation holds true: for $\mu$ almost all $x \in \X$,
\begin{equation}\label{eq:condopt}
G'\left(\int F(x,z)\,\bar{q}^x(dz)\right)F(x,y) + \bar{f}(y) \geq 0,\qquad \forall y\in \Y.
\end{equation}
In particular, if $G'(0)=-\infty$, then $\bar{q}^x(\Y)>0$ for $\mu$ almost all $x \in \X$.
Moreover, equality holds in \eqref{eq:condopt} for $\bar{q}^x$ almost all $y \in \Y$.
\end{prop}
\proof
Since
\begin{align*}
\mathcal{I}_c(\mu,\nu) &=  \int K_c\bar{f}\,d\mu - \int \bar{f}\,d\nu
 \leq \int \int \bar{f}(y)\,\bar{q}^x(dy)\mu(dx) + \int c(x,\bar{q}^x)\,\mu(dx) - \int \bar{f}\,d\nu
 = I_c^\mu[\bar{q}]
 = \mathcal{I}_c(\mu,\nu),
\end{align*}
one concludes that
\[
\int \bar{f}(y)\,\bar{q}^x(dy) + c(x,\bar{q}^x) = K_c\bar{f}(x),
\]
for $\mu$ almost every $x \in \X$. Fix some $x \in \X$ for which the equality holds. By definition of $K_c$, for all $t \in (0,1)$ and $m \in \mathcal{M}(\Y)$ such that $\int |f|\,dm <+\infty$, one has
\[
\int \bar{f}(y)\,\bar{q}^x(dy) + c(x,\bar{q}^x) \leq (1-t)\int \bar{f}(y)\,\bar{q}^x(dy) + t \int \bar{f}(y)\,m(dy) + c(x,(1-t)\bar{q}^x+tm).
\]
Sending $t \to 0$ yields
\[
\int \bar{f}(y)\,(\bar{q}^x-m)(dy)  \leq G'\left(\int F(x,z)\,\bar{q}^x(dz)\right) \int F(x,y)\,(m-\bar{q}^x)(dy)
\]
(with the convention $-\infty \times 0=0$ if $G'(0)=-\infty$ and $m=\bar{q}^x$). In particular, if $G'(0)=-\infty$, then $\bar{q}^x \neq 0$. Rearranging terms yields 
\begin{equation}\label{eq:condopt2}
\int \bar{f}(y) + G'\left(\int F(x,z)\,\bar{q}^x(dz)\right) F(x,y)\,\bar{q}^x(dy) \leq \int \bar{f}(y) + G'\left(\int F(x,z)\,\bar{q}^x(dz)\right) F(x,y)\,m(dy).
\end{equation}
If the function $\psi_x: \Y \to \R : y\mapsto \bar{f}(y) + G'\left(\int F(x,z)\,\bar{q}^x(dz)\right) F(x,y)$ takes a negative value at some point $y_o \in \Y$ then choosing $m= \lambda \delta_{y_o}$ with $\lambda>0$ arbitrary large yields a contradiction. Therefore, \eqref{eq:condopt} holds true. Since the function $\psi_x$ is nonnegative, one has $\inf_{m\in \mathcal{M}(\Y)}\int \psi_x(y)\,m(dy)$ = 0. Therefore taking the infimum over $m$ in \eqref{eq:condopt2} yields $\int \psi_x(y)\,\bar{q}^x(dy)=0$ and so $\psi_x(y)=0$ for $\bar{q}^x$ almost all $y\in \Y$.
\endproof

In a one dimension framework, we now draw from \eqref{eq:condopt} monotonicity properties of the supports of $\bar{q}^x$, $x \in \X.$
\begin{thm}\label{thm:uniqueness1d}Let $\X$ and $\Y$ be two compact subsets of $\R$, $\mu\in \mathcal{P}(\X), \nu \in \mathcal{P}(\Y)$ and suppose that $\bar{q} \in \mathcal{Q}(\mu,\nu)$ is a kernel solution of the transport problem \eqref{eq:TPintro} with cost function \eqref{eq:cost2}. For all $x \in \X$, denote by $S_x \subset \Y$ the support of $\bar{q}^x$.
Suppose also that $F: \R \times \R \to (0,+\infty)$ is twice continuously differentiable and such that
\begin{equation}\label{eq:monotonicityF}
\frac{\partial ^2 \ln F(x,y)}{\partial x \partial y} <0,\qquad \forall x,y\in \R.
\end{equation}
\begin{enumerate}
\item If $G$ is increasing on $\R^+$, then there exists $A \subset \X$ with $\mu(A)=1$ such that
\begin{equation}\label{eq:monotone1}
x_1 <x_2, x_1,x_2 \in A \Rightarrow \forall y_1 \in S_{x_1}, \forall y_2 \in S_{x_2}, y_1\leq y_2.
\end{equation}
\item If $G$ is decreasing on $\R^+$ and $G'(0)>-\infty$, then there exists $A \subset \X$ with $\mu(A)=1$ such that
\begin{equation}\label{eq:monotone2}
x_1 <x_2, x_1,x_2 \in A \Rightarrow \forall y_1 \in S_{x_1}, \forall y_2 \in S_{x_2}, y_1\geq y_2.
\end{equation}
The same conclusion holds if $G'(0)=-\infty$ provided $\Y$ is a finite set.
\end{enumerate}
If $\mu$ has no atoms and full support, then there exists a unique right-continuous map $\bar T : \X \to \Y$ which is non-decreasing when $G$ is increasing and non-increasing when $G$ is decreasing, such that any kernel solution $\tilde{q} \in \mathcal{Q}(\mu,\nu)$ of the transport problem \eqref{eq:TPintro} can be written as
\[
\tilde{q}^x(dy) = \tilde{N}(x)\delta_{\bar T(x)},
\]
for $\mu$ almost all $x \in \X$, with $\tilde{N}:\X \to \R^+$ a density with respect to $\mu.$ If $G$ is assumed to be strictly convex, the density $\tilde{N}$ is also unique.

If $F$ is such that
\[
\frac{\partial ^2 \ln F(x,y)}{\partial x \partial y} >0,\qquad \forall x,y\in \R.
\]
all conclusions are reversed: \eqref{eq:monotone1} holds when $G$ is decreasing, \eqref{eq:monotone2} holds when $G$ is increasing, and the monotonicity of $\bar T$ is the opposite of that of $G$.
\end{thm}
The existence of a kernel solution in Proposition \ref{prop:dualprimal} or Theorem \ref{thm:uniqueness1d} is, at least, granted in the two following cases: $G'(+\infty) = +\infty$ (according to Theorem \ref{thm:strongsol}) or $G'(+\infty)<+\infty$ and $\X$ is finite (according to Theorem \ref{thm:strongfinite}).
\proof
The proof is only provided when $G$ is increasing, proof for the other case being similar. According to Theorem \ref{thm:existencedualopt}, there exists a nonpositive bounded function $\bar{f}$ achieving equality in the dual formula for $\mathcal{I}_c(\mu,\nu)$. According to Proposition \ref{prop:dualprimal}, there is a set $A \subset \X$ with $\mu(A)=1$ and such that for all $x \in A$ it holds
\[
G'\left(\int F(x,z)\,\bar{q}^x(dz)\right)F(x,y) + \bar{f}(y) \geq 0,\qquad \forall y\in \Y.
\]
Denote by $\tilde{S}_x$ the set of $y\in \Y$ for which the inequality above is an equality. According to Proposition \ref{prop:dualprimal}, we know that $q^x (\tilde{S}_x)=1$ for all $x \in A.$
Condition \eqref{eq:monotonicityF} implies the following monotonicity property for $F$: if $a_1 <a_2$ and $b_1<b_2$ then
\begin{equation}\label{eq:monotonicityF2}
F(a_1,b_1)F(a_2,b_2)<F(a_1,b_2)F(a_2,b_1).
\end{equation}
Let $x_1,x_2 \in A$ such that $x_1<x_2$ and suppose that there exist $y_1 \in \tilde{S}_{x_1}$ and $y_2 \in \tilde{S}_{x_2}$ such that $y_2 <y_1$. Then, denoting by $U(x) = \int F(x,z)\,\bar{q}^x(dz)$, $x\in \X$, it holds
\begin{align*}
G'(U(x_1))F(x_1,y_1) &= - \bar{f}(y_1)\\
G'(U(x_2))F(x_2,y_2) &= - \bar{f}(y_2)\\
G'(U(x_1))F(x_1,y_2) &\geq - \bar{f}(y_2)\\
G'(U(x_2))F(x_2,y_1) &\geq - \bar{f}(y_1).
\end{align*}
Multiplying the last two inequalities (note that all quantities are nonnegative) one gets
\[
G'(U(x_1))G'(U(x_2))F(x_1,y_2)F(x_2,y_1) \geq  \bar{f}(y_1) \bar{f}(y_2) = G'(U(x_1))G'(U(x_2))F(x_1,y_1)F(x_2,y_2).
\]
Since $G$ is increasing the term $G'(U(x_1))G'(U(x_2))$ is positive and can be simplified yielding to
\[
F(x_1,y_2)F(x_2,y_1) \geq F(x_1,y_1)F(x_2,y_2),
\]
which contradicts \eqref{eq:monotonicityF2} with $a_1=x_1$, $a_2=x_2$, $b_1=y_2$ and $b_2=y_1$. Therefore, the family of sets $(\tilde{S}_x)_{x\in A}$ satisfies the following property:
\[
x_1 <x_2 \Rightarrow y_1 \leq y_2, \forall y_1 \in \tilde{S}_{x_1}, \forall y_2 \in \tilde{S}_{x_2}.
\]
Since $\bar{q}^{x}(\tilde{S}_{x})=1$, it is clear that $\tilde{S}_{x}$ is dense in $S_{x}$ and so the same property is satisfied by the family $(S_x)_{x\in A}$, which proves \eqref{eq:monotone1}.

Let us now assume that $\mu$ has no atoms. Let $\mathrm{is} (A)$ be the set of isolated points of $A$. Since this set is at most countable and $\mu$ has no atoms, then $\mu(\mathrm{is} (A))=0$. Thus, letting $A' = A \setminus \mathrm{is} (A)$, one gets $\mu(A')=1.$ Consider the maps $T_-, T_+ : A \to \R$ defined by $T_-(x) = \inf S_x$ and $T_+(x) = \sup S_x$ for all $x \in A'$. According to \eqref{eq:monotone1}, the maps $T_\pm$ are non-decreasing. Let $D \subset A'$ be the set of points where $T_-$ or $T_+$ is discontinuous. It is well-known that $D$ is at most countable. Hence, whenever $T_-(x)<T_+(x)$ then $x \in D$. Defining $A'' = A'\setminus D$, and $T(x)=T_-(x)=T_+(x)$ for $x \in A''$, one gets that $S_x = \{T(x)\}$ for all $x \in A''$ and so there exists some nonnegative number $\bar{N}(x)$ such that $\bar{q}^x = \bar{N}(x) \delta_{T(x)}$ for all $x \in A''$.
Finally, defining $\bar{T}(x) = \inf \{T(y) : y\in A'', y>x\}$, $x \in \X$, yields a non-decreasing and right-continuous extension of $T$ to the whole space $\X$.

To prove the uniqueness part, we use a classical reasoning which goes back to \cite{Br87}. Suppose that $\tilde{q} \in \mathcal{Q}(\mu,\nu)$ is another kernel solution of the transport problem. According to what precedes, there exists a density $\tilde{N}$ and a non-decreasing right-continuous map $\tilde{T}$ such that $\tilde{q}^x = \tilde{N}(x)\delta_{\tilde{T}(x)}$ for $\mu$ almost all $x \in \X$. By convexity, the nonnegative kernel $\frac{1}{2}(\bar q + \tilde{q})$ is also a solution. Therefore, there exists yet another non-decreasing right-continuous map $U$ such that
\[
\frac{1}{2}\left(\bar{N}(x)\delta_{\bar{T}(x)} + \tilde{N}(x)\delta_{\tilde{T}(x)}\right) =\frac{1}{2}(\bar {q}^x + \tilde{q}^x) = \left(\frac{1}{2}\bar{N}(x) +  \frac{1}{2}\tilde{N}(x)\right) \delta_{U(x)},
\]
for $\mu$ almost all $x \in \X.$ Therefore, $\tilde{T}(x) = \bar{T}(x) =U(x)$ for $\mu$ almost all $x \in \X$. Since $\X$ is the support of $\mu$ and $\tilde{T},\bar{T}$ are right-continuous, the equality holds true for all $x \in \X$.
Furthermore, by optimality one has
\[
\int G\left((\frac{1}{2}\bar{N}(x) +  \frac{1}{2}\tilde{N}(x)) F(x,\bar{T}(x))\right)\,\mu(dx) = \frac{1}{2}\int G\left(\bar{N}(x)F(x,\bar{T}(x))\right)\,\mu(dx) +\frac{1}{2}\int G\left(\tilde{N}(x)F(x,\tilde{T}(x))\right)\,\mu(dx)
\]
and so, assuming that $G$ is strictly convex, one gets $\bar{N}(x)=\tilde{N}(x)$ for $\mu$ almost every $x\in \X$, which completes the proof.
\endproof

\section{The particular case of conical cost functions}\label{sec:4}
This Section is devoted to the study of the transport problem \eqref{eq:TPintro}, when $c$ is a conical cost function.
We first obtain an improved duality result showing that under mild conditions there is dual attainment. Then, we obtain a Strassen-type result for a variant of the convex order involving positively $1$-homogenous functions. Finally, we prove structure results for primal and dual solutions.
\subsection{Framework}\label{sec:frameworkconical}
In the whole section, we adopt the following framework:
\begin{itemize}
\item $\X$ is a compact metrizable space,
\item $\Y$ is a compact subset of $\R^d$, equipped with some arbitrary norm $\|\,\cdot\,\|$,
 and we denote by $\mathrm{co}(\Y)$ its convex hull and by $\Z$ its conical hull, i.e
\[
\Z= \left\{\sum_{i=1}^n \lambda_i y_i : \lambda_1,\ldots,\lambda_n \in \R_+, y_1,\ldots,y_n \in \Y, n \geq1 \right \},
\]
\item the cost function $c:\X \times \mathcal{M}(\Y) \to \R$ is of the following form
\begin{equation}\label{eq:conical}
c(x,m) = F\left(x, \int y\,dm\right),\qquad x\in \X,m\in \mathcal{M}(\Y),
\end{equation}
where $F:\X\times \Z \to \R$ is lower semicontinuous on $\X\times \Z$ and convex with respect to its second variable.
\end{itemize}
When $c$ is of this form we will say that $c$ is a \emph{conical cost function}.

First let us translate Assumptions \eqref{eq:A-intro}, \eqref{eq:B-intro} and \eqref{eq:C} in this framework. Let us introduce the recession function of $F$, defined by
\[
F'_\infty(x,z) = \lim_{\lambda \to +\infty} \frac{F(x,\lambda z)}{\lambda},\qquad x\in \X, z\in \Z.
\]
\begin{itemize}
\item Assumption \eqref{eq:A-intro} is fulfilled as soon as $F$ satisfies the following condition \eqref{eq:A'}: there exists a family of continuous functions $(a_k)_{k\geq0}$ on $\X$ and a family of continuous functions $(u_k)_{k\geq0}$ on $\X$ with values in $\R^d$ such that
\begin{equation}\label{eq:A'}\tag{A'}
F(x,z) = \sup_{k\geq0} \left\{u_k(x)\cdot z + a_k(x)\right\},\qquad x\in \X, z\in \Z.
\end{equation}
If $F$ satisfies \eqref{eq:A'}, then the corresponding cost function $c$ satisfies \eqref{eq:A-intro} with $b_k(x,y) = u_k(x)\cdot y$ and the same $a_k$ for $k\geq 0$.

\item Assumption \eqref{eq:B-intro} is satisfied by $c$ if and only if $0\notin \mathrm{co}(\Y)$ and  $F$ satisfies the following condition \eqref{eq:B'}
\begin{equation}\label{eq:B'}\tag{B'}
F'_\infty(x,z) = +\infty,\qquad \forall  x\in \X, \forall z\in \Z\setminus\{0\}.
\end{equation}

\item Finally, Assumption \eqref{eq:C} holds for $c$ as soon as $F$ satisfies the following condition \eqref{eq:C'}:
\begin{align}\label{eq:C'}
\notag& \text{for all } z \in \Z, \text{ the functions } x\mapsto F(x,z) \text{ and } x\mapsto F'_\infty(x,z) \text{ are continuous on } \X \\
\tag{C'}& \text{and}\\
\notag& \text{there exists } a \in \R_+ \text{ such that } F'_\infty(x,z) \leq a \text{ for all  } x\in \X \text{ and }z\in \mathrm{co}(\Y).
\end{align}
\end{itemize}

\subsection{Duality and dual attainment for conical cost functions}

The following result improves upon the conclusion of Theorem \ref{thm:duality2} in the case of conical cost functions. Recall that a function $\varphi : \Z \to \R \cup \{+\infty\}$ defined on a cone $\Z \subset \R^d$ is positively $1$-homogenous if $\varphi(\lambda x) = \lambda \varphi(x)$, for all $\lambda \geq0$ and for all $x\in \Z$ such that $\varphi(x)<+\infty$.

\begin{thm}\label{thm:duality-conical}
Let $\mu \in \mathcal{P}(\X)$ and assume that $F$ satisfies \eqref{eq:A'} and that there exists $\lambda>1$ such that \[
M:=\sup_{y \in \Y \cup\{0\}}\int F(x ,\lambda y)\,\mu(dx) <+\infty.
\]
If $F$ satisfies \eqref{eq:C'} or if $0$ does not belong to $\mathrm{co}(\Y)$, then for any probability measure $\nu \in \mathcal{P}(\Y)$, it holds
\begin{equation}\label{eq:duality-conical}
\mathcal{I}_c(\mu,\nu) = \sup_{\varphi \in \Phi( \Z) \cap L^1(\nu)}\left\{  \int Q_{F} \varphi (x)\,\mu(dx) - \int \varphi(y)\,\nu(dy)\right\},
\end{equation}
where $\Phi( \Z)$ is the set of all lower semicontinuous, convex positively $1$-homogenous functions $\varphi :  \Z \to \R\cup \{+\infty\}$ and where
\[
Q_{F} \varphi(x) = \inf_{z \in  \Z} \left\{ \varphi(z) + F(x,z)\right\},\qquad x\in \X.
\]
Moreover, one can restrict the supremum in \eqref{eq:duality-conical} to functions  $\varphi \in \Phi( \Z) \cap L^1(\nu)$ such that
\begin{equation}\label{eq:bdbelow}
\varphi(z) \geq \frac{\mathcal{I}_c(\mu,\nu)-M}{\lambda-1},\qquad \forall z\in \mathrm{co}(\Y).
\end{equation}
In addition, there exists a function $\bar{\varphi} \in \Phi( \Z)\cap L^1(\nu)$ satisfying \eqref{eq:bdbelow} and such that
\[
\mathcal{I}_c(\mu,\nu) =  \int Q_{F} \bar{\varphi} (x)\,\mu(dx) - \int \bar{\varphi}(y)\,\nu(dy).
\]
\end{thm}
\proof
Let $\nu \in \mathcal{P}(\Y)$. Without loss of generality, we assume that $\Y = \mathrm{co}\,\mathrm{Supp}(\nu)$. Note that this is always possible since the convex hull of a compact set is itself compact. To make the proof easier to read, we also assume that $\Y= \mathrm{co}\,\mathrm{Supp}(\nu)$ has a non empty interior. If this is not the case, then one can easily adapt the arguments below using the notions of relative interior and relative boundary of a convex set.
The proof is divided into three steps.\\
\textbf{Step 1.} In this step, we show that
\begin{equation}\label{eq:duality-facile}
\sup_{\varphi \in \Phi( \Z)\cap L^1(\nu)} \int Q_F\varphi\,d\mu - \int \varphi\,d\nu \leq \mathcal{I}_c(\mu,\nu).
\end{equation}
If $\varphi \in \Phi( \Z)$ is $\nu$ integrable and $q \in \mathcal{Q}(\mu,\nu)$, then using Jensen's inequality and the positive $1$-homogeneity of $\varphi$, one gets
\[
\int \varphi\,d\nu = \int \left(\int \varphi(y)\,q^x(dy)\right) \mu(dx) \geq \int \varphi \left(\int y\,q^x(dy)\right) \mu(dx).
\]
On the other hand
\begin{align*}
\int Q_F\varphi\,d\mu& \leq \int \varphi\left( \int y\,q^x(dy) \right)+F\left(x,\int y\,q^x(dy)\right)\,\mu(dx) \\
& \leq\int \varphi\,d\nu  + \int F\left(x,\int y\,q^x(dy)\right)\,\mu(dx).
\end{align*}
Thus optimizing over $\varphi$ and over $q$ gives \eqref{eq:duality-facile}.

\noindent \textbf{Step 2.} In this step, we assume that 
\[
\left(F \text{ satisfies Assumption \eqref{eq:C'}}\right)\qquad
\text{or}
\qquad \left(F \text{ satisfies \eqref{eq:B'} and } 0 \notin \mathrm{co}(\Y)\right)
\]
and we prove that the converse inequality holds true in \eqref{eq:duality-facile} and that the supremum can be restricted to functions satisfying \eqref{eq:bdbelow}.
Recall that according to \eqref{thm:duality2}, one has
\[
\mathcal{I}_c(\mu,\nu) = \sup_{f \in \mathcal{C}_b(\Y)} \int K_cf(x)\,\mu(dx) - \int f(y)\,\nu(dy).
\]
Observe that for all $f\in \mathcal{C}_b(\Y)$
\begin{align*}
K_c f(x)& = \inf_{m \in \mathcal{M}(\Y)}\left\{ \int f(y)\,m(dy)+ F\left(x,\int y\,m(dy)\right)\right\} = Q_F \tilde{f}(x),
\end{align*}
where
\[
\tilde{f}(z) = \inf \left\{ \int f(y)\,m(dy) : m\in \mathcal{M}(\Y),  \int y\,m(dy)= z\right\},\qquad z\in  \Z.
\]
If $\tilde{f}$ takes the value $-\infty$ at some $z_o \in \Z$, then $K_c(f)(x)=Q_F(\tilde{f})(x)= - \infty$ for all $x \in \X$ and so this function $f$ can be dropped from the dual formula above for $\mathcal{I}_c(\mu,\nu)$. 
Let us now assume that $\tilde{f}$ never takes the value $-\infty$ on $\Z$. By construction $\tilde{f} \leq f$ on $\Y$ and so, by positive $1$-homogeneity, $\tilde{f}$ takes finite values on $\Z$. This function $\tilde{f}$ is the greatest convex and positively $1$-homogenous function $\varphi :  \Z \to \R$ such that $\varphi \leq f$ on $\Y$. Denote by $\bar{f}:\Z \to \R$ the lower semicontinuous envelop of $\tilde{f}$, that is the greatest lower semicontinuous function $\varphi : \Z \to \R$ such that $\varphi \leq \tilde{f}$. The function $\bar{f}$ is still convex and positively $1$-homogenous.
Moreover, we claim that $Q_F(\bar{f}) = Q_F(\tilde{f})$. The inequality $Q_F(\bar{f}) \leq Q_F(\tilde{f})$ is clear. Let us show the other direction. On $\mathrm{int}(\Z)$, it holds $\bar{f} = \tilde{f}$, and by convexity and lower semicontinuity, for all $z \in \partial \Z$ and $x\in \X$, one gets
\[
\bar{f}(z) = \lim_{t \to 1^-} \tilde{f}((1-t)z +ta) \qquad \text{and}\qquad F(x,z) =  \lim_{t \to 1^-} F(x,(1-t)z +ta)
\]
where $a \in \mathrm{int}(Z)$ is some fixed point. Since $(1-t)z +ta \in \mathrm{int}(\Z)$ for $0\leq t<1$, we get that
\[
Q_{F}(\bar{f})(x) = \inf_{z \in \mathrm{int}(\Z)} \{\bar{f}(z) + F(x,z)\} =\inf_{z \in \mathrm{int}(\Z)} \{\tilde{f}(z) + F(x,z)\} \geq Q_F(\tilde{f})(x),
\]
which completes the proof of the claim.
Since $\bar{f}$ is convex, it is bounded from below by some affine map. Thus $\bar{f}$ is $\nu$-integrable. Also, since $\bar{f} \leq f$, one has
\begin{align*}
\int K_cf(x)\,\mu(dx) - \int f(y)\,\nu(dy) &= \int Q_F(\bar{f}(x))\,\mu(dx) - \int f(y)\,\nu(dy)\\& \leq \int Q_F\bar{f}(x)\,\mu(dx) - \int \bar{f}(y)\,\nu(dy).
\end{align*}
Optimizing yields that
\[
\mathcal{I}_c(\mu,\nu) \leq \sup_{\varphi \in \Phi( \Z)\cap L^1(\nu)} \int Q_F\varphi(x)\,\mu(dx) - \int \varphi(y)\,\nu(dy),
\]
which completes the proof of \eqref{eq:duality-conical}.

Now, let us show that the supremum can be restricted to functions $\varphi \in \Phi( \Z) \cap L^1(\nu)$ satisfying \eqref{eq:bdbelow}.
If $\varphi \in \Phi( \Z) \cap L^1(\nu)$, then being $\nu$ integrable, it takes at least one finite value on the support of $\nu$. Since $\varphi$ is also lower semicontinuous, it reaches its infimum on $\mathrm{co}(\Y)$ at some point $z_0 \in \mathrm{co}(\Y)$. By definition of $Q_F\varphi$, it holds
\[
Q_F\varphi(x) \leq \varphi( \lambda z_0)+F(x,\lambda z_0), \qquad \forall x\in \X.
\]
Therefore, $\varphi$ being positively homogenous, one obtains
\[
\int Q_F\varphi(x)\,\mu(dx) - \int \varphi(y)\,\nu(dy) \leq (\lambda-1) \varphi(z_0)+\int F(x,\lambda z_0)\,\mu(dx) \leq (\lambda-1) \varphi(z_0)+ M.
\]
Thus, if $\varphi(z_0) <  \frac{\mathcal{I}_c(\mu,\nu)-M}{\lambda - 1}$, then $\int Q_F\varphi(x)\,\mu(dx) - \int \varphi(y)\,\nu(dy) <\mathcal{I}_c(\mu,\nu)$ and so one can drop such functions $\varphi$ from the supremum in \eqref{eq:duality-conical}, which completes the proof.

\noindent \textbf{Step 3.} In this step, we get rid of the Assumption \eqref{eq:B'}. More precisely, we assume that 
\[
\left(F \text{ satisfies Assumption \eqref{eq:C'}}\right)\qquad
\text{or} 
\qquad \left(0 \notin \mathrm{co}(\Y)\right)
\]
and prove dual attainment.

For all $n\geq 1$, $x \in \X$ and $z\in \Z$ define
\[
F_n(x,z) = \left\{\begin{array}{ll}  F(x,z) & \text{if } F \text{ satisfies Assumption \eqref{eq:C'}}   \\
 F(x,z) + \frac{1}{n}\|z\|^2 & \text{if } F \text{ does not satisfy Assumption \eqref{eq:C'} and }0 \notin \mathrm{co}{\Y} \end{array}\right.,
\]
and $c_n(x,m) = F_n(x,\int y\,dm)$, $x\in \X$, $m\in \mathcal{M}(\Y).$
The result of Step 2 is applicable to these functions $F_n$, $n\geq 1.$ Let 
\[
M_n = \sup_{y \in \Y \cup\{0\}}\int F_n(x ,\lambda y)\,\mu(dx) <+\infty,
\]
and observe that $M_n \to M$ as $n\to \infty$.

For all $n\geq1$, $\mathcal{I}_c(\mu,\nu) \leq \mathcal{I}_{c_n}(\mu,\nu)$.

Let $(\varphi_n)_{n\geq 1}$ be a sequence in $\Phi( \Z) \cap L^1(\nu)$ satisfying \eqref{eq:bdbelow} (with $M_n$) and such that for all $n\geq1$
\[
\mathcal{I}_{c_n}(\mu,\nu) \leq \int Q_{F_n}\varphi_{n}\,d\mu-\int \varphi_{n}\,d\nu + \frac{1}{n}.
\]
Such a sequence exists thanks to Step 2.
For all $n\geq1$, one has
\begin{equation}\label{eq:ascoli1}
Q_{F_n}\varphi_n(x) \leq \varphi_n(0) + F_n(x,0) = F(x,0),\qquad \forall x\in \X.
\end{equation}
Therefore, using the integrability assumption on $\mu$, one gets that 
\[
\sup_{n\geq 1}\int Q_{F_n}\varphi_n(x)\,\mu(dx) \leq M <\infty.
\]
Since $\int Q_{F_n}\varphi_{n}\,d\mu-\int \varphi_{n}\,d\nu \geq \mathcal{I}_c(\mu,\nu)-\frac{1}{n}$, this implies that $\sup_{n\geq1}\int \varphi_{n}\,d\nu<+\infty$.

Let us show that $(\varphi_{n})_{n\geq1}$ admits a converging subsequence. Define $\widetilde{\varphi}_n = \varphi_n + \frac{M_n-\mathcal{I}_c(\mu,\nu)}{\lambda-1}$, $n\geq1.$ For all $n\geq1$, $p\geq0$, $y_1,\ldots, y_p\in \mathrm{Supp}(\nu)$ and $\lambda_1,\ldots, \lambda_p \geq0$ such that $\sum_{i=1}^p \lambda_i=1$, it follows from Jensen's inequality and from $\widetilde{\varphi}_{n}\geq0$ on $\Y$ that
\begin{align*}
\widetilde{\varphi}_{n} \left(\sum_{i=1}^p \lambda_i \frac{\int_{B(y_i,\epsilon)} z \,\nu(dz)}{\nu(B(y_i,\epsilon))}\right) &\leq \sum_{i=1}^p \lambda_i \widetilde{\varphi}_{n} \left(\frac{\int_{B(y_i,\epsilon)} z \,\nu(dz)}{\nu(B(y_i,\epsilon))}\right)\\
 &\leq \sum_{i=1}^p \lambda_i \frac{1}{\nu(B(y_i,\epsilon))} \int_{B(y_i,\epsilon)} \widetilde{\varphi}_{n}(z)\,\nu(dz)\\
& \leq \left(\sum_{i=1}^p \lambda_i \frac{1}{\nu(B(y_i,\epsilon))}\right) \int \widetilde{\varphi}_{n}(z)\,\nu(dz),
\end{align*}
with $B(y,\epsilon)$ denoting the open ball of center $y$ and radius $\epsilon>0.$
Since $\sup_{n\geq1}\int \widetilde{\varphi}_{n}\,d\nu<+\infty$, it holds
\[
\sup_{n\geq 1} \widetilde{\varphi}_{n}(u)<+\infty
\]
for all $u$ belonging to set $C = \mathrm{co}\left\{\frac{\int_{B(y,\epsilon)} z \,\nu(dz)}{\nu(B(y,\epsilon))} : y \in \mathrm{Supp}(\nu), \epsilon>0 \right\}.$
Let us show that $C$ is dense in $\Y$. Take $y \in \Y$;  since $\Y = \mathrm{co}\, \mathrm{Supp}(\nu)$, there exists $y_1,\ldots,y_p \in \mathrm{Supp}(\nu)$ and $\lambda_1,\ldots,\lambda_p \geq 0$ such that $\sum_{i=1}^p \lambda_i=1$ and $y = \sum_{i=1}^p \lambda_i y_i$. For all $\epsilon>0$, define $y_{\epsilon} = \sum_{i=1}^p \lambda_i \frac{\int_{B(y_i,\epsilon)} z \,\nu(dz)}{\nu(B(y_i,\epsilon))} \in C$. Then, $y_\epsilon \to y$ when $\epsilon\to0$, which proves the claim.

According to \cite[Theorem 10.9 page 91]{Roc70}, one can extract from $(\widetilde{\varphi}_{n})_{n\geq1}$ a subsequence (we still denote it by $(\widetilde{\varphi}_{n})_{n\geq1}$ to avoid overloading notations) converging pointwise on $\mathrm{int} (\Y)$ to some convex function. Of course, the sequence $(\varphi_n)_{n\geq1}$ also converges pointwise on $\mathrm{int} (\Y)$. Since $ \Z = \R_+\Y$ and $\mathrm{int} ( \Z) = \R_+^* \mathrm{int} (\Y)$, the positive homogeneity of the functions $\varphi_n$ implies that $\varphi_n$ converges pointwise on $\mathrm{int} ( \Z)$. Set $\varphi(x) = \lim_{n\to \infty} \varphi_n (x)$, for all $x \in \mathrm{int} ( \Z)$. Extend $\varphi$ by setting $\varphi(a) = \liminf_{z\to a, z\in \mathrm{int}( \Z)} \varphi(z)$ whenever $a \in \partial  \Z$, so that $\varphi$ is lower semicontinuous on $ \Z$ (and still convex and positively homogenous). Obviously $\varphi$ satisfies \eqref{eq:bdbelow}.
If $a \in \partial  \Z$ and $z \in \mathrm{int} ( \Z)$, then for all $t \in (0,1)$ the point $(1-t) a + t z$ belongs to  $\mathrm{int} ( \Z)$. Therefore, letting $n \to \infty$ in the inequality
\[
\varphi_{n}((1-t) a + t z) \leq (1-t)\varphi_{n}(a)+t\varphi_{n}(z)
\]
one gets that
\[
\varphi((1-t) a + t z) \leq (1-t)\liminf_{n\to+\infty}\varphi_{n}(a)+t\varphi(z)
\]
and letting $t \to 0$ gives then
\[
\varphi(a)\leq \liminf_{n\to+\infty}\varphi_{n}(a).
\]
Then,
\begin{align*}
\mathcal{I}_c(\mu,\nu) &\leq  \limsup_{n\to +\infty}\left(\int Q_{F_n}\varphi_{n}\,d\mu-\int \varphi_{n}\,d\nu\right)\\
& \leq  \limsup_{n\to +\infty}\int Q_{F_n}\varphi_{n}\,d\mu- \liminf_{n\to +\infty}\int \varphi_{n}\,d\nu\\
& \leq  \int \limsup_{n\to +\infty} Q_{F_n}\varphi_{n}\,d\mu- \int \liminf_{n\to +\infty} \varphi_{n}\,d\nu\\
& \leq \int \limsup_{n\to +\infty} Q_{F_n}\varphi_n\,d\mu - \int \varphi \,d\nu,
\end{align*}
where the third inequality follows by Fatou lemma (note that, for all $n\geq1$, $\varphi_n\geq \inf_{k\geq 1}\frac{\mathcal{I}_c(\mu,\nu)-M_k}{\lambda-1}$ on $\Y$ and $Q_{F_n}\varphi_n \leq F(\,\cdot\,,0)$ which is $\mu$-integrable), and the last inequality because $\liminf_{n\to+\infty}\varphi_{n} \geq \varphi$ on $ \Z$.
For all $n\geq 0$, one gets that
\[
\varphi_{n}(z) \geq Q_F\varphi_{n}(x) - F_n(x,z),\qquad \forall x\in \X, \forall z \in  \Z.
\]
Therefore, letting $n \to \infty$, one also gets that
\[
\varphi(z) \geq \limsup_{n\to +\infty} Q_{F_n}\varphi_n(x)-F(x,z),\qquad \forall x \in \X,\forall z\in \mathrm{int}( \Z).
\]
So, for all $x\in \X$,
\[
 \limsup_{n\to +\infty} Q_{F_n}\varphi_n(x) \leq \inf_{z \in \mathrm{int}( \Z)} \{\varphi(z) + F(x,z)\} = \inf_{z \in \Z} \{\varphi(z) + F(x,z)\},
\]
where the last equality was proved in Step 2.
In conclusion, we have shown the existence of a function $\varphi\in \Phi( \Z)\cap L^1(\nu)$ such that
\[
\mathcal{I}_c( \mu,\nu) \leq \int Q_F\varphi\,d\mu - \int \varphi \,d\nu.
\]
Since the converse inequality is always true (according to Step 1), this completes the proof.
\endproof

\subsection{A new variant of Strassen's theorem}\label{sec:strassen}
Recall that if $\mu,\nu$ are two probability measures on $\R^d$, each having a finite first moment, $\mu$ is said to be dominated by $\nu$ in the convex order, if
\[
\int \varphi\,d\mu \leq \int \varphi\,d\nu,
\]
for all convex function $\varphi:\R^d \to \R.$ In this case, we denote this relation by $\mu \leq_c\nu$.
According to a well-known result due to Strassen \cite{Str65}, $\mu \leq_c\nu$ if and only if there exists a martingale coupling with marginals $\mu$ and $\nu$, that is a couple $(U,V)$ of random vectors with $U\sim \mu$, $V \sim \nu$ and $\E[V\mid U] = U$ a.s.

Transport problems with conical cost functions introduced above are naturally related to the following variant of the convex order. If $\mu,\nu$ are two probability measures with a finite moment of order $1$, we will say that $\mu$ is dominated by $\nu$ for the positively $1$-homogenous convex order if for all $\varphi : \R^d \to \R$ convex and positively $1$-homogenous, one has $\int \varphi\,d\mu \leq \int \varphi\,d\nu$. We will use the notation $\mu \leq _{phc} \nu$ to denote this order.

The following result generalizes Strassen's theorem to this restricted convex order. Note that if $\nu$ is compactly supported, then $\mu \leq_{phc} \nu$ does not imply that $\mu$ is also compactly supported.
\begin{thm}\label{thm:convorder}Let $\mu,\nu$ be two probability measures on $\R^d$ with $\nu$ being compactly supported. 
\begin{enumerate}
\item[(a)] If $\mu$ is compactly supported then the following are equivalent:
\begin{itemize}
\item[$(i)$] $\mu \leq _{phc} \nu$,
\item[$(ii)$] There exists a nonnegative kernel $q$ such that $\mu q = \nu$ and
\[
\int y\,q^x(dy) = x
\] for $\mu$ almost every $x$,
\item[$(iii)$] There exists a probability measure $\eta$ absolutely continuous with respect to $\mu$ with density denoted by $N$ and a couple of random vectors $(U,V)$ with $U \sim \eta$, $V \sim \nu$ such that
\begin{equation}\label{eq:martingalemod}
N(U)\E[V\mid U] = U \text{ a.s.}
\end{equation}
\end{itemize}
\item[(b)]The same conclusion holds if $\mu$ has a finite moment of order $1$ and the convex hull of the support of $\nu$ does not contain $0$. 
\end{enumerate}
\end{thm}
Also note that \eqref{eq:martingalemod} means that $(U,N(U)V)$ is a martingale. 

\begin{rem}\label{rem:dim1}
In dimension $1$, the conclusion of Theorem \ref{thm:convorder} is essentially trivial. Indeed, it is easy to see that $\mu\leq_{phc} \nu$ if and only if $\int x\,d\mu = \int x\,d\nu$, $\int [x]_+\,d\mu \leq  \int [x]_+\,d\nu$ and $\int [x]_-\,d\mu \leq  \int [x]_-\,d\nu$. By assumption, the convex hull of the support of $\nu$ does not contain $0$, so it is contained either in $(0,\infty)$ or in $(-\infty,0)$. Let us assume without loss of generality that the support of $\nu$ is contained in $(0,\infty)$. Then, $\int [x]_-\,d\mu \leq  \int [x]_-\,d\nu=0$ holds and so the support of $\mu$ is also contained in $(0,\infty)$.  Consider the nonnegative function $N(x) = \frac{x}{\int x\,d\mu} \mathbf{1}_{(0,\infty)}(x)$, which satisfies $\int N(x)\,\mu(dx)=1$, and define $\eta = N\,\mu$. Let $U,V$ be two independent random variables such that $U \sim \eta$ and $V \sim \nu$, then
\[
\E[V\mid U] = \E[V] = \int x\,d\mu = \frac{U}{N(U)} \text{a.s}
\]
holds. In higher dimension, it is not clear whether such a simple and explicit construction is available.
\end{rem}

\proof[Proof of Theorem \ref{thm:convorder}]
First, $(ii)$ and $(iii)$ are equivalent.
Suppose that $\mu$ has a finite moment of order $1$ and $\nu$ is compactly supported, and let us show that $(ii)$ implies $(i).$ Let $\varphi$ be some convex positively $1$-homogenous function; according to Jensen's inequality and positive $1$-homogeneity, one has
\begin{align*}
\int \varphi(x)\,\mu(dx) &= \int \varphi\left(\int y\,q^x(dy)\right)\,\mu(dx) \leq \iint \varphi(y)\,q^x(dy)\,\mu(dx) =  \int \varphi(y)\,\nu(dy).
\end{align*}
Now let us show that $(i)$ implies $(ii)$ when $\mu$ is compactly supported. 
Let us denote by $\X$ and $\Y$ the compact supports of $\mu$ and $\nu$ and consider the cost function $c : \X \times \mathcal{M}(\Y) \to \R_+$
\[
c(x,m) = F\left(x,\int y\,m(dy)\right)= \left\|x-\int y\,m(dy)\right\|,\qquad \forall x\in \X,\forall m\in \mathcal{M}(\Y),
\]
with $F(x,z) = \|x-z\|$, $x,z\in \R^d.$
First, $F$ satisfies Assumption \eqref{eq:A'} and $\eqref{eq:C'}$. Therefore, according to Theorem \ref{thm:duality-conical}, the following duality formula holds
\[
\mathcal{I}_c(\mu,\nu)= \sup_{\varphi \in \Phi( \Z) \cap L^1(\nu)}\left\{  \int Q\varphi (x)\,\mu(dx) - \int \varphi(y)\,\nu(dy)\right\},
\]
with $Q\varphi(x) = \inf_{z \in \Z} \{\varphi(z) +\left\|x-z\right\|\}$, $x\in \R^d$. The supremum can be restricted to $\varphi$ that are bounded from below by some constant $\kappa \in \R$.  For such functions $\varphi$, $Q\varphi$ is finite-valued, convex, and positively $1$-homogenous on $\R^d.$ Thus $\int Q\varphi\,d\mu \leq \int Q\varphi\,d\nu.$ Since $Q\varphi \leq \varphi$ on $\Z$, we conclude that $\int Q\varphi\,d\mu \leq \int \varphi\,d\nu$. Therefore, $\mathcal{I}_c(\mu,\nu)=0$. On the other hand, since $F$ satisfies Assumption \eqref{eq:C'}, it follows from Theorem \ref{thm:weaksol} and Lemma \ref{lem:recovery} that there exists some $\pi \in \Pi(\mathrm{Supp}(\mu),\nu)$ such that
\[
\bar{I}_c^\mu[\pi] = \bar{\mathcal{I}}_c(\mu,\nu) = \mathcal{I}_c(\mu,\nu)=0.
\]
Since $c'_\infty(x,m) =  \left\|\int y\,m(dy)\right\|$, $m\in \mathcal{M}(\Y)$, one thus gets
\begin{equation}\label{eq:0}
\int  \left\|x-\frac{d\eta^{ac}}{d\mu}(x)\int y\,\pi^x(dy)\right\|\,\mu(dx) + \int  \left\|\int y\,\pi^x(dy)\right\|\,\eta^s(dx)=0,
\end{equation}
where $\eta=\eta^{ac}+\eta^s$ denotes the first marginal of $\pi$. Let us define $q_1^x(dy) = \frac{d\eta^{ac}}{d\mu}(x)\pi^x(dy)$, $x\in \X$,  $\nu_1 = \mu q_1$ and $\nu_2 = \nu - \nu_1$. It follows from \eqref{eq:0} that $\int y\,q_1^x(dy)=x$ for $\mu$ almost all $x$. Moreover, $\int y\,\pi^x(dy)=0$ for $\eta^s$ almost all $x$. Therefore,
\[
\int y\,\nu_2(dy) = \int \left(\int y \pi^x(dy)\right) \eta(dx) -  \int \left(\int y \pi^x(dy)\right) \eta^{ac}(dx) =  \int \left(\int y \pi^x(dy)\right) \eta^s(dx)=0.
\]
The non-negative kernel $q$ defined by $q^x(dy) = q_1^x(dy)+\nu_2(dy)$, $x \in \X$, is such that $q \in \mathcal{Q}(\mu,\nu)$ and $\int y\,q^x(dy) = x$ for $\mu$ almost all $x \in \X$, which completes the proof of part (a). 

Finally, let us show that $(i)$ implies $(ii)$ when $\mu$ has a finite first moment and the convex hull of the support of $\nu$ does not contain $0$.
Let us construct a sequence of compactly supported probability measures $(\mu_n)_{n\geq 1}$ converging to $\mu$ in the weak topology and such that $\mu_n \leq_c \mu$ for all $n\geq1$. One can, for instance, obtain such a sequence as follows. Consider $C_n = [-n,n]^d$ and write $\R^d \setminus C_n = \bigcup_{1\leq k\leq K_n} D_{n,k}$, where $(D_{n,k})_{1\leq k\leq K_n}$ are disjoints convex subsets of $\R^d$. Then define
\[
\mu_n(dx)= \mathbf{1}_{C_n}(x)\,\mu(dx) + \sum_{k=1}^{K_n} \mu(D_{n,k})\delta_{x_{n,k}}(dx)
\]
where $x_{n,k} = \frac{1}{\mu(D_{n,k})}\int_{D_{n,k}} x \,\mu(dx)$ if $\mu(D_{n,k}) >0$ and any point in $D_{n,k}$ otherwise.
If $f:\R^d \to \R$ is a convex function, then it follows from Jensen's inequality that
\[
\int f(x)\,\mu_n(dx) = \int_C f(x)\,\mu(dx) +  \sum_{k=1}^{K_n} \mu(D_{n,k})f(x_{n,k}) \leq \int_C f(x)\,\mu(dx) +  \sum_{k=1}^{K_n} \int_{D_{n,k}} f(x)\,\mu(dx) = \int f(x)\,\mu(dx)
\]
and so $\mu_n \leq_c \mu$.
Also, for any bounded continuous function $f : \R^d \to \R$, one has that
\[
\left|\int f\,d\mu_n - \int f\,d\mu\right| = \left|\int_{\R^d \setminus C_n} f\,d\mu_n - \int_{\R^d \setminus C_n} f\,d\mu\right|  \leq 2 \|f\|_\infty \left(1 - \mu (C_n)\right) \to 0,
\]
as $n\to \infty$, and so $(\mu_n)_{n\geq1}$ converges to $\mu$ in the weak topology.

Since $\mu_n \leq_c \mu$ and $\mu \leq_{phc} \nu$, $\mu_n \leq_{phc} \nu$. By construction, $\mu_n$ has a compact support, and so, according to part (a), there exists a nonnegative kernel $q_n$ such that $\mu_nq_n = \nu$ and $\int y\,q_n^x(dy) = x$ for $\mu_n$ almost all $x \in \R^d$. For all $n\geq 1$, write for all $x\in \R^d$, $q_n^x = N_n(x)p_n^x$ where $p_n^x$ is a probability measure and denote $\eta_n(dx) = N_n(x)\mu_n(dx)$ and $\pi_n(dxdy) = \eta_n(dx)p_n^x(dy).$ Let us show that the sequence $(\eta_n)_{n\geq1}$ is tight. By assumption, $\int \|x\|\,\mu(dx) <+\infty$; thus by the de la Vall\'ee Poussin theorem, there exists some non-decreasing convex function $\alpha : \R_+ \to \R_+$ such that $\alpha(0)=0$, $\alpha(x)/x \to +\infty$ as $x \to +\infty$ and $\int \alpha(\|x\|)\,\mu(dx)<+\infty.$ The function $x\mapsto \alpha(\|x\|)$ being convex, we thus get that
\[
\sup_{n\geq 1}\int \alpha(\|x\|)\,\mu_n(dx) \leq \int \alpha(\|x\|)\,\mu(dx) :=M
\]
and, since $\int y\,q_n^x(dy) = x$ for $\mu_n$ almost all $x$,
\[
\sup_{n\geq 1}\int \alpha\left(\left\|\int y\,q_n^x(dy)\right\|\right)\,\mu_n(dx) \leq M.
\]
But,
\[
\int \alpha\left(\left\|\int y\,q_n^x(dy)\right\|\right)\,\mu_n(dx) = \int \alpha\left(N_n(x) B_n(x)\right)\,\mu_n(dx),
\]
where $B_n(x) =\left\|\int y\,p_n^x(dy)\right\|$. Since $\int y \,p_n^x(dy)$ belongs to the convex hull of the support of $\nu$ which is a compact convex set not containing zero, there exists some $b>0$ independent of $n$ such that $B_n(x) \geq b$ for all $n\geq 1$ and $x \in \R^d$. Since $\alpha$ is non-decreasing, then
\[
\sup_{n\geq 1}\int \alpha\left(bN_n(x)\right)\,\mu_n(dx) \leq M.
\]
Set $\alpha^*(t) = \sup_{s\geq 0} \{st - \alpha(s)\}$, $t \geq 0$, and note that $\alpha^*$ is non-decreasing, finite valued and vanishes at $0$. If $f : \R^d \to \R_+$ is a nonnegative function, then using Young's inequality $st \leq \alpha(s) + \alpha^*(t)$, for all $s,t \geq 0$ it is easily seen that for all $u>0$, it holds
\[
\int uf\,d\eta_n = \int ufN_n\,d\mu_n \leq \frac{1}{b} \int \alpha^*(uf)\,d\mu_n +\frac{1}{b} \int \alpha (b N_n)\,d\mu_n
\]
and so
\begin{equation}\label{eq:borne-eta_n}
\int f\,d\eta_n \leq \frac{1}{bu} \int \alpha^*(uf)\,d\mu_n + \frac{M}{bu}.
\end{equation}
In particular, if $f = \mathbf{1}_{A}$ with $A$ a measurable set, then
\[
\eta_n(A) \leq  \frac{\max(M ; 1)}{b} \psi(\mu_n(A)),
\]
where
\[
\psi(t) = \inf_{u>0} \left\{\frac{\alpha^*(u)}{u}t + \frac{1}{u} \right\},\qquad t\geq0.
\]
Since the sequence $(\mu_n)_{n\geq 1}$ converges to $\mu$, it is tight and so for all $\varepsilon>0$, there exists a compact set $K_\varepsilon$ such that $\sup_{n\geq 1}\mu_n(\R^d \setminus K_\varepsilon) \leq \varepsilon$. Choosing $A = \R^d\setminus K_\varepsilon$ one thus sees that
\[
\sup_{n \geq0}\eta_n(\R^d \setminus K_\varepsilon) \leq  \frac{\max(M ; 1)}{b} \psi(\varepsilon).
\]
First, note that $\psi\left(\frac{1}{\alpha(s)}\right) = \frac{s}{\alpha(s)}$, as soon as $\alpha(s)>0$. This implies in particular that $\psi(t) \to 0$ as $t\to 0$. So the bound above shows that the sequence $(\eta_n)_{n\geq 1}$ is tight. Therefore, according to Prokhorov's theorem, one can extract from the sequence $(\eta_n)_{n\geq 1}$ a subsequence converging to some probability measure $\eta$ on $\R^d$. For notational convenience this subsequence will still be denoted by $(\eta_n)_{n\geq 1}$. Letting $n \to \infty$ in \eqref{eq:borne-eta_n}, one sees that for all bounded continuous and nonnegative function $f$, the inequality
\[
\int f\,d\eta \leq \frac{1}{bu} \int \alpha^*(uf)\,d\mu + \frac{M}{bu}
\]
holds. If $A$ is a compact subset of $\R^d$, then considering the sequence $f_k(x) = \left[1 - kd(x,A)\right]_+$, $k\geq0$, of bounded continuous nonnegative functions which converges monotonically to $\mathbf{1}_A$, one sees that
\[
\eta(A) \leq  \frac{1}{bu} \alpha^*(u)\mu(A) + \frac{M}{bu},\qquad \forall u>0.
\]
Since both $\eta$ and $\mu$ are inner regular, this inequality is actually true for any measurable set $A$ of $\R^d.$ In particular, if $\mu(A)=0$ then it follows that $\eta(A)=0$ and so $\eta$ is absolutely continuous with respect to $\mu$. Since $\eta_n$ converges, the sequence $\pi_n \in \Pi(\eta_n,\nu)$ is also tight. One can thus assume without loss of generality that it converges to some $\pi \in \Pi(\eta,\nu)$ in the weak topology. If $u$ is some compactly supported function on $\R^d$, then
\[
\iint u(x)y \,\pi_n(dxdy) = \int u(x) \int y\,q_n^x(dy)\mu_n(dx) = \int u(x) x \,\mu_n(dx)
\]
holds, and letting $n\to \infty$ gives
\[
\iint u(x)y \,\pi(dxdy) = \int u(x) x \,\mu(dx)
\]
which reads, writing $\pi(dxdy) = N(x)\mu(dx)p^x(dy)$,
\[
\int u(x)\left(N(x) \int y\,p^x(dy)-x\right) \mu(dx) = 0.
\]
Since this holds for all compactly supported continuous functions $u$, one concludes that for $\mu$ almost all $x \in \R^d$
\[
\int y\,q^x(dy)=x
\]
with $q^x = N(x)p^x$, which completes the proof of part (b).
\endproof

\subsection{Study of a particular class of nonpositive conical transport problems}
In this Section, we consider a \emph{nonpositive} cost function
\[
F: \X\times  \Z \to \R_-
\]
such that $F(x,\lambda z)/\lambda \to 0$ as $\lambda \to +\infty$ for all $x \in \X$, $z \in \Z \setminus\{0\}$. The following result shows in particular that under mild additional assumptions on $F$ any weak solution is strong.

\begin{thm}\label{thm:duality-eco}
Assume that
\begin{itemize}
\item $0$ does not belong to $\mathrm{co}(\Y)$,
\item $F$ is nonpositive, satisfies Assumption \eqref{eq:A'}, and is continuous on $\X\times  \Z$,
\item For all $x\in \X$ and $y\in \Y$, 
\[
\frac{F(x,\lambda y)}{\lambda} \to 0
\]
holds as $\lambda \to +\infty$. 
\item
\[
\sup_{x\in \X,y\in \Y} F(x,\lambda y) \to -\infty
\]
holds as $\lambda \to +\infty$.
\end{itemize}
Then,
\begin{equation}\label{eq:duality-eco}
\mathcal{I}_c(\mu,\nu) = \sup_{\varphi \in \Phi^+( \Z) \cap L^1(\nu)}\left\{  \int Q_{F} \varphi (x)\,\mu(dx) - \int \varphi(y)\,\nu(dy)\right\},\qquad \forall \nu \in \mathcal{P}(\Y)
\end{equation}
where $\Phi^+( \Z)$ is the set of all nonnegative, lower semicontinuous, convex, and positively $1$-homogenous functions $\varphi :  \Z \to \R_+\cup \{+\infty\}$ and where
\[
Q_{F} \varphi(x) = \inf_{z \in  \Z} \left\{ \varphi(z) + F(x,z)\right\},\qquad x\in \X.
\]
Furthermore, the supremum in \eqref{eq:duality-eco} is attained at some function $\bar{\varphi} \in \Phi^+( \Z) \cap L^1(\nu)$, which is positive on $ \Z \setminus \{0\}$.
Finally, any weak solution for the transport problem between $\mu$ and $\nu$ is a strong solution.
\end{thm}
\begin{rem}Note that the Assumption $\sup_{x\in \X,y\in \Y} F(x,\lambda y) \to -\infty$ as $\lambda \to \infty$ actually implies that $0 \notin \mathrm{co}(\Y)$.
\end{rem}
The assumptions of Theorem \ref{thm:duality-eco} are for instance satisfied by $F:\X \times \R_+^d \to (-\infty,0]$ of the following form
\begin{equation}\label{eq:exampleCK}
F(x,z) = - \|A(x) z\|_\sigma^\eta,\qquad x\in \X, z\in  \R_+^d,
\end{equation}
where $0<\eta<1$,  $A : \X \to M_{>0}(\R^d)$ a continuous function taking values in the space $M_{>0}(\R^d)$ of $d\times d$ matrices with positive entries, and for $0<\sigma\leq 1$, the $\sigma$-``norm'' is defined by
\[
\|z\|_\sigma = \left(\sum_{i=1}^d |z_i|^\sigma\right)^{1/\sigma},\qquad z\in \R^d.
\]
Since $\|\,\cdot\,\|_\sigma$ satisfies the following reverse triangle inequality on $\R_+^d$:
\[
\|z_1+z_2\|_\sigma \geq \|z_1\|_\sigma + \|z_2\|_\sigma,\qquad \forall z_1,z_2 \in \R_+^d.
\]
This implies that the function $F$ defined by \eqref{eq:exampleCK} is convex with respect to its second variable. Since $0<\eta<1$, for every $x \in \X$ and $z \in \R_+^d$, $F(x,\lambda z)/\lambda \to 0$ as $\lambda \to +\infty$. Finally, if $\Y$ is a compact subset included in $(0,\infty)^d$, then $\sup_{x\in \X, y\in \Y} F(x,\lambda y) \to - \infty$ as $\lambda \to +\infty$.
\begin{rem}
Cost functions of the form \eqref{eq:exampleCK} are considered in \cite{CK21} to represent minus the output of a firm $x$ when it hires a worker of type $z$. In this context, the variable $\varphi$ appearing in the dual formulation of $\mathcal{I}_c(\mu,\nu)$ corresponds to a wage function and is thus naturally nonnegative.
\end{rem}

\proof
Applying Theorem \ref{thm:duality-conical} to the cost function $c$ yields
\[
\mathcal{I}_c(\mu,\nu) = \sup_{\varphi \in \Phi( \Z) \cap L^1(\nu)}\left\{  \int Q_{F} \varphi (x)\,\mu(dx) - \int \varphi(y)\,\nu(dy)\right\},\qquad \forall \nu \in \mathcal{P}(\Y).
\]
Let us show that the supremum can be restricted to nonnegative functions. Indeed, let $y_0 \in \mathrm{co}(\Y)$ be such that $\varphi (y_0) = \inf_{\mathrm{co}(\Y)} \varphi$ and assume that $\varphi(y_0)<0.$ Since $F \leq0$, for all $\lambda \geq0$ one has
\[
Q_{F} \varphi(x) \leq  \inf_{\lambda >0} \varphi (\lambda y_0)= \inf_{\lambda >0} \lambda \varphi ( y_0) = -\infty.
\]
Therefore, such functions $\varphi$ can be dropped from the supremum.
According to Theorem \ref{thm:duality-conical}, we also know that the supremum in \eqref{eq:duality-eco} is reached at some $\bar{\varphi}\in \Phi^+( \Z) \cap L^1(\nu)$. Let us show that this function $\bar{\varphi}$ is positive over $ \Z \setminus \{0\}$. Consider again $y_0 \in \mathrm{co}(\Y)$ such that $\bar{\varphi} (y_0) = \inf_{\mathrm{co}(\Y)} \varphi$ and set $a = \bar{\varphi} (y_0) \geq0$.
Define, for all $u\geq0$,
\[
\psi(u) = \inf_{\lambda >0} \left\{\lambda u + \sup_{x\in \X,y\in \Y} F(x,\lambda y)\right\}.
\]
Observe that
\[
Q_{F} \bar{\varphi}(x) \leq  \inf_{\lambda >0} \left\{ \bar{\varphi} (\lambda y_0) + F(x,\lambda y_0) \right\} =  \inf_{\lambda >0} \left\{ \lambda a + F(x, \lambda y_0) \right\} \leq \psi(a),
\]
where we used that, since $y_0 \in  \mathrm{co}(\Y)$, we have $F(x,\lambda y_0) \leq \sup_{y\in \Y} F(x,\lambda y).$
Therefore,
\[
-\infty <\mathcal{I}_c(\mu,\nu) = \int Q_{F} \bar{\varphi} (x)\,\mu(dx) - \int \bar{\varphi}(y)\,\nu(dy) \leq \psi(a)-a \leq  \psi(a).
\]
By assumption $\psi(0)=-\infty$, so $a>0$.
Therefore $\bar{\varphi}$ is positive on $\mathrm{co}(\Y)$. Since $ \Z = \R_+ \mathrm{co}(\Y)$, we conclude that $\bar{\varphi}$ is positive on $ \Z \setminus \{0\}$.

We now show that the transport problem admits only strong solutions. According to Theorem \ref{thm:weaksol} (applied to the cost function $c$ which satisfies condition \eqref{eq:C}), we know there exists $\bar{\pi} \in \Pi(\mathrm{Supp}(\mu),\nu)$ such that
\begin{equation}\label{eq:duality-eco-proof}
\mathcal{I}_c(\mu,\nu) = \overline{\mathcal{I}}_c(\mu,\nu) = \bar{I}_c^\mu[\bar{\pi}],
\end{equation}
where, for all $\pi \in \Pi(\mathrm{Supp}(\mu),\nu)$,
\[
\bar{I}_c^\mu[\pi] = \int F\left(x,\frac{d\pi^{ac}_1}{d\mu}(x)\int y\,\pi^x(dy)\right)\,\mu(dx)
\]
(because, by assumption, $F'_\infty(x,z)=0$ for all $x,z$).
Denoting by $S(x) = \frac{d\bar{\pi}^{ac}_1}{d\mu}(x)\int y\,\bar{\pi}^x(dy)$, one has
\begin{align*}
\mathcal{I}_c(\mu,\nu) & =  \int Q_F\bar{\varphi}(x)\,\mu(dx) - \int \bar{\varphi}(y)\,\nu(dy)\\
& \leq  \int \bar{\varphi}(S(x)) +F(x,S(x))\,\mu(dx) - \int \bar{\varphi}(y)\,\nu(dy)\\
&\stackrel{(i)}{\leq} \int \left(\int \bar{\varphi}(y)\,\bar{\pi}^x(dy) \right) \bar{\pi}_1^{ac}(dx)+\int F(x,S(x))\,\mu(dx)- \int \bar{\varphi}(y)\,\nu(dy)\\
& \stackrel{(ii)}{=} \bar{I}_c^\mu[\bar{\pi}] - \int\left(\int \bar{\varphi}(y)\,\bar{\pi}^x(dy) \right) \bar{\pi}_1^{s}(dx)\\
& \stackrel{(iii)}{=} \mathcal{I}_c(\mu,\nu) - \int\left(\int \bar{\varphi}(y)\,\bar{\pi}^x(dy) \right) \bar{\pi}_1^{s}(dx),
\end{align*}
where
\begin{itemize}
\item $(i)$ comes from the positive $1$-homogeneity of $\bar{\varphi}$ and Jensen's inequality,
\item $(ii)$ comes from the definition of $ \bar{I}_c^\mu[\bar{\pi}] $ and the fact that
\[
\int \bar{\varphi}(y)\,\nu(dy) = \int\left(\int \bar{\varphi}(y)\,\bar{\pi}^x(dy) \right) \bar{\pi}_1(dx) = \int\left(\int \bar{\varphi}(y)\,\bar{\pi}^x(dy) \right) \bar{\pi}_1^{ac}(dx)+\int\left(\int \bar{\varphi}(y)\,\bar{\pi}^x(dy) \right) \bar{\pi}_1^{s}(dx),
\]
\item and $(iii)$ comes from \eqref{eq:duality-eco-proof}.
\end{itemize}
We conclude that
\[
 \int\left(\int \bar{\varphi}(y)\,\bar{\pi}^x(dy) \right) \bar{\pi}_1^{s}(dx) \leq 0.
\]
Because $\int \bar{\varphi}(y)\,\bar{\pi}^x(dy) \geq a>0$, $\bar{\pi}_1^{s}=0$ (no singular part) is the only possible option.
Therefore, $\bar{\pi} \in \Pi(\ll \mu,\nu)$ and $\bar{I}_c^\mu[\bar{\pi}] = I_c^\mu[\bar{\pi}] = \mathcal{I}_c(\mu,\nu)$
and so $\bar{\pi}$ is a strong solution.
\endproof

\subsection{Structure of solutions for conical cost functions}
In this subsection, we assume that $F$ is a function satisfying Assumption \eqref{eq:A'} and that $c$ is the conical cost function associated to $F$ and defined by \eqref{eq:conical}.

The following result shows that the transport cost $\mathcal{I}_c(\mu,\nu)$ is a shortest ``transport distance'' between $\mu$ and the set of probability measures dominated by $\nu$ in the order $\leq_{phc}$.
\begin{thm}\label{thm:primal-structure}
Let $\mu \in \mathcal{P}(\X)$ and $\nu \in \mathcal{P}(\Y)$ be such that $\mathcal{I}_c(\mu,\nu)<+\infty$, and assume that the convex hull of the support of $\nu$ does not contain $0$.
Then, the following identity holds
\begin{equation}\label{eq:Ic-cvxorder}
\mathcal{I}_c(\mu,\nu)= \inf_{\gamma \leq_{phc} \nu} \mathcal{T}_{F}(\mu,\gamma),
\end{equation}
where and $\mathcal{T}_F$ denotes the classical transport cost associated to the cost function $F$:
\[
\mathcal{T}_F(\mu,\gamma) = \inf_{\pi \in \Pi(\mu,\gamma)} \iint F(x,z)\,\pi(dxdz),\qquad \forall \mu \in \mathcal{P}(\X),\forall \gamma \in \mathcal{P}( \Z).
\]
Moreover, suppose that $\bar{q}$ is a nonnegative kernel solution to the transport problem \eqref{eq:TP}, consider the map $\bar{S}$ defined by
\[
\bar{S}(x) = \int y\,\bar{q}^x(dy),\qquad x\in \X,
\]
and denote by $\bar{\nu}$ the image of $\mu$ under the map $\bar{S}.$
Then the following holds:
\begin{itemize}
\item the probability measure $\bar{\nu}$ is dominated by $\nu$ in the positively $1$-homogenous convex order,
\item one has that
\[
\mathcal{I}_c(\mu,\nu) = \int F(x, \bar{S}(x))\,\mu(dx) = \inf_{\gamma \leq_{phc} \nu} \mathcal{T}_{F}(\mu,\gamma)
\]
\end{itemize}
\end{thm}
Indeed, the map $\bar{S}$ provides an optimal transport map between $\mu$ and $\bar{\nu}$ for the cost $\mathcal{T}_F$.
The proof below adapts the proof of \cite[Proposition 1.1]{GJ20} and \cite[Lemma 6.1]{BBP19} to our setting.
\proof
Let $q$ be a nonnegative kernel such that $\mu q = \nu$.
Thanks to Jensen's inequality, for all positively $1$-homogenous convex function $\varphi : \R^d \to \R$, it holds
\[
\int \varphi(S(x))\,\mu(dx) \leq \iint \varphi(y)\,q^x(dy)\,\mu(dx) = \int \varphi(y)\,\nu(dy),
\]
where $S(x) = \int y\,q^x(dy)$, $x\in \X$, and so $S_\# \mu$ is dominated by $\nu$ in the positively homogenous convex order.
Therefore,
\[
\int F(x, S(x))\,\mu(dx) \geq \inf_{\gamma \leq_{phc} \nu} \inf_{\pi \in \Pi(\mu,\gamma)} \int F(x,z)\,\pi(dxdz).
\]
Optimizing over $q$ shows that
\[
\inf_{\gamma \leq_{phc} \nu} \mathcal{T}_F(\mu,\gamma) \leq \mathcal{I}_c(\mu,\nu).
\]
Let us prove the converse inequality. Let $\gamma \leq_{phc} \nu$ and $\pi \in \Pi(\mu,\gamma)$.
Since $\gamma \leq_{phc} \nu$, Theorem \ref{thm:convorder} shows that there exists a nonnegative kernel $(r^z)_{z\in \R^d}$ such that $\int r^z(dy)\,\gamma(dz) = \nu(dy)$ and $\int y \,r^z(dy) = z$ for $\gamma$ almost all $z$.
Write
\[
\pi(dxdz) = \mu(dx) p^x(dz)
\]
where $p$ is a probability kernel, and consider the nonnegative kernel $q$ defined by
\[
q^x(dy) = \int p^x(dz)r^z(dy),
\]
which satisfies $\mu q = \nu$. Moreover, for $\mu$ almost all $x$:
\[
\int y q^x (dy) = \int p^x(dz) \int y r^z(dy) = \int z p^x(dz).
\]
Thus,
\begin{align*}
\int F(x,z)\,\pi(dxdz) = \iint F(x,z)\,\mu(dx) p^x(dz) &\geq \int F\left(x,\int z\,p^x(dz)\right)\,\mu(dx)\\
& =  \int F\left(x,\int y q^x (dy) \right)\,\mu(dx) \geq \mathcal{I}_c(\mu,\nu).
\end{align*}
Optimizing over $\pi$ and over $\gamma \leq_{phc}\nu$ gives that
\[
\inf_{\gamma \leq_{phc} \nu} \mathcal{T}_F(\mu,\gamma) \geq \mathcal{I}_c(\mu,\nu)
\]
which proves  \eqref{eq:Ic-cvxorder}. Now if $\bar{q}$ is a strong solution, then $\bar{\nu}=\bar{S}_\# \mu \leq_{phc} \nu$ and
\[
\mathcal{I}_c(\mu,\nu) = \int F(x,\bar{S}(x))\,\mu(dx),
\]
which completes the proof.
\endproof

The following result considers the particular case of dimension $1$.
\begin{prop}\label{prop:dim1}
Let $\nu \in \mathcal{P}(\R_+)$ be a compactly supported probability measure with $m = \int y\,\nu(dy)$ and denote by
\[
C_m = \left\{\gamma \in \mathcal{P}_1(\R_+) : \int x\,\gamma(dx) = m\right\}.
\]
Then the following holds
\begin{itemize}
\item $C_m = \{\gamma \in \mathcal{P}_1(\R) : \gamma \leq_{phc} \nu\}$,
\item if $\mu \in \mathcal{P}(\X)$ is such that
\[
\mathcal{I}_c(\mu,\nu) = \min_{\gamma \in C_m} \mathcal{T}_F(\mu,\gamma),
\]
then, there exists a map $S_m : \X \to \R_+$ transporting $\mu$ onto $\gamma_m \in C_m$ such that
\[
\mathcal{I}_c(\mu,\nu)= \mathcal{T}_F(\mu,\gamma_m) = \int F(x,S_m(x))\,\mu(dx),
\]
\item the nonnegative kernel $\bar{q}$ defined by
\[
\bar{q}^x (dy) = \frac{S_{m}(x)}{m}\,\nu(dy),\qquad x\in \X,
\]
is a strong solution of the transport problem \eqref{eq:TP}.
\end{itemize}
\end{prop}
This result tells us that, in dimension one, once the solutions for the transport problem between $\mu$ and Dirac masses $\delta_m$, $m\geq0$, are known, then optimal solutions can be deduced for general $\nu$ on $\R_+$.
\proof
According to Remark \ref{rem:dim1}, for any compactly supported probability measure $\nu$ on $\R_+$,
\[
\{\gamma \in \mathcal{P}_1(\R) : \gamma \leq_{phc} \nu\} = \left\{\gamma \in \mathcal{P}(\R_+) :  \int x\,\gamma(dx) = \int x\,\nu(dx)\right\} = C_{m}.
\]
According to Theorem \ref{thm:primal-structure}, 
\[
\mathcal{I}_c(\mu,\nu) = \inf_{\gamma \leq_{phc} \nu} \mathcal{T}_F(\mu,\gamma) =  \inf_{\gamma \in C_{m}} \mathcal{T}_F(\mu,\gamma)
\]
holds. By assumption, this last infimum is reached at some point $\gamma_m' \in C_m$.
Let $\pi \in \Pi(\mu,\gamma_{m}')$ be an optimal coupling for $\mathcal{T}_F(\mu,\gamma_{m}')$ and write $\pi(dxdy) = \mu(dx)p^x(dy)$.
By Jensen's inequality, 
\[
\mathcal{T}_F(\mu,\gamma_{m}') = \iint F(x,y)p^x(dy)\mu(dx) \geq \int F(x, S_m(x))\,\mu(dx),
\]
holds (where $S_m(x)= \int yp^x(dy)$). Denoting by $\gamma_m = (S_m)_\# \mu$, one sees that
$\gamma_m \leq_c \gamma_m'$ and in particular, $\gamma_m  \in C_m$. Therefore, one gets
\[
\inf_{\gamma \in C_m} \mathcal{T}_F(\mu,\gamma)=\mathcal{T}_F(\mu,\gamma_{m}') \geq \int F(x, S_m(x))\,\mu(dx) \geq \mathcal{T}_F(\mu,\gamma_m) \geq \inf_{\gamma \in C_m} \mathcal{T}_F(\mu,\gamma).
\]
This proves that $\gamma_m'$ can be replaced by $\gamma_m$ and that $S_m$ is an optimal transport map (for the cost $F$) between $\mu$ and $\gamma_m$.
The nonnegative kernel $\bar{q}$ defined in Proposition \ref{prop:dim1} satisfies
\[
\int f(y)\,\bar{q}^x(dy)\mu(dx) =\frac{1}{m} \iint f(y)S_m(x)\,\mu(dx)\nu(dy) = \int f(y)\,\nu(dy) \frac{\int y \,\gamma_m(dy)}{m} =  \int f(y)\,\nu(dy)
\]
and so $\mu \bar{q} = \nu$. Moreover, for all $x$, $\int y \,\bar{q}^x(dy) = S_m(x)$ and so
\[
\int F(x,S_m(x))\, \mu(dx) = \int F\left(x, \int y \,\bar{q}^x(dy) \right)\,\mu(dx),
\]
which shows that $\bar{q}$ is a strong solution.
\endproof

The next result establishes a variant of Theorem \ref{thm:primal-structure} involving the classical convex order $\leq_c$.
\begin{thm}\label{thm:primal-structure-bis}
Let $\mu \in \mathcal{P}(\X)$ and $\nu \in \mathcal{P}(\Y)$, and assume that $\bar{\pi} \in \Pi(\ll\mu,\nu)$ is a strong solution to the transport problem \eqref{eq:TP} for the conical cost function $c$ defined by \eqref{eq:conical}. Let $\bar{\eta}(dx) = \bar{N}(x)\,\mu(dx)$ be the first marginal of $\bar{\pi}$ and consider the map $\bar{T}$ defined by
\[
\bar{T}(x) = \int y\,\bar{\pi}^x(dy),\qquad x\in \X,
\]
and denote by $\tilde{\nu}$ the image of $\bar{\eta}$ under the map $\bar{T}.$
Consider the function $G : \X \times \Y \to \R$ defined for $x\in \X$ and $y\in \Y$ by
\[
G(x,y) = \frac{F(x, \bar{N}(x) y)}{\bar{N}(x)},
\]
if $\bar{N}(x)>0$, and $G(x,y)=0$ (or any other arbitrary value) otherwise.
The following holds:
\begin{itemize}
\item the probability measure $\tilde{\nu}$ is dominated by $\nu$ in the convex order,
\item one has that
\[
\mathcal{I}_c(\mu,\nu) = \int G(x, \bar{T}(x))\,\bar{\eta}(dx) = \inf_{\gamma \leq_c \nu} \mathcal{T}_{G}(\bar{\eta},\gamma),
\]
where we denote by $\mathcal{T}_G$ the Monge-Kantorovich optimal transport cost associated to the cost function $G$:
\[
\mathcal{T}_G(\eta,\gamma) = \inf_{\pi \in \Pi(\eta,\gamma)} \iint G(x,z)\,\pi(dxdz),\qquad \forall \eta \in \mathcal{P}(\X),\forall \gamma \in \mathcal{P}(\Y).
\]
\end{itemize}
\end{thm}
In other words, the probability measure $\tilde{\nu}$ turns out to be the closest point to $\bar{\eta}$ among the set $\{\gamma \in \mathcal{P}(\Y) : \gamma \leq_c \nu\}$ for the transport ``distance'' $\mathcal{T}_G$. Moreover, the map $\bar{T}$ (which is sometimes called the barycentric projection of the coupling $\bar{\pi}$) provides an optimal transport map between $\bar{\eta}$ and $\tilde{\nu}$ for the cost $\mathcal{T}_G$.
\begin{rem}
Suppose that $\bar{\pi}$ is a strong solution to the transport problem \eqref{eq:TP}, denote by $\bar{\eta} = \bar{N}\,\mu$ its first marginal, and consider the nonnegative kernel $\bar{q}^x = \bar{N}(x)\bar{\pi}^x$, $x\in \X$. Then, for all $x\in \X$, $\bar{S}(x) = \bar{N}(x) \bar{T}(x)$. However, $\bar{\nu}$ and $\tilde{\nu}$ are in general two different probability measures, so that the conclusions of Theorems \ref{thm:primal-structure} and \ref{thm:primal-structure-bis} are not equivalent. Nevertheless, for all positively homogenous function $\varphi$ one has that $\int \varphi(y)\,\bar{\nu}(dy) = \int \varphi(y)\,\tilde{\nu}(dy).$
\end{rem}

\proof
Thanks to Jensen's inequality, for all convex function $\varphi : \R^d \to \R$,
\[
\int \varphi(\bar{T}(x))\,\bar{\eta}(dx) \leq \iint  \varphi(y)\,\bar{\pi}^x(dy)\,\bar{\eta}(dx) = \int \varphi(y)\,\nu(dy),
\]
and so $\tilde{\nu} = \bar{T}_\# \bar{\eta}$ is dominated by $\nu$ in the convex order.
Therefore,
\[
\mathcal{I}_c(\mu,\nu)  = \int G(x, \bar{T}(x))\,\bar{\eta}(dx) \geq \inf_{\gamma \leq_c \nu} \inf_{\pi \in \Pi(\bar{\eta},\gamma)} \int G(x,z)\,\pi(dxdz).
\]
Let us prove the converse inequality. Let $\gamma \leq_c \nu$ and $\pi \in \Pi(\bar{\eta},\gamma)$.
Since $\gamma \leq_c \nu$, there exists a martingale coupling $m$ between $\gamma$ and $\nu$, that is a probability measure $m \in \Pi(\gamma,\nu)$ such that $m(dzdy) = \gamma(dz)m^z(dy)$ and $\int y \,m^z(dy)=z$ for $\gamma$ almost every $z$. Write
\[
\pi(dxdz) = \bar{\eta}(dx) p^x(dz) = \gamma(dz)r^z(dx)
\]
and consider the coupling $\bar{\pi}$ defined by
\[
\bar{\pi}(dxdy) = \int r^z(dx)m^z(dy)\gamma(dz).
\]
It is easily seen that $\bar{\pi} \in \Pi(\bar{\eta},\nu).$ Also,
\[
\bar{\pi}(dxdy) = \bar{\eta}(dx) \int p^x(dz)m^z(dy)
\]
and so $\bar{\pi}^x (dy) = \int p^x(dz)m^z(dy)$. Moreover
\[
\int y \bar{\pi}^x (dy) =  \int \left(\int y m^z(dy)\right) p^x(dz) = \int z \,p^x(dz).
\]
Thus,
\begin{align*}
\iint G(x,z)\,\pi(dxdz) &= \iint  G(x,z)\,\bar{\eta}(dx) p^x(dz) \geq \int G\left(x,\int z\,p^x(dz)\right)\,\bar{\eta}(dx)\\ &=  \int G\left(x,\int y \,\bar{\pi}^x (dy) \right)\,\bar{\eta}(dx) =  I_c^\mu[\bar \pi] \geq \mathcal{I}_c(\mu,\nu)
\end{align*}
which completes the proof.
\endproof

The articulation between primal and dual optimizers is analyzed thanks to the next result.

\begin{thm}\label{thm:articulation} Let $\mu \in \mathcal{P}(\X)$ and $\nu \in \mathcal{P}(\Y)$ be such that $\mathcal{I}_c(\mu,\nu)<+\infty$ and assume that $\bar{q}$ is a kernel solution and $\bar{\varphi} \in \Phi(\Z)\cap L^1(\nu)$ a dual optimizer:
\[
\mathcal{I}_c(\mu,\nu) = \int F\left(x,\int y\,\bar{q}^x(dy)\right)\,\mu(dx)= \int Q_F\bar{\varphi}(x)\,\mu(dx) - \int \bar{\varphi}(y)\,\nu(dy).
\]
Define $\bar{S}(x) = \int y\,\bar{q}^x(dy)$, $x\in \X$.
\begin{enumerate}
\item For $\mu$ almost every $x \in \X$, 
\[
Q_F\bar{\varphi}(x) = \bar{\varphi}(\bar{S}(x))+F(x,\bar{S}(x)).
\]
\item If $M$ denotes the set of $x \in \X$ for which the support $K(x)\subset \Y$ of $\bar{q}^x$ contains at least two points, then for $\bar{\eta}$ almost $x \in M$, the function $\bar{\varphi}$ is affine on the convex hull of $K(x)$: there exist $u_x\in \R^d$ and $v_x \in \R$ such that $\bar{\varphi}(z) = u_x\cdot z+v_x$ for all $z \in \mathrm{co}(K(x))$.
\item If $F$ is strictly convex with respect to its second variable, then the map $\bar{S}(x) = \int y \,\bar{q}^x(dy)$, $x\in \X$, is $\mu$-almost surely unique among all strong solutions $\bar{q}$ of the transport problem.
\item If $F$ is strictly convex  with respect to its second variable and if for all $x \in \X$ there exist $A_x \in \R$ and $M_x>0$ such that $\bar{\varphi}(z)+F(x,z) \geq A_x + M_x \|z\|$ for all $z \in \Z$, then for all $x \in \X$ the map $\bar{\varphi}^* \square F^*(x,\,\cdot\,)$ is differentiable in a neighborhood of $0$ and it holds
\[
\bar{S}(x) = \nabla\left( \bar{\varphi}^* \square F^*(x,\,\cdot\,) \right)(0)
\]
for $\mu$ almost all $x$.
\end{enumerate}
\end{thm}
In the result above we denoted
\[
\bar{\varphi}^*(u) = \sup_{z\in \Z} \{z\cdot u - \bar{\varphi}(z)\}, \qquad u \in \R^d,
\]
and for $x \in \X$
\[
F^*(x,u) = \sup_{z\in \Z} \{z\cdot u - F(x,z)\}, \qquad u \in \R^d,
\]
the Fenchel-Legendre transforms of the functions $\bar{\varphi}$ and $F(x,\,\cdot\,)$ extended by $+\infty$ outside $\Z$.
Moreover, we recall that the infimum convolution between $\bar{\varphi}^* $ and $F^*(x,\,\cdot\,)$ is defined by
\[
\bar{\varphi}^* \square F^*(x,\,\cdot\,)(u) = \inf_{v \in \R^d} \{\bar{\varphi}^*(v)+F^*(x,u-v)\}, \qquad u \in \R^d.
\]
\begin{rem}\label{rem:support}
Denoting by $C(\bar{\varphi})=\{u \in \R^d : u\cdot z \leq \bar{\varphi}(z), \forall z\in \Z\}$, we see that
\[
\bar{\varphi}^* = \chi_{C(\bar{\varphi})},
\]
where $\chi_{C(\bar{\varphi})}(u)=0$ if $u \in C(\bar{\varphi})$ and $+\infty$.
So, it holds that
\begin{equation}\label{eq:infconv}
\bar{\varphi}^* \square F^*(x,\,\cdot\,)(u) = \inf_{v \in  C(\bar{\varphi})} \{F^*(x,u-v)\}, \qquad u \in \R^d.
\end{equation}
\end{rem}
\proof
By optimality of $\bar{\varphi}$ and $\bar{q}$, it holds
\begin{align*}
\mathcal{I}_c(\mu,\nu) & =  \int Q_F\bar{\varphi}(x)\,\mu(dx) - \int \bar{\varphi}(y)\,\nu(dy)\\
& \leq  \int \bar{\varphi}(\bar{S}(x)) +F(x,\bar{S}(x))\,\mu(dx) - \int \bar{\varphi}(y)\,\nu(dy)\\
& \leq \int \int \bar{\varphi}(y)\,\bar{q}^x(dy) +F(x,\bar{S}(x))\,\mu(dx)- \int \bar{\varphi}(y)\,\nu(dy)\\
& = \int F(x,\bar{S}(x))\,\mu(dx)\\
& = \mathcal{I}_c(\mu,\nu),
\end{align*}
where the first inequality comes from the definition of $Q_F\bar{\varphi}$, the second from $\bar{\varphi}$ being convex and positive $1$-homogenous. Analyzing the equality cases completes the proof of $(1)$ and $(2)$.
Now, assume that $F$ is strictly convex with respect to its second variable, and consider $\bar{r}$ another minimizing nonnegative kernel and define $\bar{U}(x) = \int y \,\bar{r}^x(dy)$, $x\in \X$.
According to what precedes, for $\mu$ almost all $x\in \X$, the points $\bar{S}(x)$ and $\bar{U}(x)$ minimize the function
\[
z \mapsto \bar{\varphi}(z) + F(x,z),\qquad z\in \Z.
\]
This function being strictly convex, then $\bar{S}(x)=\bar{U}(x)$ and so $\bar{S}=\bar{U}$ $\mu$ a.e., which proves $(3)$.
Let us now prove $(4)$. Consider the function $H : \X \times \R^d \to \R\cup\{+\infty\}$ defined by $H(x,z)=\bar{\varphi}(z) + F(x,z)$, $x\in \X$, $z\in \R^d$ (with $\bar{\varphi}$ and $F(x,\,\cdot\,)$ extended by $+\infty$ outside $\Z$).
For a given $x \in \X$, observe that $\bar{S}(x)$ minimizes the convex function $H(x,\,\cdot\,)$ if and only if $0 \in \partial H(x,\,\cdot\,)(\bar{S}(x))$, where $\partial H(x,\,\cdot\,)(z)$ denotes the subdifferential of the function $H(x,\,\cdot\,)$ at the point $z$. By the well-known conjugation relation of subdifferentials, one has that
\[
0 \in \partial H(x,\,\cdot\,)(\bar{S}(x)) \Leftrightarrow \bar{S}(x) \in \partial H^*(x,\,\cdot\,)(0)
\]
(see e.g \cite[Corollary E.1.4.4]{HUL01}). Moreover, since $H(x,\,\cdot\,)$ is a sum of two convex lower-semicontinuous functions, its Fenchel-Legendre transform is given as follows:
\[
H^*(x,u) = \bar{\varphi}^*\square F^*(x,\,\cdot\,)(u),\qquad u\in \R^d
\]
(see e.g \cite[Theorem E.2.3.2]{HUL01} with $F(x,\,\cdot\,)$ being finite over $\Z$ which contains the relative interior of the domain of $\bar{\varphi}$). The assumed lower bound on $H$ easily implies that, for all $x \in \X$, the function $H^*(x,\,\cdot\,)$ takes finite values in a neighborhood of $0$. Since $H(x,\,\cdot\,)$ is strictly convex for all $x \in \X$, it follows from \cite[Theorem E.4.1.1]{HUL01} that the function $H^*(x,\,\cdot\,)$ is continuously differentiable on the interior of its domain. In particular, it is continuously differentiable in a neighborhood of $0$, and so, for all $x\in \X$, $\partial H^*(x,\,\cdot\,)(0) = \{ \nabla\left( \bar{\varphi}^* \square F^*(x,\,\cdot\,) \right)(0)\}$, which completes the proof.
\endproof

\begin{cor}\label{cor:barS} Let $F : \X \times \Z \to \R$ be a cost function satisfying assumption \eqref{eq:A'} and strictly convex with respect to its second variable. Let $\mu \in \mathcal{P}(\X)$ and $\nu \in \mathcal{P}(\Y)$ be such that the convex hull of the support of $\nu$ does not contain $0$.\\ Assume further that $F$ satisfies
\begin{itemize}
\item Assumption \eqref{eq:B'}\\
or
\item $F:\X\times \Z \to \R_-$ is a nonpositive function satisfying the assumptions of Theorem \ref{thm:duality-eco}.
\end{itemize}
Then, for any kernel solution $\bar{q}$ and any dual optimizer $\bar{\varphi}$ of the transport problem \eqref{eq:TPintro}, the map $\bar{\varphi}^* \square F^*(x,\,\cdot\,)$ is differentiable in a neighborhood of $0$, for all $x \in \X$, and then:
\[
\bar{S}(x)=\int y\,\bar{q}^x(dy) = \nabla\left( \bar{\varphi}^* \square F^*(x,\,\cdot\,) \right)(0),
\]
for $\mu$ almost all $x\in \X$.
\end{cor}

\proof
Assumption \eqref{eq:B'} is known to be equivalent to the $1$-coercivity of $F(x,\,\cdot\,)$, that is 
\[
\lim_{z \in \Z, \|z\| \to +\infty}\frac{F(x,z)}{\|z\|} \to + \infty.
\]
Therefore, since the convex function $\bar{\varphi}$ admits at least one affine minorant, for every $x\in \X$ there exist $A_x \in \R$ and $M_x >0$ such that $\bar{\varphi}(z)+F(x,z) \geq A_x + M_x \|z\|$ for all $z \in \Z$. We conclude using Item $(4)$ of Theorem \ref{thm:articulation}.
Similarly, if $F$ satisfies the assumptions of Theorem \ref{thm:duality-eco}, then for all $x\in \X$ and $z \in \Z$ it holds that $F(x,\lambda z)/\lambda \to 0$ as $\lambda \to +\infty$. This is actually equivalent to
\[
\lim_{z \in \Z, \|z\| \to +\infty}\frac{F(x,z)}{\|z\|} =0.
\]
Let us briefly sketch the proof. Let $z_k \in \Z$, $k\geq 0$, be a sequence such that $\lambda _k = \|z_k\| \to +\infty$ monotonically, as $k\to +\infty$. Define $u_k = z_k/\lambda_k$, $k\geq0$. By compactness, one can assume without loss of generality that $u_k \to u \in \Z$ as $k\to +\infty$. Then, by convexity, for all $k\geq n$
\[
\frac{F(x,z_k)-F(x,0)}{\lambda_k} = \frac{F(x,\lambda_ku_k)-F(x,0)}{\lambda_k} \geq  \frac{F(x,\lambda_nu_k)-F(x,0)}{\lambda_n}
\]
Thus letting $k\to +\infty$, one gets that
\[
\liminf_{k\to \infty}\frac{F(x,z_k)}{\|z_k\|} \geq \frac{F(x,\lambda_nu)-F(x,0)}{\lambda_n}
\]
and letting $n\to +\infty$ gives that
\[
\liminf_{k\to \infty}\frac{F(x,z_k)}{\|z_k\|} \geq F'_\infty(x,u)=0.
\]
Since $F$ is nonpositive, this proves the claim. Now, if $\bar{\varphi}$ is some dual optimizer, we know by Theorem \ref{thm:duality-eco} that $\bar{\varphi}>0$ on $\Z \setminus \{0\}$. Thus, denoting by $M = \inf_{\|u\|=1, u \in \Z}\bar{\varphi}(u)>0$, one sees that $\bar{\varphi}(z) \geq M \|z\|$, for all $z \in \Z$. And so
\[
\limsup_{\|z\| \to +\infty}\frac{\bar{\varphi}(z) + F(x,z)}{\|z\|} \geq M>0,
\]
and we conclude using Item $(4)$ of Theorem \ref{thm:articulation}.
\endproof

Let us emphasize a particularly simple case related to Brenier's theorem \cite{Br87,Br91}.
In the following result, adapted from \cite[Theorem 1.2]{GJ20}, we assume that $\X \subset \R^d$ is a compact subset of $\R^d$ and we consider the cost function $c_2 : \X \times \mathcal{M}(\Y) \to \R_+$ defined by
\[
c_2(x,m) = \frac{1}{2} \left\|x- \int y\,m(dy)\right\|_2^2,\qquad x\in \X, m\in \mathcal{M}(\Y),
\]
which corresponds to the function $F_2:\R^d \times \R^d\to \R^+ : (x,z) \mapsto \frac{1}{2}  \|x-z\|_2^2$, where $\|\,\cdot\,\|_2$ is the standard Euclidean norm.

\begin{thm}\label{thm:Brenier}
Let $\mu \in \mathcal{P}(\X)$ and $\nu \in \mathcal{P}(\Y)$ be such that the convex hull of the support of $\nu$ does not contain $0$.
Then, there exists a unique probability measure $\bar{\nu} \in \mathcal{P}(\Z)$ such that
\begin{equation}\label{eq:barnu}
\mathcal{I}_{c_2}(\mu,\nu) = \frac{1}{2}\inf_{\eta \leq_{phc} \nu} W_2^2(\mu,\eta) = \frac{1}{2} W_2^2(\mu,\bar{\nu}).
\end{equation}
Moreover, there exists a closed convex set $C \subset \R^d$ such that, for any nonnegative kernel $\bar{q}$ minimizing $\mathcal{I}_c(\mu,\nu)$, one has
\[
\bar{S}(x) = \int y\,\bar{q}^x(dy) = x-p_C(x),
\]
for $\mu$ almost every $x$, where $p_C:\R^d \to \R^d$ is the orthogonal projection onto $C$.
The probability $\bar{\nu}$ is the image of $\mu$ under the map $x\mapsto x-p_C(x)$.
\end{thm}
\proof
The cost function $c_2$ clearly satisfies Assumption \eqref{eq:B-intro} and so, according to Theorem \ref{thm:strongsol}, the transport problem \eqref{eq:TPintro} between $\mu$ and $\nu$ admits kernel solutions. According to Theorem \ref{thm:duality-conical}, it also admits dual optimizers. Let $\bar{\varphi} \in \Phi(\Z)\cap L^1(\nu)$ be a dual optimizer (extended by $+\infty$ outside $\Z$). As observed in Remark \ref{rem:support}, $\bar{\varphi}^*=\chi_{C}$ for some closed convex set $C$, and according to \eqref{eq:infconv} we have
\[
\bar{\varphi}^*\square F_2(x,\,\cdot\,)(u) = \inf_{v \in C} F_2^*(x,u-v) =  \inf_{v \in C} \left\{\frac{1}{2}\|u-v\|_2^2+(u-v)\cdot x\right\} = -\frac{\|x\|_2^2}{2} +  \frac{1}{2}d_C^2(x+u),
\]
where $d_C(a) = \inf_{v\in C}\|a-v\|_2$. It is well known that $d_C^2$ is differentiable over $\R^d$ and that
\[
\nabla \left(\frac{1}{2}d^2_C\right)(a) = a-p_C(a),\qquad \forall a\in \R^d.
\]
Therefore,
\[
\nabla \left(\bar{\varphi}^*\square F_2(x,\,\cdot\,)\right)(0) = x-p_C(x), \qquad \forall x\in \X.
\]
According to Corollary \ref{cor:barS}, for any kernel solution $\bar{q}$, 
\[
\int y\,\bar{q}^x(dy) =  x-p_C(x),
\]
holds for $\mu$ almost all $x\in \R^d.$ According to Theorem \ref{thm:primal-structure}, we conclude that the probability measure $\bar{\nu}$ defined as the push forward of $\mu$ under $x\mapsto  x-p_C(x)$ satisfies \eqref{eq:barnu}. The uniqueness of $\bar{\nu}$ is obtained as in \cite[Proposition 1.1]{GJ20}.
\endproof
\begin{rem}Let us give a geometric justification that $x-p_C(x)$ belongs to $\Z$ for all $x \in \R^d$. Denoting by $\bar{\varphi}$ the dual optimizer used in the proof (extended by $+\infty$ outside $\Z$), one has $\Z_{\bar{\varphi}}:=\mathrm{cl}\, \mathrm{dom} (\bar{\varphi}) \subset \Z$. But, since $\bar{\varphi}$ is the support function of $C = \{u : u\cdot x \leq \bar{\varphi}(x), \forall x \in \Z\}$, one has according to \cite[Proposition C.2.2.4]{HUL01}
\[
C_\infty^\circ = \Z_{\bar{\varphi}},
\]
where
\begin{itemize}
\item $C_\infty$ denotes the asymptotic cone of $C$, defined by
\[
C_\infty = \bigcap_{t>0} \frac{C-x_o}{t},
\]
with $x_o \in C$ some arbitrary point,
\item $C_\infty^\circ$ denotes the polar cone of $C_\infty$ defined by
\[
C_\infty^\circ = \{x \in \R^d : x\cdot y \leq 0, \forall y \in C_\infty\}.
\]
\end{itemize}
By definition of the orthogonal projection on $C$, one can write
\[
(x-p_C(x))\cdot (a-p_C(x)) \leq 0,\qquad \forall x \in \R^d,\forall a \in C.
\]
In particular, taking $a = p_C(x)+d$, with $d \in C_\infty$ leads to $(x-p_C(x))\cdot d \leq 0$ for all $d \in C_\infty$ and so $x-p_C(x) \in C_\infty^\circ = \Z_{\bar{\varphi}} \subset \Z$.
\end{rem}

\appendix

\addcontentsline{toc}{section}{Appendix}

\numberwithin{equation}{section}

\setcounter{equation}{0}  

\section{Proofs of some technical results}\label{sec:Appendix}
\subsection{Proof of Proposition \ref{prop:barI_csci}} \label{app:1}
The proof of Proposition \ref{prop:barI_csci} is adapted from \cite{AFP00} (paragraph 2.6).
First, let us see how the recession function $c'_\infty$ can be expressed when $c$ satisfies Assumption \eqref{eq:A-intro}.
\begin{lem}\label{lem:recession}
Under Assumption \eqref{eq:A-intro}, one has
\[
c'_\infty (x,m) = \sup_{k\in \N} \int b_k(x,y)\,m(dy),\qquad x\in \X, m\in \mathcal{M}(\Y).
\]
\end{lem}
\proof
Since $c(x,\,\cdot\,)$ is convex, the function $\lambda \mapsto \frac{c(x,\lambda m)-c(x,0)}{\lambda}$ is non-decreasing on $[0,\infty).$ Therefore, for all $x\in \X$ and $m\in \mathcal{M}(\Y)$, then
\begin{equation}\label{eq:non-decreasing}
c'_\infty (x,m) = \lim_{\lambda \to \infty}\frac{c(x,\lambda m)-c(x,0)}{\lambda}=\sup_{\lambda > 0}\frac{c(x,\lambda m)-c(x,0)}{\lambda}.
\end{equation}
Thus, using Assumption \eqref{eq:A-intro}, one has
\begin{align*}
c'_\infty (x,m) &= \sup_{\lambda > 0} \sup_{k\geq0}\frac{\lambda \int b_k(x,y)\,m(dy) + a_k(x)-c(x,0)}{\lambda} \\
& =  \sup_{k\geq0} \left\{ \int b_k(x,y)\,m(dy) +\sup_{\lambda> 0}\frac{a_k(x)-c(x,0)}{\lambda}\right\}\\
& = \sup_{k\geq0} \int b_k(x,y)\,m(dy),
\end{align*}
where the last equality comes from $c(x,0) = \sup_{k\geq0} a_k(x)$, and so for a fixed $k$, $a_k(x)-c(x,0)\leq0$. Hence, the function $\lambda\mapsto \frac{a_k(x)-c(x,0)}{\lambda}$ is non-decreasing on $[0,\infty).$
\endproof

We will also need the following lemma
\begin{lem}\label{lem:intsup}Let $\lambda$ be a finite measure on $\X$. If $\psi_0,\psi_1,\ldots,\psi_n:\X\to \R$ are $\lambda$ integrable functions with $\psi_0\geq0$, then
\[
\int \sup_{0\leq k\leq n} \psi_k(x)\,\lambda(dx) = \sup_{(f_0,\ldots,f_n)\in \mathcal{F}_n} \int \sum_{k=1}^n \psi_k(x)f_k(x)\,\lambda(dx)
\]
where $\mathcal{F}_n$ denotes the set of $(n+1)$-tuples $(f_0,\ldots,f_n)$ of continuous functions $f_0,\ldots,f_n: \X\to [0,1]$ such that $f_0(x)+\cdots+f_n(x)\leq 1.$
\end{lem}
\proof
See the proof of Proposition 9.4 of \cite{GRST17}.
\endproof

\proof[Proof of Proposition \ref{prop:barI_csci}.]
Without loss of generality, one can assume that the functions $b_0$ and $a_0$ involved in \eqref{eq:A-intro} are nonnegative.
If this is not the case, consider the cost function
\[
\tilde{c}(x,m) = \sup_{k\geq 0} \left\{\int \tilde{b}_k(x,y)\,m(dy) + \tilde{a}_k(x)\right\},\qquad x\in \X, m\in \mathcal{M}(\Y),
\]
with
$\tilde{b}_k = b_k - r_0$ and $\tilde{a}_k=a_k - s_0$, with $r_0 = \min_{x\in \X,y\in \Y} b_0(x,y)$ and $s_0 = \min_{x\in \X} a_0(x)$. It holds $\tilde{a}_0\geq0$, $\tilde{b}_0 \geq0$ and one sees that $\tilde{c}(x,m)= c(x,m)-r_0m(\Y)-s_0$, $x\in \X$, $m\in \mathcal{M}(\Y),$ and so $\bar{I}_{\tilde{c}}^\mu[\pi] = \bar{I}_c^\mu[\pi] - (s_0+r_0)$, and  $(\mu,\pi) \mapsto \bar{I}_c^\mu[\pi]$ is lower semicontinuous if and only if $(\mu,\pi) \mapsto \bar{I}_{\tilde{c}}^\mu[\pi]$ is.

For $n\geq0$, define for all $\mu \in \mathcal{P}(\X)$ and $\pi \in \mathcal{P}(\X\times \Y)$
\[
J_n^\mu[\pi] = \sup_{(f_0,\ldots,f_n) \in \mathcal{F}_n} \left\{ \sum_{k=0}^n \iint  b_k(x,y)f_k(x)\,\pi(dxdy)+\sum_{k=0}^n\int a_k(x)f_k(x)\,\mu(dx)  \right\},
\]
where $\mathcal{F}_n$ is defined in Lemma \ref{lem:intsup} above.
Then consider the functional $J^\mu$ defined by
\begin{equation}\label{eq:barI}
J^\mu[\pi] = \sup_{n\geq0} J_n^\mu[\pi],\qquad \pi \in \mathcal{P}(\X\times \Y).
\end{equation}
For each $n\geq0$, the functional $(\mu,\pi)\mapsto J_n^\mu[\pi]$ is lower semicontinuous as a supremum of continuous functionals. Similarly, the functional $(\mu,\pi)\mapsto J^\mu[\pi]$ being the supremum of lower semicontinuous functionals is itself lower semicontinuous.

For all $n\geq0$, write $c_n(x,m)= \sup_{0\leq k\leq n} \left\{\int b_k(x,y)\,m(dy) + a_k(x)\right\}$, $x\in \X$, $m\in \mathcal{M}(\Y)$.
According to Lemma \ref{lem:recession}, it holds that $c'_{n,\infty}(x,m)= \sup_{0\leq k\leq n} \left\{\int b_k(x,y)\,m(dy)\right\}$, $x\in \X$, $m\in \mathcal{M}(\Y)$.

By monotone convergence,
\[
\bar{I}_c^\mu[\pi] = \sup_{n\geq0} \bar{I}_{c_n}^\mu[\pi],\qquad \forall \mu \in \mathcal{P}(\X),\forall \pi \in \mathcal{P}(\X\times \Y).
\]
Let us show that for any $\mu \in \mathcal{P}(\X)$ and $\pi \in  \mathcal{P}(\X\times \Y)$, one has that $J_n^\mu[\pi] = \bar{I}_{c_n}^\mu[\pi]$ for all $n\geq0$. This will immediately imply that $\bar{I}_c^\mu[\pi]=J^\mu[\pi]$ for all $\mu \in \mathcal{P}(\X)$ and $\pi \in  \mathcal{P}(\X\times \Y)$ and show the lower semicontinuity of $\bar{I}_c^{\cdot}[\,\cdot\,].$

Fix $\mu \in \mathcal{P}(\X)$; for all $\pi \in \mathcal{P}(\X\times \Y)$, then
\[
\bar{I}_{c_n}^\mu[\pi] = \int \sup_{0\leq k\leq n} \psi_k(x)\,\mu(dx) +  \int \sup_{0\leq k\leq n} \varphi_k(x)\,\pi_1^s(dx)
\]
with, for all $0\leq k\leq n$,
\[
\psi_k(x) = \frac{d\pi_1^{ac}}{d\mu}(x)\int b_k(x,y)\,\pi^x(dy) + a_k(x),\qquad x\in \X
\]
and
\[
\varphi_k(x) = \int b_k(x,y)\,\pi^x(dy),\qquad x\in \X.
\]
Let $A \subset \X$ be a Borel subset such that $\mu(A)=0$ and $\pi_1^s (\X\setminus A)=0$.
Define
\[
F_k(x) = \left\{\begin{array}{ll} \psi_k(x) & \text{if } x\in \X\setminus A   \\ \varphi_k(x)  & \text{if } x\in A  \end{array}\right..
\]
Then
\[
\bar{I}_{c_n}^\mu[\pi] = \int \sup_{0\leq k\leq n} F_k(x)\,(\mu(dx) +\pi_1^s(dx)).
\]
According to Lemma \ref{lem:intsup} above (with $\lambda = \mu +\pi_1^s$), one has
\[
\bar{I}_{c_n}^\mu[\pi] =\int \sup_{0\leq k\leq n} F_k(x)\,(\mu(dx) +\pi_1^s(dx))= \sup_{(f_0,\ldots,f_n)\in \mathcal{F}_n} \int \sum_{k=1}^n F_k(x)f_k(x)\,(\mu(dx) +\pi_1^s(dx))= J_n^\mu[\pi],
\]
which completes the proof. The proof that $(\lambda, \mu, \pi) \mapsto \bar{I}_{c_\lambda}^\mu[\pi]$ is also lower semicontinuous is easily adapted and left to the reader.
\endproof

\subsection{Proof of Lemma \ref{lem:recovery}}\label{App:2} The proof of Lemma \ref{lem:recovery} below is inspired from the proof of \cite[Theorem C.12]{LV09} dealing with entropy-type functionals on the space of probability measures.
\begin{lem}\label{lem:subadd}
If $c : \X \times \mathcal{M}(\Y) \to \R$ is convex with respect to its second variable, then for any $x\in \X$ and $m_1,m_2 \in \mathcal{M}(\Y)$, we have that
\[
c(x,m_1+m_2) \leq c(x,m_1)+c'_\infty(x,m_2).
\]
\end{lem}
\proof
Let $\theta \in (0,1)$; using the convexity of $c(x,\,\cdot\,)$ and \eqref{eq:non-decreasing}, one gets
\begin{align*}
c(x,m_1+m_2) & \leq \theta c\left(x,\frac{m_1}{\theta}\right)+ (1-\theta) c\left(x,\frac{m_2}{1-\theta}\right)\\
& =  \theta c\left(x,\frac{m_1}{\theta}\right)+ (1-\theta) \left[c\left(x,\frac{m_2}{1-\theta}\right)- c(x,0)\right] + (1-\theta) c(x,0)\\
& \leq \theta c\left(x,\frac{m_1}{\theta}\right)+ c'_\infty(x,m_2) + (1-\theta) c(x,0),
\end{align*}
and the result follows by letting $\theta \to 1.$
\endproof

\begin{lem}\label{lem:boundC}
If $c : \X \times \mathcal{M}(\Y) \to \R$ is convex with respect to its second variable and satisfies Assumption \eqref{eq:C}, then there exists $a,b\geq0$ such that
\[
c(x,m) \leq b + am(\Y),\qquad \forall x\in \X, \forall m \in \mathcal{M}(\Y).
\]
\end{lem}
\proof
Using \eqref{eq:non-decreasing}, one gets
\[
c(x,m) \leq c(x,0) + m(\Y)c'_\infty\left(x, \frac{m}{m(\Y)}\right) \leq c(x,0)+ am(\Y),
\]
where $a$ is the constant appearing in Assumption \eqref{eq:C}.
Since $c(\,\cdot\,,0)$ is continuous on the compact space $\X$ it is upper bounded by some constant $b\geq0$, which completes the proof.
\endproof

\proof[Proof of Lemma \ref{lem:recovery}]
Let $\pi \in \Pi(\mathrm{Supp}(\mu),\nu)$ and denote by $\eta \in \mathcal{P}(\mathrm{Supp (\mu)})$ its first marginal.
According to Lemma \ref{lem:descr}, $\pi_n\to \pi$ for the weak topology, and the first marginal of $\pi_n$ is $\eta_n = K_n\eta$ and is thus absolutely continuous with respect to $\mu$.
Since $\pi_n\in \Pi(\ll\mu,\,\cdot\,)$, it holds $I_c^\mu[\pi_n] = \bar{I}_c^\mu[\pi_n]$ and since $\bar{I}_c^\mu$ is lower semicontinuous, one gets that
\[
\liminf_{n\to \infty} I_c^\mu[\pi_n] = \liminf_{n\to \infty} \bar{I}_c^\mu[\pi_n] \geq \bar{I}_c^\mu[\pi].
\]
Now we prove that $\limsup_{n\to \infty}I_c^\mu[\pi_n] \leq \bar{I}_c^\mu[\pi]$.
Observe that $\pi_n(dxdy) = q_n^x(dy)\mu(dx)$ with
\[
q_n^x(dy) = \int K_n(x,z)\pi^z(dy)\eta(dz) = \int K_n(x,z)\pi^z(dy)\eta^{ac}(dz) +\int K_n(x,z)\pi^z(dy)\eta^{s}(dz) := q_n^{ac,x} + q_n^{s,x}
\]
where $\eta = \eta^{ac} + \eta^s$ is the decomposition of $\eta$ into absolutely continuous and singular parts (with respect to $\mu$).
According to Lemma \ref{lem:subadd}, it holds
\[
I_c^\mu[\pi_n] = \int c(x,q_n^x)\,\mu(dx) =  \int c(x, q_n^{ac,x} + q_n^{s,x} )\,\mu(dx)  \leq \int c(x, q_n^{ac,x})\,\mu(dx) + \int c'_\infty\left(x,q_n^{s,x} \right)\,\mu(dx).
\]
Write $\eta^{ac}(dz) = h(z)\,\mu(dz)$ and let us bound the first term. Since $\int K_n(x,z)\,\mu(dz)=1$, Jensen inequality yields
\begin{align*}
\int c\left(x, q_n^{ac,x}\right)\,\mu(dx) &=\int c\left(x, \int K_n(x,z)\pi^z(\,\cdot\,)h(z)\mu(dz) \right)\,\mu(dx)\\
& \leq \iint K_n(x,z) c\left(x, \pi^z(\,\cdot\,)h(z) \right)\,\mu(dx)\mu(dz)\\
& = \int  (K_n C_z)(z) \mu(dz),
\end{align*}
where
\[
C_z(x) = c\left(x, \pi^z(\,\cdot\,)h(z) \right),\qquad x,z\in \X.
\]
By assumption, the function $x\mapsto c\left(x, \pi^z(\,\cdot\,)h(z) \right)$ is continuous on $\X$. Therefore, according to Lemma \ref{lem:LV}, one gets that $K_nC_z (u)\to C_z (u)$ for any $u\in \X$ (and even uniformly in $u$) as $n\to \infty$.
Also, according to Lemma \ref{lem:boundC}, it holds $c\left(x, \pi^z(\,\cdot\,)h(z) \right) \leq b +a h(z)$ and so $K_nC_z(z) \leq b + ah(z)$, which is $\mu$ integrable. Therefore, according to the dominated convergence theorem, we get that
\begin{equation}\label{eq:limsup1}
\limsup_{n\to \infty}\int c\left(x, q_n^{ac,x}\right)\,\mu(dx) \leq \lim_{n\to \infty} \int  (K_n C_z)(z) \mu(dz)=  \int  C_z(z) \mu(dz) = \int c\left(z, \pi^z(\,\cdot\,)h(z) \right)\,\mu(dz).
\end{equation}
Now, let us bound the second term. Using the convexity and $1$-homogeneity of the function $m \mapsto c'_\infty(x,m)$ one gets
\begin{align*}
\int c'_\infty\left(x,q_n^{s,x} \right)\,\mu(dx) &= \int c'_\infty\left(x,\int K_n(x,z)\pi^z(\,\cdot\,)\eta^{s}(dz)\right)\mu(dx) \\
&\leq \iint  c'_\infty\left(x,K_n(x,z)\pi^z(\,\cdot\,)\right)\mu(dx)\eta^s(dz) \\
& = \iint  K_n(x,z)c'_\infty\left(x,\pi^z(\,\cdot\,)\right)\mu(dx)\eta^s(dz) \\
& =  \int  (K_nD_z)(z)\,\eta^s(dz),
\end{align*}
where
\[
D_z(x) = c'_\infty\left(x,\pi^z(\,\cdot\,)\right),\qquad \forall x,z \in \X.
\]
By assumption, the function $x\mapsto c'_\infty\left(x, \pi^z \right)$ is continuous on $\X$. Therefore, according to Lemma \ref{lem:LV}, one gets that $K_nD_z (u)\to D_z (u)$ for any $u\in \X$ (and even uniformly in $u$) as $n\to \infty$.
By assumption, it holds $c'_\infty\left(x, \pi^z \right) \leq a $ and so $K_nD_z(z) \leq a$. Therefore, according to the dominated convergence theorem, 
\begin{equation}\label{eq:limsup2}
\limsup_{n\to \infty}\int c\left(x, q_n^{s,x}\right)\,\mu(dx) \leq \lim_{n\to \infty}\int  (K_nD_z)(z)\,\eta^s(dz) =   \int  D_z(z) \eta^s(dz) = \int c'_\infty(z,\pi^z)\eta^s(dz)
\end{equation}
holds. Adding \eqref{eq:limsup1} and \eqref{eq:limsup2}, one gets that $\limsup_{n\to \infty}I_c^\mu[\pi_n] \leq \bar{I}_c^\mu[\pi]$, which completes the proof.
\endproof

\subsection*{Acknowledgments} The authors warmly thank the two anonymous referees for the careful reading and suggestions that greatly helped us to improve the quality of this work.

\bibliographystyle{amsplain}

\providecommand{\bysame}{\leavevmode\hbox to3em{\hrulefill}\thinspace}
\providecommand{\MR}{\relax\ifhmode\unskip\space\fi MR }
\providecommand{\MRhref}[2]{%
  \href{http://www.ams.org/mathscinet-getitem?mr=#1}{#2}
}
\providecommand{\href}[2]{#2}

\end{document}